\documentclass[10pt,a4paper]{amsart}
\numberwithin{equation}{section}

\usepackage{cite}                            
\usepackage{amssymb}
\usepackage{amsmath}
\usepackage{graphicx}
\usepackage{color}
\usepackage{subfigure}
\usepackage{xcolor}
\usepackage{rotating}
\usepackage{array,multirow}
\usepackage{adjustbox}

\makeatletter
\def\author@andify{%
	\nxandlist {\unskip ,\penalty-1 \space\ignorespaces}%
	{\unskip {} \@@and~}%
	{\unskip \penalty-2 \space \@@and~}%
}
\makeatother
\graphicspath{{pics/}}
\renewcommand{\rho}{u}
\renewcommand{\j}{v}

\newcommand{\dt}{\Delta t}

\newcommand{\e}{{\varepsilon}}

\newcommand{\be}{\begin{equation}}
\newcommand{\ee}{\end{equation}}
\newcommand{\bea}{\begin{eqnarray}}
\newcommand{\eea}{\end{eqnarray}}
\newcommand{\bean}{\begin{eqnarray*}}
\newcommand{\eean}{\end{eqnarray*}}

\def\ba{\begin{array}{l}\displaystyle}
\def\ea{\end{array}}
\def\epsi{{\varepsilon}}
\def\epsilon{{\varepsilon}}
\def\talpha{{\tilde\alpha}}
\def\tepsi{{\tilde\varepsilon}}

\newtheorem{theorem}{Theorem}
\definecolor{PBlue}{rgb}{0.0,0.0,0.69}
\definecolor{PGreen}{rgb}{0.0,0.69,0.0}
\definecolor{PMagenta}{rgb}{0.69,0.0,0.69}
\definecolor{PRed}{rgb}{0.8,0.0,0.0}

\newcommand{\revised}[1]{\textcolor{black}{#1}}	

\newcommand{\comment}[1]{\textcolor{black}{#1}}

\title[IMEX multistep methods for hyperbolic system with relaxation]{Implicit-Explicit multistep methods for hyperbolic systems with multiscale relaxation}\thanks{The work of G. Albi has been partially supported by the project Computer Science for Industry 4.0, 'MIUR Departments of Excellence 2018-2022'. The authors acknowledge the partial support of GNCS-INdAM 2019 project Numerical
approximation of hyperbolic problems and applications.}
\author{Giacomo Albi} \thanks{Computer Science Department, University of Verona, ({\tt giacomo.albi@univr.it}).}
\author{Giacomo Dimarco} \thanks{Mathematics Department, University of Ferrara, Italy ({\tt giacomo.dimarco@unife.it}).}
\author{Lorenzo Pareschi}\thanks{Mathematics Department, University of Ferrara, Italy ({\tt lorenzo.pareschi@unife.it}).}

\begin{document}

	\begin{abstract}
		We consider the development of high order space and time numerical methods based on Implicit-Explicit (IMEX) multistep time integrators for hyperbolic systems with relaxation. More specifically, we consider hyperbolic balance laws in which the convection and the source term may have very different time and space scales. As a consequence the nature of the asymptotic limit changes completely, passing from a hyperbolic to a parabolic system. From the computational point of view, standard numerical methods designed for the fluid-dynamic scaling of hyperbolic systems with relaxation present several drawbacks and typically lose efficiency in describing the parabolic limit regime. In this work, in the context of Implicit-Explicit linear multistep methods we construct high order space-time discretizations which are able to handle all the different scales and to capture the correct asymptotic behavior, independently from its nature, without time step restrictions imposed by the fast scales. Several numerical examples confirm the theoretical analysis.
	\end{abstract}
		\maketitle
	 \medskip
	

	\tableofcontents

\section{Introduction}
The goal of the present work is to develop high order numerical methods based on IMEX linear multistep (IMEX-LM) schemes for hyperbolic systems with relaxation \cite{CLL, Liu, Natalini}. These systems often contain multiple space-time scales which may differ by several orders of magnitude. In fact, the various parameters characterizing the models permit to describe different physical situations, like flows which pass from compressible to incompressible regimes or flows which range from rarefied to dense states. This is the case, for example, of kinetic equations close to the hydrodynamic limits \cite{BGL, Ce, CIP, Vill}. In such regimes these systems can be more conveniently described in terms of macroscopic equations since these reduced systems permit to describe all the features related to the space-time scale under consideration \cite{BGL, Vill}. However, such macroscopic models can not handle all the possible regimes one is frequently interested in. For such reason one has to resort to the full kinetic models. They permit to characterize a richer physics but on the other hand they are computationally more expensive and limited by the stiffness induced by the scaling under consideration \cite{DP15}.

The prototype system we will use in the rest of the paper is the following \cite{BPR17, NP2} 
\be
\left\{  
\begin{array}{l} 
\displaystyle  
\partial_t u + \partial_x v =0, \\
\displaystyle    
\partial_t v + \frac{1}{\epsi^{2\alpha}} \partial_x p(u) = -\frac{1}{\epsi^{1+\alpha}} \left( v-f(u ) \right),\quad \alpha\in[0,1]
\\
\end{array}
\right. 
\label{I61}
\ee
where $\varepsilon$ is the scaling factor and $\alpha$ characterizes the different type of asymptotic limit that can be obtained. The condition $p^\prime(u)>0$ should be satisfied for hyperbolicity to hold true since the eigenvalues of (\ref{I61}) are given by $\pm \sqrt{p^\prime(u)}/\epsi^{\alpha}$. Note that, except for the case $\alpha=0$, the eigenvalues are unbounded for small values of $\epsi$.

\revised{System \eqref{I61} is obtained from a classical $(2\times 2)$ $p$-system with relaxation under the space-time scaling $t\to t/\e^{1+\alpha}$, $x\to x/\e$ and by the change of variables $ v=\tilde v/\varepsilon^\alpha$, $ f(u)=\tilde f(u)/\varepsilon^\alpha$, where $\tilde v$ is the original unknown and $\tilde f(u)$ the original flux associated to the variable $u$ in the non rescaled $p$-system.} For $\alpha = 0$, system (\ref{I61}) reduces to the usual hyperbolic scaling 
\be
\left\{  
\begin{array}{l} 
	\displaystyle  
	\partial_t u + \partial_x v =0, \\
	\displaystyle    
	\partial_t v +\partial_x p(u) = -\frac{1}{\epsi} \left( v-f(u ) \right),
	\\
\end{array}
\right. 
\label{I61_hyp}
\ee
whereas for $\alpha=1$ yields the so-called diffusive scaling
\be
\left\{  
\begin{array}{l} 
	\displaystyle  
	\partial_t u + \partial_x v =0, \\
	\displaystyle    
	\partial_t v +\frac1{\epsi^2}\partial_x p(u) = -\frac{1}{\epsi^2} \left( v-f(u ) \right).
	\\
\end{array}
\right. 
\label{I61_par}
\ee
More in general, thanks to the Chapman-Enskog expansion \cite{Ce}, for small values of $\epsi$ we get from \eqref{I61} the following nonlinear convection-diffusion equation
\be
\left\{
\begin{array}{l}
\displaystyle
v= f(u) - \epsi^{1-\alpha} \partial_x p(u)+\epsi^{1+\alpha}f'(u)^2 \partial_x u + {\mathcal{O}}(\epsi^2), \\
\displaystyle
\partial_t u + \partial_x f(u )= \epsi^{1+\alpha} \partial_{x}\Bigg[\left(\frac{p'(u)}{\epsi^{2\alpha}}-f'(u)^2\right)\partial_x u\Bigg]+{\mathcal{O}}(\epsi^2). 
\\
\end{array}
\right.
\label{I7b}
\ee
In the limit $\varepsilon\rightarrow 0$, for $\alpha \in [0,1)$,  we are led to the conservation law 
\be
\left\{
\begin{array}{l}
\displaystyle
v= f(u ), \\
\displaystyle
\partial_t u + \partial_x f(u )= 0, 
\\
\end{array}
\right.
\label{I7}
\ee
while, when $\alpha=1$, in the asymptotic limit we obtain the following advection-diffusion equation
\be
\left\{  
\begin{array}{l} 
	\displaystyle  
	v= f(u ) -\partial_{x} p(u ),\\
	\displaystyle    
	\partial_t u + \partial_x f(u )= \partial_{xx} p(u ). 
	\\
\end{array}
\right. 
\label{I4}
\ee
Note that, the main stability condition for system (\ref{I7b}) corresponds to
\be
 f^\prime (u ) ^2 < \frac{p^\prime(u)}{\epsi^{2\alpha}}, 
\label{I5}
\ee
and it is always satisfied in the limit $\epsi\rightarrow 0$ when $\alpha > 0$, whereas for $\alpha=0$ the function $p(u)$ and $f(u)$ must satisfy the classical sub-characteristic condition \cite{CLL, Liu}. 

The space-time scaling just discussed, in classical kinetic theory, is related to the hydrodynamical limits of the Boltzmann equation. In particular, for $\alpha=0$ it corresponds to the compressible Euler scaling, whereas for $\alpha \in (0,1)$ to the incompressible Euler limit. In the case $\alpha=1$ dissipative effects become non-negligible and we get the incompressible Navier-Stokes scaling. \revised{We refer to \cite[Chapter~11]{CIP} and \cite{Vill} for further details and the mathematical theory behind the hydrodynamical limits of the Boltzmann equation. Moreover, we refer to \cite{LT, MR} for theoretical results on the diffusion limit of a system like \eqref{I61_par}}. 

The development of numerical methods to solve hyperbolic systems with stiff source terms has attracted many researches in the recent past \cite{CJR, GPT, GPT, GPR, Pe, PR, BR, Klar, JL, JPT, JPT2, DP2, Lem2}. The main computational challenge is related to the presence of the different scales that require a special care to avoid loss of stability and spurious numerical solutions. In particular, in diffusive regimes the schemes should be capable to deal with the very large characteristic speeds of the system \revised{by avoiding a CFL condition of the type $\Delta t = {\mathcal O}(\varepsilon^{\alpha})$}. A particular successful class of schemes is represented by the so-called \emph{asymptotic-preserving} (AP) schemes which aims at preserving the correct asymptotic behavior of the system without any loss of efficiency due to time step restrictions related to the small scales \cite{DP,Jin,JP,JF}. 





In the large majority of these works, the authors focused on the specific case $\alpha=0$, 
where a hyperbolic to hyperbolic scaling is studied, or to the  case $\alpha=1$, where a hyperbolic to parabolic scaling is analysed. 
Very few papers have addressed the challenging multiscale general problem for the various possible values of $\alpha \in [0,1]$ (see \cite{NP,JP,BPR17,GPR}) and all of them refer to one-step IMEX methods in a Runge-Kutta setting. We refer to \cite{Ak1, Ak2, Ascher2, DPLMM, FHV, HR, RSSZ1, RSSZ2} for various IMEX-LM methods developed in the literature and we mention that comparison between IMEX Runge-Kutta methods and IMEX-LM methods have been presented in \cite{HR}.
 
In the present work, following the approach recently introduced by Boscarino, Russo and Pareschi in \cite{BPR17}, we analyze the construction of IMEX-LM for such problems that work uniformly independently of the choices of $\alpha$ and $\epsi$. \revised{By this, we mean that the schemes are designed in such a way as to be stable for all different ranges of the scaling parameters independently of the time step. At the same time, they should ensure high order in space and time and should be able to accurately describe the various asymptotic limits. Moreover, whenever possible, the above described properties must be achieved without the need of an iterative solver for non linear equations.} 

In \cite{BPR17}, using a convenient partitioning of the original problem, the authors developed IMEX  Runge-Kutta schemes for a system like (\ref{I61}) which work uniformly with respect to the scaling parameters. Here, we extend these results to IMEX-LM, \revised{ previous results for IMEX-LM methods refer to the case $\alpha=0$ (see \cite{DPLMM, HR})}. Among others, there are two main reasons to consider the development of such schemes. First, in contrast to the IMEX Runge-Kutta case, for IMEX-LM it is relatively easy to construct schemes up to fifth order in time and typically they show a more uniform behavior of the error with respect to the scaling parameters. Second, thanks to the use of BDF methods, it is possible consider only one evaluation of the source term per time step independently of the scheme order. The latter feature is particularly significant in term of computational efficiency for kinetic equations where often the source term represents the most expensive part of the computation \cite{DP,DPLMM}. 

The rest of the paper is organized as follows. In Section 2, we discuss the discretization of these multiscale problems and motivate our partitioning choice of the system by analyzing a simple first order IMEX scheme. Next, in Section 3, we introduce the general IMEX-LM methods and discuss the asymptotic-preserving properties of our approach. \revised{Two classes of schemes are considered, AP-explicit and AP-implicit accordingly to the way the diffusion term in the limit equation is treated}. In Section 4, we perform a linear stability analysis for \revised{$2\times 2$ linear systems} in the case of IMEX-LM methods based on a backward-differentiation formula (BDF). \revised{ In Section 5, the space discretization is briefly discussed. Several numerical examples are reported in Section 6 which confirm the theoretical findings}. Final considerations and future developments are discussed at the end of the article.

\section{First order IMEX discretization}
\label{sec:1}
\begin{figure}\centering
	{\includegraphics[width=7cm]{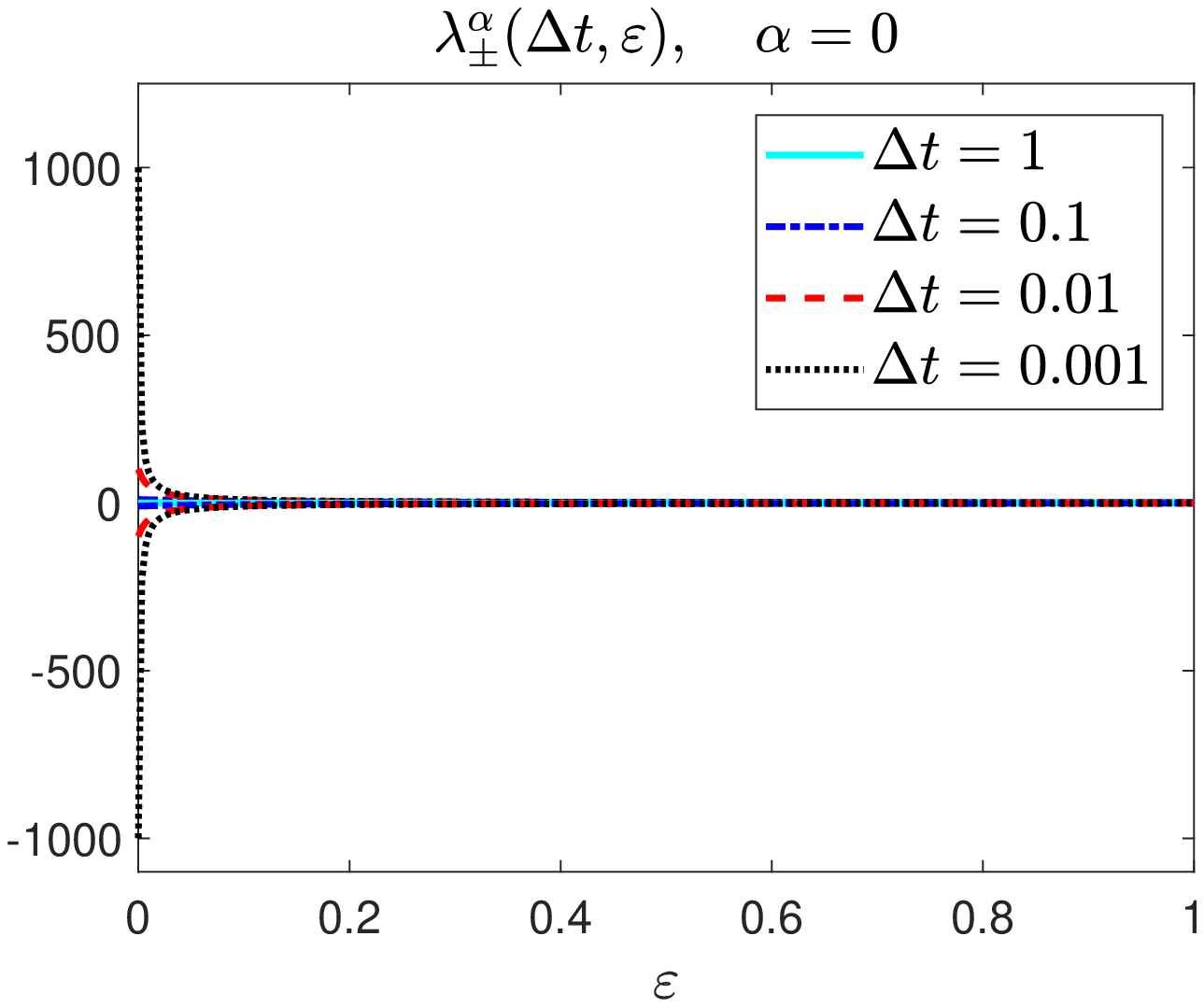}}\hskip .5cm
		{\includegraphics[width=7cm]{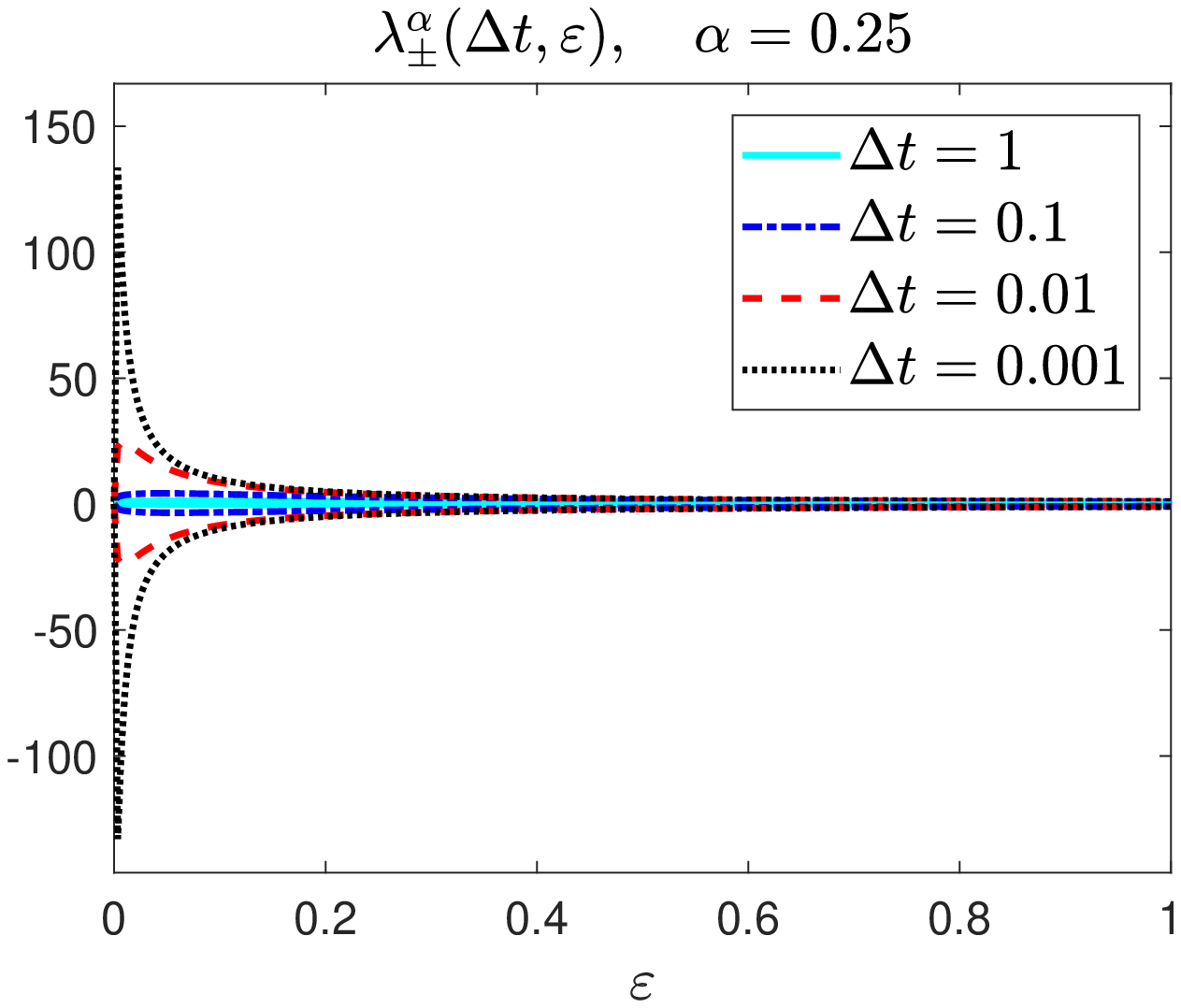}}
	\vskip .5cm
	{\includegraphics[width=7cm]{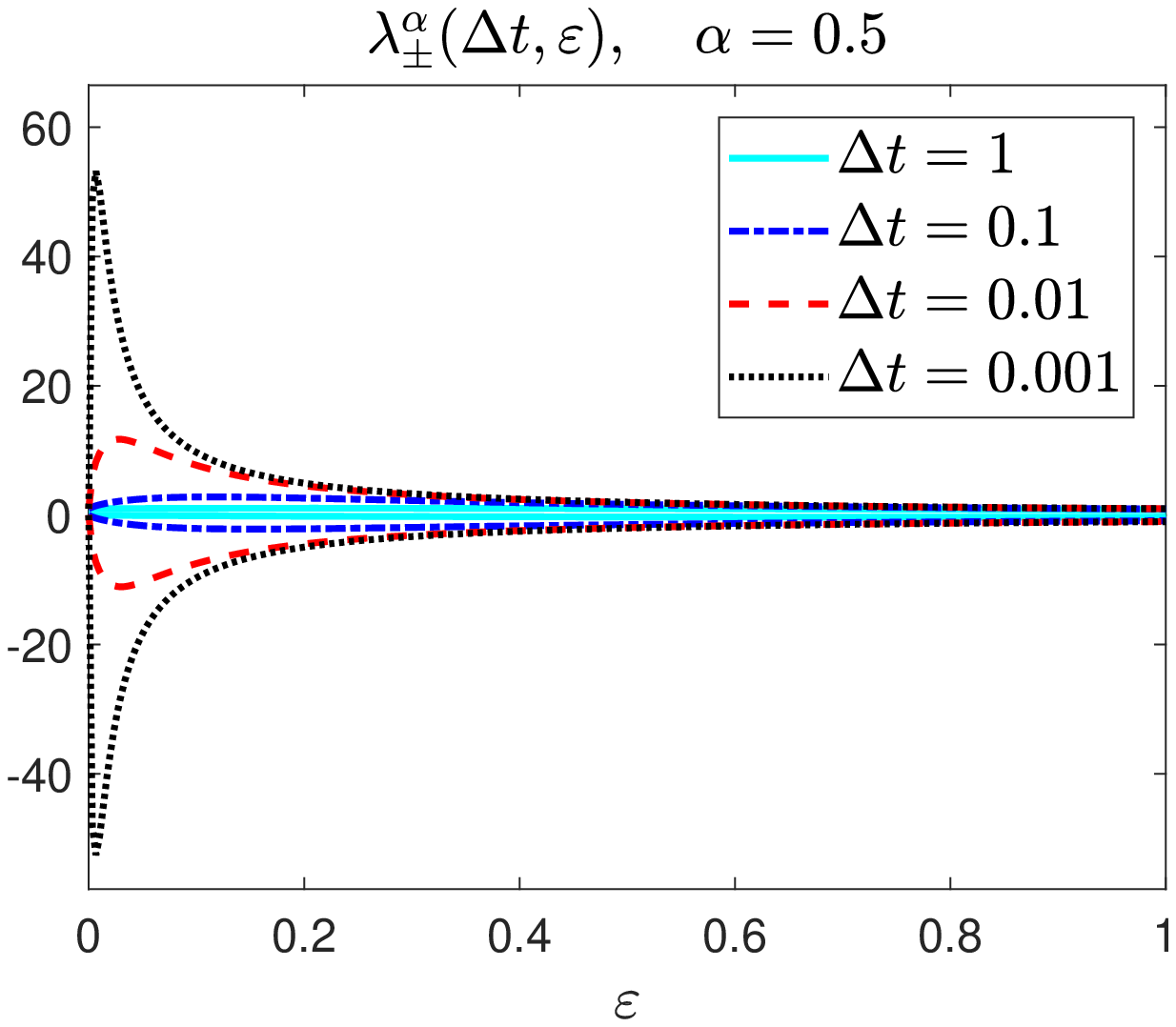}}\hskip .5cm
		{\includegraphics[width=7cm]{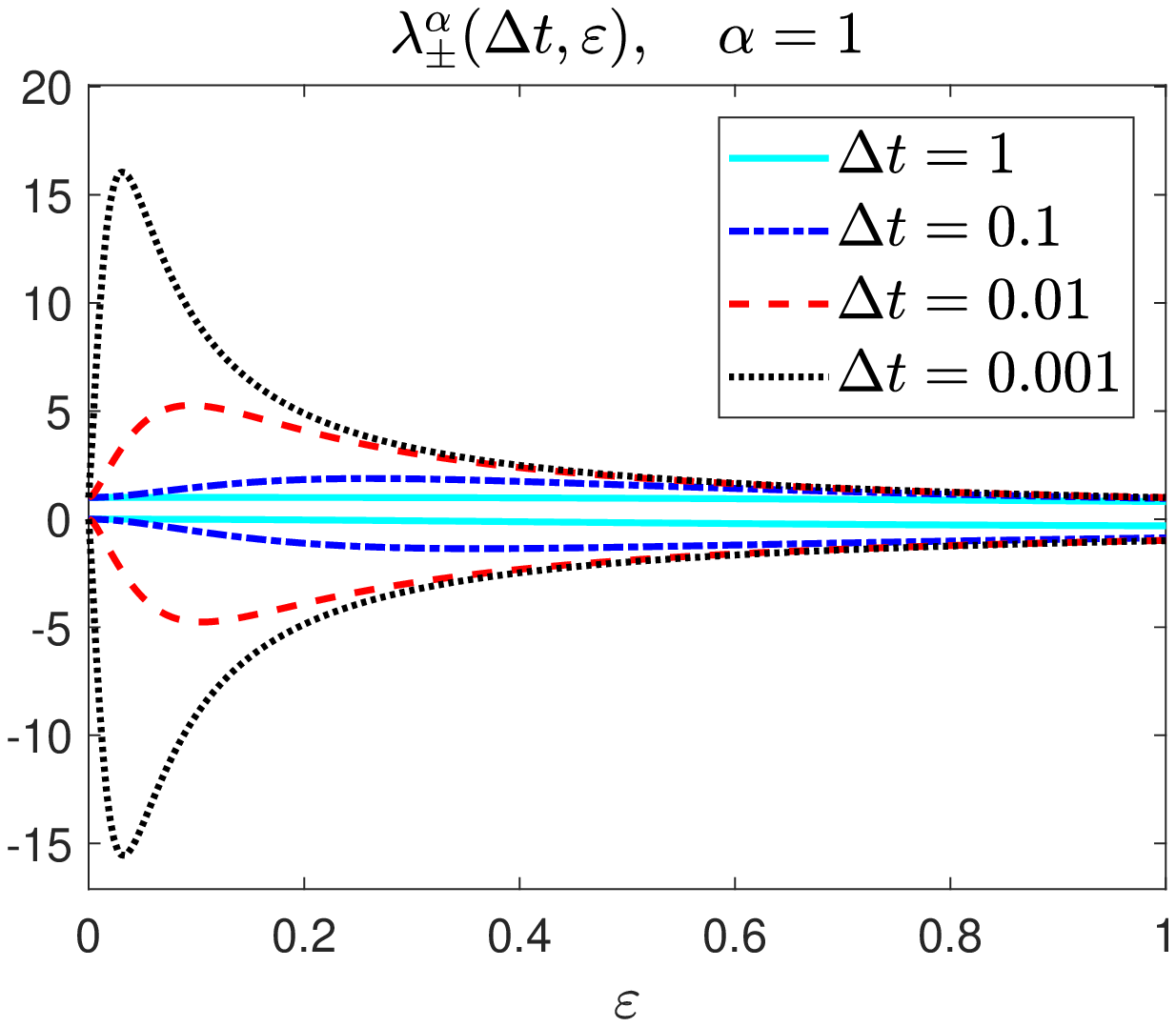}}
	\caption{Eigenvalues of the modified system \eqref{eq:SP2bisa} as a function of $\epsi$ for different values of the time step $\Delta t$ and choices of $\alpha$. }\label{fig:lambda_alpha}
\end{figure}

In this part, we discuss a first order IMEX time-discretization of the relaxation system (\ref{I61}) and we analyze its relationship with a reformulated system in which the eigenvalues are bounded for any value of the scaling parameter $\varepsilon$. To this aim, following \cite{BPR17}, we consider the following implicit-explicit first order partitioning of system \eqref{I61} 
\be
\label{eq:SP1}
\begin{aligned}
\frac{u^{n+1}-u^n}{\dt} & = - \partial_x v^{n+1}, \\
\varepsilon^{1+\alpha}\frac{v^{n+1}-v^n}{\dt} & = -\left(\varepsilon^{1-\alpha}\partial_x p(u^n) + v^{n+1} - f(u^n)\right). 
\end{aligned}
\ee
One can notice that in system (\ref{eq:SP1}) besides its implicit form, the second equation can be solved explicitly by inversion of the linear term $v^{n+1}$. This gives 
\be
   v^{n+1} = \frac{\epsi^{1+\alpha}}{\epsi^{1+\alpha}+\dt}v^n - \frac{\dt}{\epsi^{1+\alpha}+\dt}\left(\epsi^{1-\alpha}\partial_x p(u^n)-f(u^n)\right). 
\ee
Then, making use of the above relation and inserting it in the first equation, one gets 
\be
\label{eq:SP2u}
\begin{aligned}
\frac{u^{n+1}-u^n}{\dt} + \frac{\epsi^{1+\alpha}}{\epsi^{1+\alpha}+\dt}v^n_x  + \frac{\dt}{\epsi^{1+\alpha}+\dt}\partial_x f(u^n) & = \frac{\dt\,\epsi^{1-\alpha}}{\epsi^{1+\alpha}+\dt}\partial_{xx} p(u^n),
\end{aligned}
\ee
while a simple rewriting of the second equation gives
\be
\label{eq:SP2n}
\begin{aligned}
\frac{v^{n+1}-v^n}{\dt} + \frac{\epsi^{1-\alpha}}{\epsi^{1+\alpha}+\Delta t}\partial_x p(u^n) & = - \frac{1}{\epsi^{1+\alpha}+\Delta t}\left(v^n - f(u^n)\right).
\end{aligned}
\ee
Therefore, the IMEX scheme can be recast in an equivalent fully explicit form. \revised{Similarly to the continuous case,} depending on the choice of $\alpha$, as $\epsi\to 0$, we have different limit behaviors. For $\alpha\in[0,1)$ we obtain
\be
\label{eq:SP4b}
\begin{aligned}
	\frac{u^{n+1}-u^n}{\dt} +  \partial_x f(u^n) & = 0,
\end{aligned}
\ee
whereas in the case $\alpha=1$ we get
\be
\label{eq:SP4}
\begin{aligned}
	\frac{u^{n+1}-u^n}{\dt} +  \partial_x f(u^n) & = \partial_{xx} p(u^n).
\end{aligned}
\ee
\revised{
For small values of $\dt$, the scheme (\ref{eq:SP2u})-(\ref{eq:SP2n}) corresponds up to first order in time to the system
\be
\label{eq:SP2bisa}
\begin{aligned}
	\partial_t u + \frac{\epsi^{1+\alpha}}{\epsi^{1+\alpha}+\dt}\partial_x v  + \frac{\dt}{\epsi^{1+\alpha}+\dt}\partial_x f(u) & = \frac{\dt\,\epsi^{1-\alpha}}{\epsi^{1+\alpha}+\dt}\partial_{xx} p(u),\\
	\partial_t v + \frac{\epsi^{1-\alpha}}{\epsi^{1+\alpha}+\Delta t}\partial_x p(u) & = - \frac{1}{\epsi^{1+\alpha}+\Delta t}\left(v - f(u)\right),
\end{aligned}
\ee
where the following Taylor expansion has been employed at $t=t^n=n\Delta t$
\[
\displaystyle \frac{u^{n+1} - u^n}{\Delta t} = \partial_t u\big|_{t=t^n} + \mathcal{O}(\Delta t), \quad \frac{v^{n+1} - v^n}{\Delta t} = \partial_t v\big|_{t=t^n} + \mathcal{O}(\Delta t).
\]}
The main feature of system (\ref{eq:SP2bisa}) is that its left-hand side has bounded characteristic speeds. These are given by
\be\label{eq:eigen_alpha_1}
\lambda^{\alpha}_{\pm}(\dt,\epsi) = \frac{1}{2}\left(\gamma(1-\theta_{\alpha})\pm\sqrt{\gamma^2(1-\theta_\alpha)^2+4\epsi^{-2\alpha}\theta_\alpha^2}\right),
\ee
with 
$$\theta_{\alpha}(\dt,\epsi) := \frac{\epsi^{1+\alpha} }{\epsi^{1+\alpha} + \dt},$$ and where, for simplicity,
\revised{we considered $f'(u)=\gamma$, $\gamma \in \mathbb{R}$ and $p'(u)=1$ so that
\be
\begin{aligned}\label{eq:simplifications}
\partial_x f(u) = f'(u)\partial_x u = \gamma \partial_x u, \qquad
\partial_x p(u) = p'(u)\partial_x u = \partial_x u. 
\end{aligned}
\ee}
If we fix $\epsi$ and send $\Delta t \to 0$ we obtain the usual characteristic speeds of the original hyperbolic system, i.e. 
\[
\lambda^{\alpha}_{\pm}(0,\epsi)=\pm \frac1{\epsi^\alpha},
\]
while for a fixed $\dt$, the characteristic speeds $\lambda^\alpha_+$ and $\lambda^\alpha_-$ are respectively decreasing and increasing functions of $\varepsilon$ and, as $\epsi\to 0$, they converge to 
\be
\label{char_lim}
\lambda^{\alpha}_{\pm}(\Delta t,{0}) =\frac{1}{2} \left(\gamma\pm |\gamma|\right).
\ee
In Figure \ref{fig:lambda_alpha}, we show the shape of the eigenvalues \eqref{eq:eigen_alpha_1} for $\gamma=1$ and different values of the scaling parameter $\alpha$ and the time step $\Delta t$. We observe that when $\varepsilon$ grows \revised{the absolute value of} the eigenvalues diminish accordingly and when $\varepsilon$ diminishes the eigenvalues grow but remain bounded by the finite time step. 

Thus, for a given $\Delta t$, if we denote by $\Delta x$ the space discretization parameter, from the left hand side of \eqref{eq:SP2bisa} we expect the hyperbolic CFL condition $\Delta t \le \Delta x/|\gamma|$ in the limit $\epsi\to 0$. On the other hand, the stability restriction coming from the parabolic term requires $\Delta t= {\mathcal O}(\Delta x^2)$ when $\alpha=1$. 


In the next Section, we will show how to generalize the above arguments to the case of high order IMEX multistep methods.

{

\section{\revised{AP-explicit and AP-implicit IMEX-LM methods}}
\revised{In this part, we focus our attention on order $s$, $s$-step IMEX-LM methods with $s\geq 2$ (see Appendix A.1 for derivation and order conditions). First we discuss methods that yield a fully explicit discretization in the limit $\varepsilon\to 0$. For clarity of presentation, we separate the discussion of the diffusive case $\alpha=1$ from the general case $\alpha\in[0,1)$. In the second part, we discuss IMEX-LM discretizations which deal 
with the stiffness caused by the parabolic term in the asymptotic limit.}  

\subsection{\revised{AP-explicit methods in the diffusive case: $\alpha=1$}}

In this case, we can write the $s-$step IMEX-LM for the original hyperbolic system (\ref{I61}) as follows
\be
\label{eq:SPs}
\begin{aligned}
u^{n+1}& = -\sum_{j=0}^{s-1}a_ju^{n-j} -{\dt} \sum_{j=-1}^{s-1}c_j\partial_x v^{n-j}, \\
v^{n+1}& = -\sum_{j=0}^{s-1}a_jv^{n-j} -\frac{\dt}{\varepsilon^2}\left(\sum_{j=0}^{s-1}b_j \partial_x p(u^{n-j}) + \sum_{j=-1}^{s-1}c_j v^{n-j} - \sum_{j=0}^{s-1}b_j f(u^{n-j})\right),\\
\end{aligned}
\ee
where we introduced the following coefficients
\be
\label{eq:SPs_coeff}
\begin{aligned}
{Explicit}\qquad &a^T = (a_0,a_1,\ldots, a_{s-1})\\
\qquad &b^T = (b_0,b_1,\ldots, b_{s-1})\\
{Implicit}\qquad &c_{-1} \neq 0,\quad c^T = (c_0,c_1,\ldots, c_{s-1}).\\
\end{aligned}
\ee
\revised{Methods for which $c_j=0$, $j=0,\ldots,s-1$ are referred to as implicit-explicit backward differentiation formula, IMEX-BDF in short. Another important class of LM is represented by implicit-explicit Adams methods, for which $a_0=-1$, $a_j=0$, $j=1,\ldots,s-1$.} We refer to Appendix \ref{app:imex_multistep} for a brief survey of some IMEX multistep methods, and to~\cite{Ak1, Ak2, Ascher2, DPLMM, HR, FHV, RSSZ1, RSSZ2} for further details and additional schemes. 

In what follows, we rely on the equivalent vector-matrix notation,
\be
\label{eq:SPs_vec}
\begin{aligned}
u^{n+1} & =-a^T\cdot U - \dt c^T\cdot {\partial_x V} - \dt c_{-1}\partial_x v^{n+1}, \\
v^{n+1}& =-a^T\cdot V -\frac{\dt}{\varepsilon^2}\left(b^T \cdot \partial_x p(U) + c^T \cdot V  + c_{-1}v^{n+1} - b^T\cdot f(U)\right),
\end{aligned}
\ee
where $U = (u^{n},\ldots,u^{n-s+1})^T$, $V = (v^{n},\ldots,v^{n-s+1})^T$, $\partial_x p(U) = (\partial_x p(u^{n}),\ldots,\partial_x p(u^{n-s+1}))^T $ and $f(U) = (f(u^{n}),\ldots,f(u^{n-s+1}))^T $  are $s$-dimensional vectors.

Similarly to the one step scheme \eqref{eq:SP1}, we proceed by rewriting the multistep methods in the fully explicit vector form. To that aim, we observe that the second equation in \eqref{eq:SPs_vec} can be explicitly solved in terms of $v$ due to the linearity in the relaxation part. This gives
\be
\label{eq:SPs_vec_ex1}
\begin{aligned}
v^{n+1} &=-\frac{\epsi^2}{\epsi^2 + \dt c_{-1}}a^T\cdot V -\frac{\dt}{\epsi^2 + \dt c_{-1}}\left(b^T \cdot \partial_x p(U) + c^T \cdot V  - b^T\cdot f(U)\right)\\
& = -a^T\cdot V +\left(a^T - \frac{\epsi^2a^T}{\epsi^2 + \dt c_{-1}}-\frac{\dt c^T}{\epsi^2 + \dt c_{-1}}\right)\cdot V - \frac{\dt b^T}{\epsi^2 + \dt c_{-1}}\cdot\left(\partial_x p(U)-f(U)\right),
\end{aligned}
\ee
and thus
\be
\label{eq:SPs_vec_ex1b}
	v^{n+1}
	= -a^T\cdot V -\frac{\dt}{\epsi^2 + \dt c_{-1}}\left(c^T-c_{-1}a^T\right)\cdot V - \frac{\dt }{\epsi^2 + \dt c_{-1}}b^T\cdot\left(\partial_x p(U)-f(U)\right).
\ee
Substituting equation (\ref{eq:SPs_vec_ex1b}) into the first one of \eqref{eq:SPs_vec} leads to 
\be
\label{eq:SPs_vec_ex2}
\begin{aligned}
u^{n+1}&=-a^T\cdot U -\dt \left(c^T-c_{-1}a^T\right)\cdot {\partial_x V} +\\ &\,\, +\frac{\dt^2 c_{-1}}{\epsi^2 + \dt c_{-1}}\left(c^T-c_{-1}a^T\right)\cdot {\partial_x V} +\frac{\dt^2 c_{-1} }{\epsi^2 + \dt c_{-1}}b^T\cdot\partial_x \left(\partial_x p(U)-f(U)\right)\\
& = -a^T\cdot U -\dt \frac{\epsi^2 }{\epsi^2 \revised{+}\dt c_{-1}}\left(c^T-c_{-1}a^T\right)\cdot {\partial_x V} +\frac{\dt^2 c_{-1} }{\epsi^2 + \dt c_{-1}}b^T\cdot\partial_x \left(\partial_x p(U)-f(U)\right).
\end{aligned}
\ee
Hence, we obtain the following system 
\be
\label{eq:SPs_vec_ex}
\begin{aligned}
\frac{u^{n+1} +a^T\cdot U}{\Delta t} &= - \frac{\epsi^2 \left(c^T-c_{-1}a^T\right)}{\epsi^2 + \dt c_{-1}}\cdot {\partial_x V} +\frac{\dt c_{-1} }{\epsi^2 + \dt c_{-1}}b^T\cdot\partial_x \left(\partial_x p(U)-f(U)\right),\\
\frac{v^{n+1}+a^T\cdot V}{\Delta t} &=  - \frac{\left(c^T-c_{-1}a^T\right)}{\epsi^2 + \dt c_{-1}}\cdot V - \frac{1}{\epsi^2 + \dt c_{-1}}b^T\cdot\left(\partial_x p(U)-f(U)\right),
\end{aligned}
\ee
which is the generalization of system \eqref{eq:SP2n} to an $s$-step IMEX scheme. 

\revised{
Our aim now is to show that, similarly to the simple first order method analyzed in Section \ref{sec:1}, the above discretization for a small value of $\Delta t$ corresponds, up to first order in time, to a modified hyperbolic problem where the characteristic speeds are bounded even in the limit $\varepsilon\to 0$.
More precisely, see also \cite{ARS}, for a smooth solution in time, by Taylor series expansion about $t=t^{n}$, we have
\begin{eqnarray*}
U &=& e\,u\big|_{t=t^{n}} - \Delta t J \partial_t u \big|_{t=t^{n}} + \ldots + (-1)^s\frac{\Delta t^{s}}{s!} J^s \partial^s_t u\big|_{t=t^{n}}+ {\mathcal O}(\Delta t^{s+1})\\
V &=& e\,v\big|_{t=t^{n}} - \Delta t J \partial_t v \big|_{t=t^{n}} + \ldots + (-1)^s\frac{\Delta t^{s}}{s!} J^s \partial^s_t v\big|_{t=t^{n}}+ {\mathcal O}(\Delta t^{s+1}),
\end{eqnarray*}
where $e$ is a vector of ones in $\mathbb{R}^{s}$, $J=(0,\ldots,s-1)^T$, $\partial^q_t$, $q=1,\ldots,s$ denotes the $q$-derivative and the vector powers must be understood component-wise. Similarly, we can expand $u^{n+1}$ and $v^{n+1}$ about $t=t^n$.
} 
\revised{
Therefore, we obtain
\begin{eqnarray}
\nonumber
\frac{u^{n+1} +a^T\cdot U}{\Delta t} &=& \frac{(1+a^T\cdot e)}{\Delta t} u\big|_{t=t^{n}}+ (1-a^T\cdot J)\partial_t u \big|_{t=t^{n}}+{\mathcal O}(\Delta t),\\[-.2cm]
\label{eq:IMEXest2}
\\[-.2cm]
\nonumber
\frac{v^{n+1} +a^T\cdot V}{\Delta t} &=& \frac{(1+a^T\cdot e)}{\Delta t} v\big|_{t=t^{n}}+ (1-a^T\cdot J)\partial_t v \big|_{t=t^{n}}+{\mathcal O}(\Delta t).
\end{eqnarray}}
\revised{From the order conditions we have 
\[
1+a^T\cdot e =0,\qquad 1-a^T\cdot J=  b^T\cdot e = c_{-1}+c^T\cdot e =: \beta, 
\]
and subsequentely
\[
\left(c^T-c_{-1}a^T\right)\cdot e  =c^T\cdot e + c_{-1} = \beta.
\]
Thus, scheme \eqref{eq:SPs_vec_ex} for small values of $\dt$ can be considered as a first order approximation  of the following modified system
\be\label{modsys2}
\begin{aligned}
\partial_t u &+ \frac{\epsi^2 }{\epsi^2 + \dt c_{-1}} \partial_x v  +\frac{\dt c_{-1} }{\epsi^2 + \dt c_{-1}} \partial_x f(u) = \frac{\dt c_{-1} }{\epsi^2 + \dt c_{-1}}\ \partial_{xx} p(u)
\\
\partial_t v &+ \frac{1}{\epsi^2 + \dt c_{-1}} \partial_x p(u) =  - \frac{1}{\epsi^2 + \dt c_{-1}}\left(v - f(u)\right),
\end{aligned}
\ee
where the factor $\beta$ simplifies in all the terms.}

\revised{
Note that, the above system has exactly the same structure as system \eqref{eq:SP2bisa} (with $c_{-1}=1$ and $\alpha=1$) which was derived from to the first order time discretization. As a consequence, under the same simplification assumptions \eqref{eq:simplifications} on the fluxes $f(u)$ and $p(u)$, the eigenvalues of the hyperbolic part correspond to} 
\be
\Lambda_{\pm}(\dt,\epsi) = \frac{1}{2}\left(\gamma(1-\theta_1)\pm\sqrt{\gamma^2(1-\theta_1)^2+4\epsi^{-2}\theta_1^2}\right),
\label{eq:eigenvalues}
\ee
where
\be
\theta_1(\dt,\epsi) := \frac{\epsi^2 }{\epsi^2 + \dt c_{-1}}.
\ee
Thus, the bounds for the characteristic velocities are the same as for the first order scheme and we get the limit cases
\[
\Lambda_{\pm}(\dt,0) = \frac{1}{2}\left(\gamma\pm|\gamma|\right),\qquad 
\Lambda_{\pm}(0,\epsi) =\pm\frac{1}{\epsi}.
\]

\subsubsection{Asymptotic preserving property for $\alpha=1$.}
Now, we study the capability of the schemes (\ref{eq:SPs_vec}) to become a consistent discretization of the limit system (\ref{I4}). To that aim, letting $\varepsilon\rightarrow 0$, in the reformulated scheme (\ref{eq:SPs_vec_ex}), we get from the first equation
\be
\label{eq:SPs_vec_ex_eps0u}
\begin{aligned}
	\frac{u^{n+1} +a^T\cdot U}{\Delta t} &= b^T\cdot\partial_x \left(\partial_x p(U)-f(U)\right),\\
\end{aligned}
\ee
which \revised{corresponds} to the explicit multistep scheme applied to the limiting convection-diffusion equation (\ref{I4}). \revised{For this reason, from now on, we refer to this class of IMEX-LM schemes as {\em AP-explicit} methods.}
On the other hand, we have for the second equation
\[
\begin{aligned}
	c_{-1}v^{n+1}+c^T\cdot V &= -b^T\cdot\left(\partial_x p(U)-f(U)\right),
\end{aligned}
\]
\revised{or equivalently
\be
\label{eq:SPs_vec_ex_eps0v}
v^{n+1} = -\frac{c^T}{c_{-1}}\cdot V-\frac{b^T}{c_{-1}}\cdot\left(\partial_x p(U)-f(U)\right).
\ee
}
Let us observe that in order to have at time $t^{n+1}$ a consistent projection over the asymptotic limit a condition over the states $U=(u^n,..,u^{n-s+1})^T$ and $V=(v^n,..,v^{n-s+1})^T$ should be imposed. This can be resumed by saying that \textit{the vector of the initial data should be well prepared to the asymptotic state}. This means that for the first variable
\be\label{well_p1}
u^{n-j}=\bar u^{n-j}+\tilde u^{n-j}_\varepsilon, \qquad \lim\limits_{\varepsilon\rightarrow 0}\tilde u^{n-j}_\varepsilon=0, \qquad j=0,...,s-1,
\ee
where $\bar u^{n-j}$ is a consistent solution of the limit system \eqref{I4} while $\tilde u^{n-j}_\varepsilon$ is a perturbation that disappears in the limit. An analogous relation should hold true for the second variable $v$, i.e.
\be\label{well_p2}
v^{n-j}=\bar v^{n-j}+\tilde v^{n-j}_\varepsilon, \qquad \lim\limits_{\varepsilon\rightarrow 0}\tilde v^{n-j}_\varepsilon=0, \qquad j=0,...,s-1,
\ee
\revised{
where $\bar v^{n-j}=f(u^{n-j})-\partial_x p(u^{n-j})$, $j=0,\ldots,s-1$ is a consistent projection of the asymptotic limit,  while $\tilde v^{n-j}_\varepsilon$ is a perturbation that disappears in the limit. Under this assumption, as a consequence of the order conditions, relation \eqref{eq:SPs_vec_ex_eps0v} is an $s$-order approximation of the asymptotic limit $v=f(u)-\partial_x p(u)$ and therefore, at subsequent time steps, the numerical solution is guaranteed to satisfy \eqref{well_p1}-\eqref{well_p2}. If such conditions are not imposed on the initial values, then the numerical solution may present a spurious initial layer and deterioration of accuracy is observed.  
In particular,
for IMEX-BDF methods, expression \eqref{eq:SPs_vec_ex_eps0v} simplifies to
\be
\label{eq:SPs_vec_ex_eps0v2}
v^{n+1} = -\frac{b^T}{c_{-1}}\cdot\left(\partial_x p(U)-f(U)\right).
\ee
This shows that, even for non well prepared initial data in $v$ but only in $u$ we obtain an $s$-order approximation of the equilibrium state and subsequently of the numerical solution for all times. This stronger AP property of IMEX-BDF methods is satisfied also in the asymptotic limits analyzed for the various schemes in the sequel of the manuscript. 
We refer to \cite{DPLMM,DP,PR} for more detailed discussions.
}

\subsection{\revised{AP-explicit methods in the general case: $\alpha\in [0,1)$.}}
The IMEX-LM scheme in vector form reads
\be
\label{eq:SPs_vec_a}
\begin{aligned}
u^{n+1} & =-a^T\cdot U - \dt c^T\cdot {\partial_x V} - \dt c_{-1}\partial_x v^{n+1}, \\
v^{n+1}&  =-a^T\cdot V -\frac{\dt}{\varepsilon^{2\alpha}}b^T \cdot \partial_x p(U) - \frac{\dt}{\varepsilon^{1+\alpha}}\left(c^T \cdot V  + c_{-1}v^{n+1} - b^T\cdot f(U)\right).
\end{aligned}
\ee
Again, we rewrite the second equation by solving it in terms of $v^{n+1}$ as follows
\begin{align*}
v^{n+1} &=-a^T\cdot V-\frac{\dt}{\epsi^{1+\alpha} +\dt c_{-1}}\left(c^T-c_{-1}a^T\right)\cdot V+\\&\,\, +\frac{\dt }{\epsi^{1+\alpha} + \dt c_{-1}}b^T\cdot f(U) - \frac{\dt \epsi^{1-\alpha} }{\epsi^{1+\alpha} + \dt c_{-1}}b^T\cdot \partial_x p(U),
\end{align*}
and using this solution in the first equation, we obtain the explicit scheme 
\revised{
\be
\label{eq:SPs_vec_ex2_alpha}
\begin{aligned}
	\frac{u^{n+1} +a^T\cdot U}{\Delta t} &= \frac{-\epsi^{1+\alpha}\left(c^T-c_{-1}a^T\right)}{\epsi^{1+\alpha} + \dt c_{-1}}\cdot {\partial_x V}-\frac{\dt c_{-1}b^T}{\epsi^{1+\alpha} + \dt c_{-1}}\cdot (\partial_x f(U)-\epsi^{1-\alpha}\partial_{xx} p(U)),\\
	\frac{v^{n+1}+a^T\cdot V}{\Delta t} &=  \frac{-\left(c^T-c_{-1}a^T\right)}{\epsi^{1+\alpha} + \dt c_{-1}}\cdot V+\frac{b^T}{\epsi^{1+\alpha} + \dt c_{-1}}\cdot \left(f(U)-\epsi^{1-\alpha} \partial_x p(U)\right).
\end{aligned}
\ee}
\revised{
 Considering, the same Taylor approximations \eqref{eq:IMEXest2} the above schemes corresponds up to first order in time to the modified system
 \be
\label{eq:SPs_vec_ex3}
\begin{aligned}
	\partial_t u + \frac{\epsi^{1+\alpha}}{\epsi^{1+\alpha}+\dt c_{-1}}\partial_x v  + \frac{\dt c_{-1}}{\epsi^{1+\alpha}+\dt c_{-1}}\partial_x f(u) & = \frac{\dt\,c_{-1} \epsi^{1-\alpha}}{\epsi^{1+\alpha}+\dt c_{-1}}\partial_{xx} p(u),\\
	\partial_t v + \frac{\epsi^{1-\alpha}}{\epsi^{1+\alpha}+\Delta t c_{-1}}\partial_x p(u) & = - \frac{1}{\epsi^{1+\alpha}+\Delta t c_{-1}}\left(v - f(u)\right).
\end{aligned}
\ee
Clearly, system \eqref{eq:SPs_vec_ex3} has again the same structure as \eqref{eq:SP2bisa}. As a consequence, under the same simplification assumptions \eqref{eq:simplifications}, the eigenvalues of the hyperbolic part correspond to 
\be\label{eq:eigen_alpha}
\Lambda^{\alpha}_{\pm}(\dt,\epsi) = \frac{1}{2}\left(\gamma(1-\theta_{\alpha})\pm\sqrt{\gamma^2(1-\theta_\alpha)^2+4\epsi^{-2\alpha}\theta_\alpha^2}\right),
\ee
where $\theta_{\alpha}$ is defined as follows
\be
\theta_{\alpha}(\dt,\epsi) := \frac{\epsi^{1+\alpha} }{\epsi^{1+\alpha} + \dt c_{-1}},
\ee
and  the bounds for the characteristic velocities for $\epsi = 0$ and $\dt = 0$ are 
\[
\Lambda^{\alpha}_{\pm}(\dt,0) = \frac{1}{2}\left(\gamma\pm|\gamma|\right),\qquad \Lambda^{\alpha}_{\pm}(0,\epsi) =\pm\frac{1}{\epsi^{\alpha}}.
\]}

\subsubsection{Asymptotic preserving property for $\alpha\in[0,1)$.}
We consider now the analogous asymptotic preserving property proved for the schemes (\ref{eq:SPs_vec}) in the case of the schemes (\ref{eq:SPs_vec_a}). Namely, we want to show that (\ref{eq:SPs_vec_a}) becomes a consistent discretization of the limit system (\ref{I7}) when $\varepsilon\rightarrow 0$. Taking scheme (\ref{eq:SPs_vec_ex2_alpha}), we get from the first equation
\be
\label{eq:SPs_vec_ex_eps0u_a}
\begin{aligned}
	u^{n+1}=- a^T\cdot U- \Delta t b^T\cdot \partial_x f(U),\\
\end{aligned}
\ee
which is a standard explicit multistep discretization of the asymptotic hyperbolic limit, i.e. of equation (\ref{I7}). On the other hand, we have for the second equation
\[
	c_{-1}v^{n+1}+c^T\cdot V = b^T\cdot f(U),
\]
\revised{
or equivalentely 
\be
\label{eq:SPs_vec_ex_eps0v_a}
v^{n+1} = -\frac{c^T}{c_{-1}}\cdot V + \frac{b^T}{c_{-1}}\cdot f(U).
\ee
As a consequence of the order conditions, equation \eqref{eq:SPs_vec_ex_eps0v_a} defines an $s$-order consistent approximation of the asymptotic limit $v=f(u)$ provided that the vector of the initial data is well prepared. These conditions are the analogous of (\ref{well_p1}) and (\ref{well_p2}) except that now
$\bar v^{n-j}=f(u^{n-j})$, $j=0,\ldots,s-1$.}

\revised{For IMEX-BDF methods \eqref{eq:SPs_vec_ex_eps0v_a} reduces to
\be
\label{eq:SPs_vec_ex_eps0v_a2}
v^{n+1} = \frac{b^T}{c_{-1}}\cdot f(U),
\ee
so that, even for non well prepared initial data in $v$ but only in $u$ we obtain an $s$-order approximation of the equilibrium state and subsequently of the numerical solution for all times.} 


\subsection{\revised{Removing the parabolic stiffness: AP-implicit methods}}
Although, the schemes developed in the previous Section overcome the stiffness related to the scaling factor $\varepsilon$, there is another stiffness that may appear in the equations close to the asymptotic limit. In fact, as shown previously, all the schemes originate a fully explicit scheme in the limit. 

\revised{In diffusive regimes, this typically leads to the time step restriction $\Delta t = {\mathcal O}(\Delta x^2)$ when 
\be
\label{eq:diffreg}
\epsi^{1-\alpha}\Delta t c_{-1} /(\epsi^{1+\alpha}+\Delta t c_{-1})={\mathcal O}(1),
\ee
(see the diffusion coefficient in equation \eqref{eq:SPs_vec_ex3}).
Therefore, for small $\varepsilon$ and in the case of $\alpha\simeq 1$, the main stability restriction is due to the second order term of the Chapmann--Enskog expansion (\ref{I7b}). Note that, beside the case $\alpha=1$ and $\epsi\to 0$ where we obtain a parabolic problem in the limit, the above time step limitation may occur also in transient regimes for $\alpha\neq 1$ as soon as \eqref{eq:diffreg} holds true.}

For this reason, we modify the partitioning of the system taking also $\partial_x p(u)$ implicit in the second equation as follows
\begin{eqnarray}
\nonumber
u^{n+1} & =&-a^T\cdot U - \dt c^T\cdot {\partial_x V} - \dt c_{-1}\partial_x v^{n+1}, \\[-.1cm]
\label{eq:SPs_vec1}
\\[-.2cm]
\nonumber
v^{n+1}&  =&-a^T\cdot V  - \frac{\dt}{\varepsilon^{1+\alpha}}\left(c^T \cdot V  + c_{-1}v^{n+1} - b^T\cdot f(U)\right)-\frac{\dt}{\varepsilon^{2\alpha}}\left(c^T \cdot \partial_x p(U) + c_{-1} \partial_x p(u^{n+1})\right).
\end{eqnarray}
We can still solve the second equation in $v$, to get
  \begin{align*}
v^{n+1}  &=-a^T\cdot V-\frac{\dt}{\epsi^{1+\alpha} + \dt c_{-1}}\left(c^T-c_{-1}a^T\right)\cdot V+\\& \,\,+\frac{\dt }{\epsi^{1+\alpha} + \dt c_{-1}}b^T\cdot f(U)- \frac{\dt \epsi^{1-\alpha} }{\epsi^{1+\alpha} + \dt c_{-1}}\left(c^T \cdot \partial_x p(U) + c_{-1} \partial_x p(u^{n+1})\right),
\end{align*}
which, inserted into the first equation of (\ref{eq:SPs_vec1}) yields the IMEX formulation
\be
\label{eq:SPs_vec_ex2_par}
\begin{aligned}
\frac{u^{n+1} +a^T\cdot U}{\Delta t} &= - \frac{\epsi^{1+\alpha}\left(c^T-c_{-1}a^T\right)}{\epsi^{1+\alpha} + \dt c_{-1}}\cdot {\partial_x V} -\frac{\dt c_{-1} }{\epsi^{1+\alpha} + \dt c_{-1}}b^T\cdot \partial_x f(U) \\
&\qquad\quad\quad+\frac{\epsi^{1-\alpha}\dt c_{-1} }{\epsi^{1+\alpha} + \dt c_{-1}}\left(c^T\cdot \partial_{xx} p(U)+ \revised{c_{-1}}\partial_{xx} p(u^{n+1})\right),\\
\frac{v^{n+1}+a^T\cdot V}{\Delta t} &=\quad  - \frac{\left(c^T-c_{-1}a^T\right)}{\epsi^{1+\alpha} + \dt c_{-1}}\cdot V +\frac{1}{\epsi^{1+\alpha} + \dt c_{-1}}b^T\cdot f(U)\\
&\qquad\quad-\frac{\epsi^{1-\alpha}}{\epsi^{1+\alpha} + \dt c_{-1}}\left(c^T\cdot \partial_x p(U)+\revised{c_{-1}} \partial_x p(u^{n+1})\right).
\end{aligned}
\ee
\revised{Note that, except for the case in which $p(u)$ is linear, in general, the first equation in \eqref{eq:SPs_vec_ex2_par} requires the adoption of a suitable solver for nonlinear problems to compute $u^{n+1}$.
By the same arguments of in the previous Sections, for small values of $\Delta t$, the scheme \eqref{eq:SPs_vec_ex2_par} corresponds up to first order to the modified system \eqref{eq:SPs_vec_ex3}. 
Thus, under the same simplification assumptions \eqref{eq:simplifications}, the eigenvalues of the hyperbolic part are given by (\ref{eq:eigen_alpha}).} 

\subsubsection{Asymptotic preserving property for the AP-implicit methods.}
Finally, we conclude our analysis by studying the asymptotic preserving property. \revised{As we will see the main difference is that in the asymptotic limit the diffusive terms are integrated implicitly. For this reason we refer to this class of IMEX-LM schemes as {\em AP-implicit} methods}. We first consider the case $\alpha=1$.
Taking the reformulated IMEX scheme (\ref{eq:SPs_vec_ex2_par}) and letting $\varepsilon\rightarrow 0$ with $\alpha=1$, gives
\be
\label{eq:SPs_vec_ex2ll}
\begin{aligned}
	\frac{u^{n+1} +a^T\cdot U}{\Delta t} &= -b^T\cdot \partial_x f(U) + c^T\cdot \partial_{xx} p(U)+ c_{-1}\partial_{xx} p(u^{n+1}),
\end{aligned}
\ee
which correspond to the IMEX multistep scheme applied to the limiting convection diffusion problem where the diffusion term is treated implicitly \cite{Ak1,Ak2}. \revised{Note that, for the second equation we have
\[
	c_{-1}v^{n+1}+c^T\cdot V = b^T\cdot f(U)-c^T\cdot \partial_{x} p(U)-c_{-1} \partial_{x} p(u^{n+1}),
\]
or equivalently
\be
\label{eq:SPs_vec_ex_eps0v_aa}
	v^{n+1}=-\frac{c^T}{c_{-1}}\cdot V +\frac{b^T}{c_{-1}}\cdot f(U)-\frac{c^T}{c_{-1}}\cdot \partial_{x} p(U)- \partial_{x} p(u^{n+1}),
\ee
which, under the assumption of well prepared initial values, as a consequence of the order conditions,  corresponds to a $s$-order approximation of the equilibrium projection $v=f(u)-\partial_x p(u)$. 
}

\revised{On the other hand, in the case $\alpha\in[0,1)$, we get the same asymptotic limit \eqref{eq:SPs_vec_ex_eps0u_a} of the AP-explicit method (see Section 3.2.1) 
\be
\label{eq:SPs_vec_ex_eps0u_apar}
\begin{aligned}
	u^{n+1}=- a^T\cdot U- \Delta t b^T\cdot \partial_x f(U).
\end{aligned}
\ee
We will not discuss further this limit system, but we emphasize that \eqref{eq:SPs_vec_ex_eps0u_apar} is obtained as the limit of the implicit-explicit scheme \eqref{eq:SPs_vec_ex2_par} whereas \eqref{eq:SPs_vec_ex_eps0u_a} is obtained as the limit of the explicit scheme \eqref{eq:SPs_vec_ex2_alpha}.
}

\section{Linear stability analysis}
Monotonicity properties for IMEX-LM have been previously studied in \cite{Ascher2,DPLMM,FHV,HR,HRS}. Due to the well-known difficulties in extending the usual stability analysis for linear systems to the implicit-explicit setting most results are limited to the single scalar equation. \revised{In our case, however, the schemes are specifically designed to deal with systems in the form \eqref{I61}, and we need therefore to tackle the stability properties in such a case.
Here we show that, in the case of IMEX-BDF methods, we can generalize some of the stability results for the single scalar equations to linear multiscale systems of the form}
\begin{eqnarray}
\nonumber
\partial_t u+\partial_x v&=&0,\\[-.3cm]
\label{eq:tp0}
\\[-.3cm]
\nonumber
\partial_t v+  \frac1{\varepsilon^{2\alpha}} \partial_x u&=&-\frac1{\varepsilon^{1+\alpha}} (v-\gamma u),
\end{eqnarray}
where $\varepsilon>0$, $\gamma >0$ and $\alpha\in [0,1]$. Note that, for small values of $\varepsilon$, the above system   when $\alpha=1$ reduces to the convection-diffusion equation $\partial_t u + \gamma \partial_x u = \partial_{xx} u $, whereas when $\alpha=0$, if $\gamma <1$, yields the simple advection equation $\partial_t u+\gamma \partial_x u=0$. To simplify notations, in the sequel we will assume $\gamma=1$ and $\alpha>0$. The case $\alpha=0$ is rather classical and follows straightforwardly from our analysis. Under these assumptions, the Chapman-Enskog expansion for small values of $\varepsilon$ gives the limiting convection-diffusion equation
\be
\partial_t u + \partial_x u = \epsi^{1-\alpha} \partial_{xx} u + {\mathcal O}(\epsi^{1+\alpha}).
\label{eq:cel}
\ee

In Fourier variables we get 
\begin{eqnarray}
\nonumber
{\hat u}'&=&-i\xi{\hat v},\\[-.3cm]
\label{eq:tp1f}
\\[-.3cm]
\nonumber
{\hat v}'&=&-\frac{i\xi}{\varepsilon^{2\alpha}}{\hat u}-\frac1{\varepsilon^{1+\alpha}}({\hat v}-{\hat u}),
\end{eqnarray}
where, for example, $\xi = \sin(2k)/\Delta x$ if we use central differences and $k$ is the frequency of the corresponding Fourier mode.

The change of variables $y={\hat u}$, $z=\varepsilon^{\alpha} {\hat v}$, $\lambda_I=i\xi/\varepsilon^{\alpha}$, $\lambda_R=1/\varepsilon^{1+\alpha}$, $\lambda=\lambda_I+\lambda_R\in\mathbb{C}$ transforms the system into the problem
\begin{eqnarray}
\nonumber
y'&=&-\lambda_I z,
\\[-.3cm]
\label{eq:tp1fd}&&\qquad\qquad\qquad\qquad\qquad\qquad \lambda\in \mathbb{C},
\\[-.3cm]
\nonumber
z'&=&-(\lambda_I -\lambda_R\varepsilon^{\alpha}) y-\lambda_R z.
\end{eqnarray}
Let us note that the above problem is equivalent to the second order differential equation
\be
y'' = - \lambda_R y'+\lambda_I(\lambda_I-\lambda_R \varepsilon^{\alpha})y .
\label{eq:so}
\ee
\subsection{AP-explicit methods}
We then apply an IMEX-LM method to system \eqref{eq:tp1fd} as follows
\begin{eqnarray}
\label{eq:first}
y^{n+1} &=& -\sum_{j=0}^{s-1} a_j y^{n-j} -\lambda_I \Delta t \sum_{j=-1}^{s-1} c_j z^{n-j}\\
\label{eq:second}
z^{n+1} &=& -\sum_{j=0}^{s-1} a_j z^{n-j} -(\lambda_I -\lambda_R\varepsilon^{\alpha}) \Delta t \sum_{j=0}^{s-1} b_j y^{n-j} -\lambda_R \Delta t \sum_{j=-1}^{s-1} c_j z^{n-j}.
\end{eqnarray}
In the case of IMEX-BDF methods the first equation \eqref{eq:first} permits to write 
\[
z^{n+1} = -\frac1{\Delta t c_{-1}\lambda_I}\left(y^{n+1} +\sum_{j=0}^{s-1} a_j y^{n-j}\right) 
\]
and more in general for $j=0,\ldots,s-1$
\be
z^{n-j} = -\frac1{\Delta t c_{-1}\lambda_I}\left(y^{n-j} +\sum_{h=0}^{s-1} a_h y^{n-j-h-1}\right). 
\label{eq:zn1}
\ee
Thus, by direct substitution into the second equation \eqref{eq:second} we obtain a discretization to \eqref{eq:so} in the form
\begin{eqnarray*}
	\left(y^{n+1} +\sum_{j=0}^{s-1} a_j y^{n-j}\right)\left(1+\lambda_R\Delta t c_{-1}\right) &=& -\sum_{j=0}^{s-1} a_j\left(y^{n-j} +\sum_{h=0}^{s-1} a_h y^{n-j-h-1}\right)\\
	&+&\lambda_I(\lambda_I -\lambda_R\varepsilon^{\alpha}) \Delta t^2 c_{-1} \sum_{j=0}^{s-1} b_j y^{n-j}. 
\end{eqnarray*}
Finally, we can rewrite the resulting scheme as
\begin{eqnarray}
\nonumber
y^{n+1} &=& -\sum_{j=0}^{s-1} a_j y^{n-j} - \frac1{1+\lambda_R\Delta t c_{-1}}\sum_{j=0}^{s-1} a_j\left(y^{n-j} +\sum_{h=0}^{s-1} a_h y^{n-j-h-1}\right)\\[-.25cm]
\label{eq:scstab}
\\[-.25cm]
\nonumber
&+&\frac{\lambda_I(\lambda_I -\lambda_R\varepsilon^{\alpha})  \Delta t^2 c_{-1}}{1+\lambda_R\Delta t c_{-1}} \sum_{j=0}^{s-1} b_j y^{n-j}. 
\end{eqnarray}
Note that, the above IMEX-LM in the limit $\varepsilon\to 0$ for $\alpha=1$ leads to the reduced scheme
\be
y^{n+1} = -\sum_{j=0}^{s-1} a_j y^{n-j} -\Delta t (i+\xi)\xi \sum_{j=0}^{s-1} b_j y^{n-j},
\ee
which corresponds to an explicit LMM for the convection-diffusion equation.

The characteristic equation for scheme \eqref{eq:scstab} reads
\be
\varrho(\zeta) + \frac1{1+z_R c_{-1}} \sigma_1(\zeta)-\frac{z_I(z_I-z_R\varepsilon^{\alpha}) c_{-1}}{1+z_R c_{-1}}\sigma_2(\zeta)=0,
\label{eq:char}
\ee
with $z_R = \lambda_R\Delta t$, $z_I = \lambda_I\Delta t$ and
\[
\varrho(\zeta) = \zeta^s + \sum_{j=0}^{s-1} a_j \zeta^{n-j},\quad \sigma_1(\zeta)=\sum_{j=0}^{s-1} a_j\left(\zeta^{n-j} +\sum_{h=0}^{s-1} a_h \zeta^{n-j-h-1}\right),\quad \sigma_2(\zeta)=\sum_{j=0}^{s-1} b_j \zeta^{n-j}.
\]
Stability corresponds to the requirement that all roots of \eqref{eq:char} have modulus less or equal one and that all multiple roots have modulus less than one.
\begin{figure}\centering
	{\includegraphics[width=6.5cm]{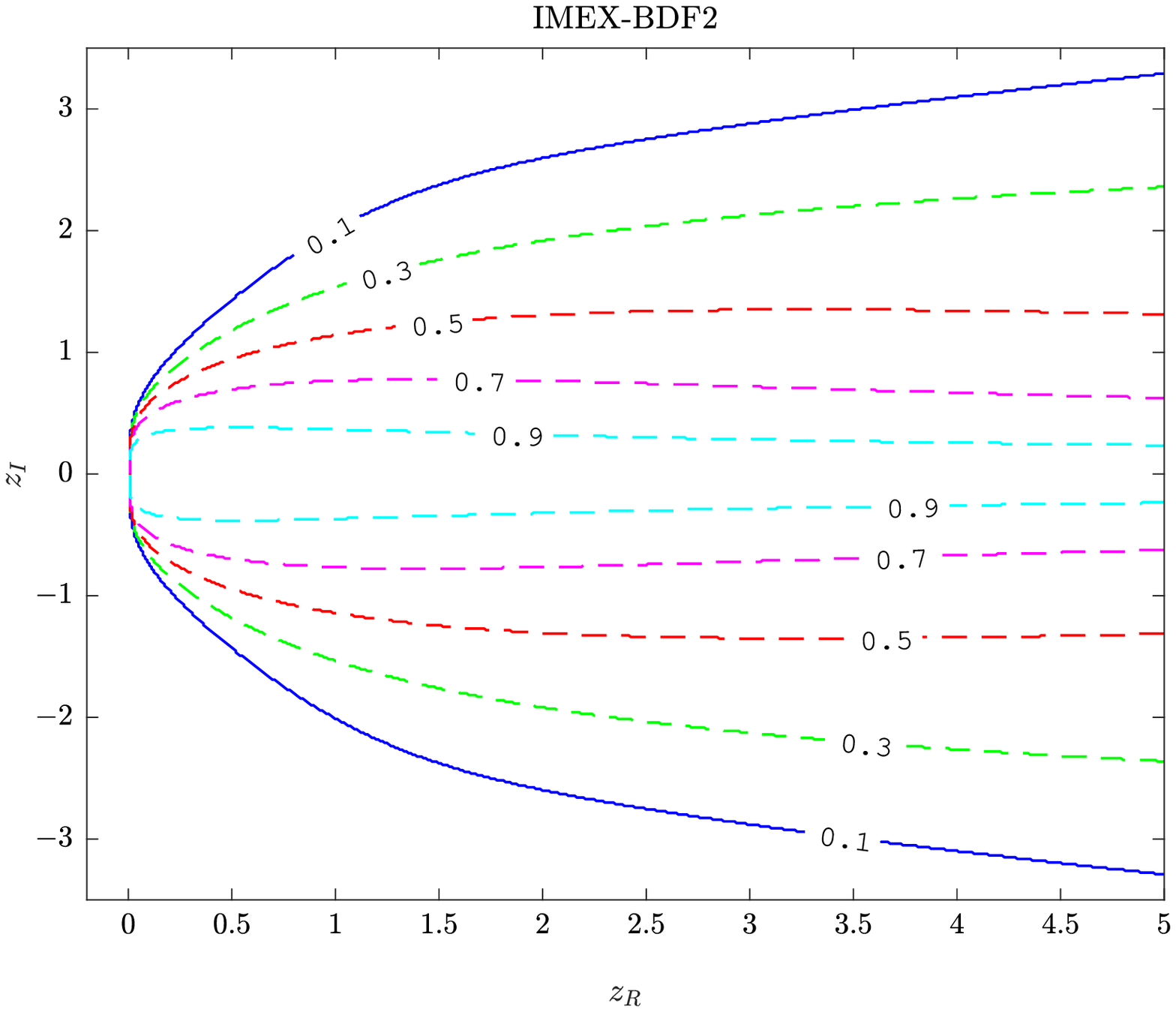}}
	\hspace{+0.15cm}
	{\includegraphics[width=6.5cm]{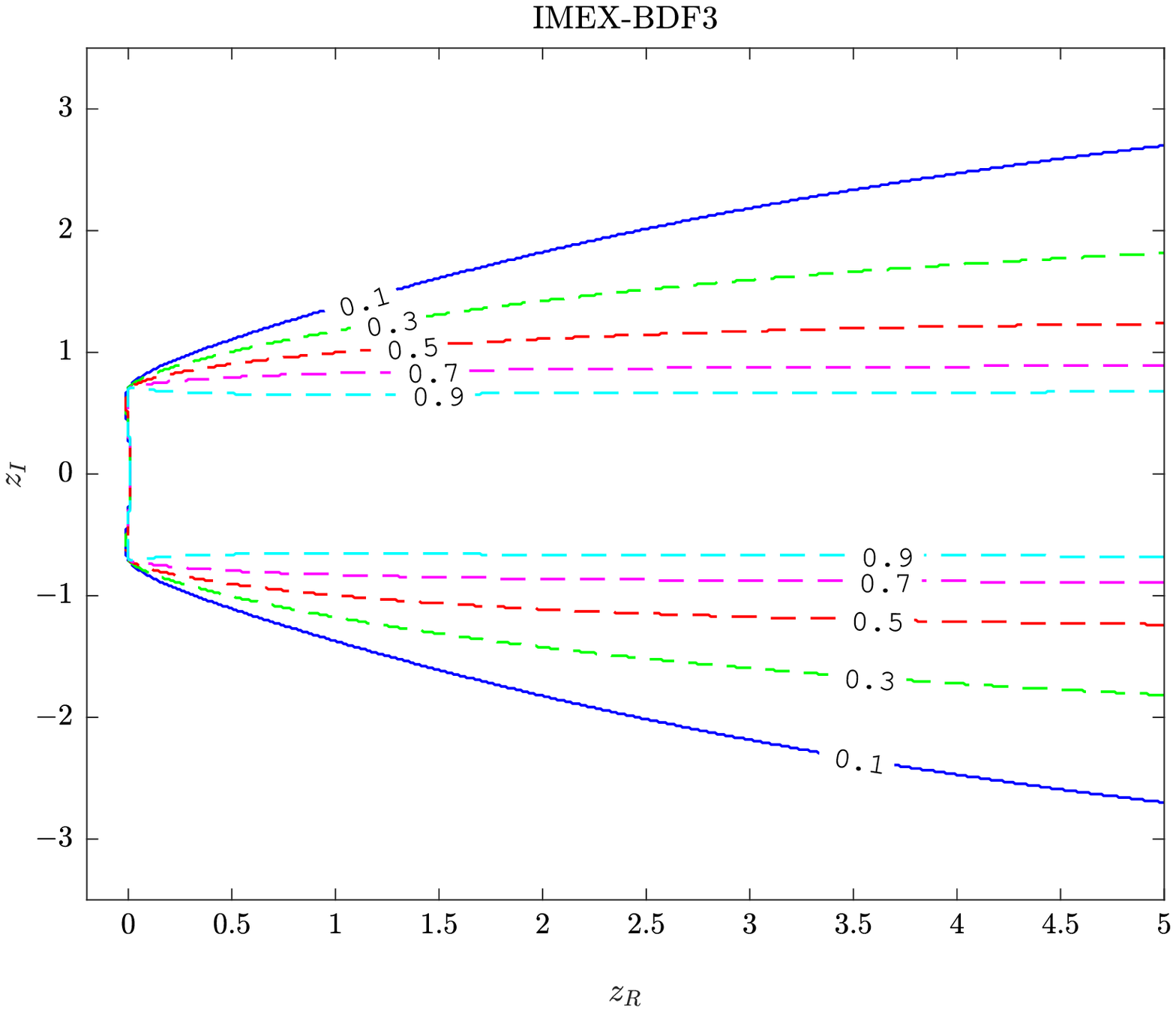}}
	\vspace{1cm}
	\\
	{\includegraphics[width=6.5cm]{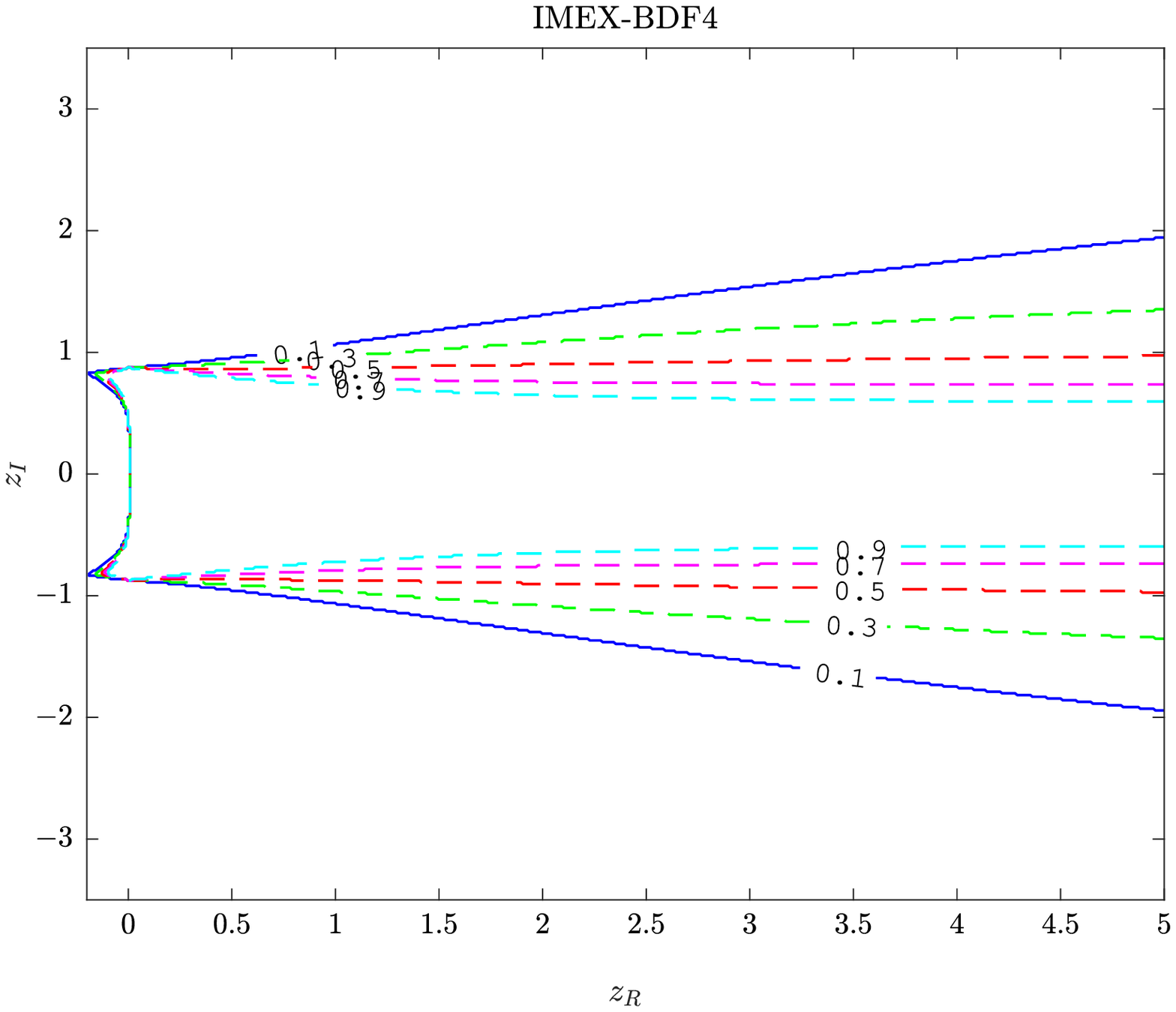}}
	\hspace{+0.15cm}
	{\includegraphics[width=6.5cm]{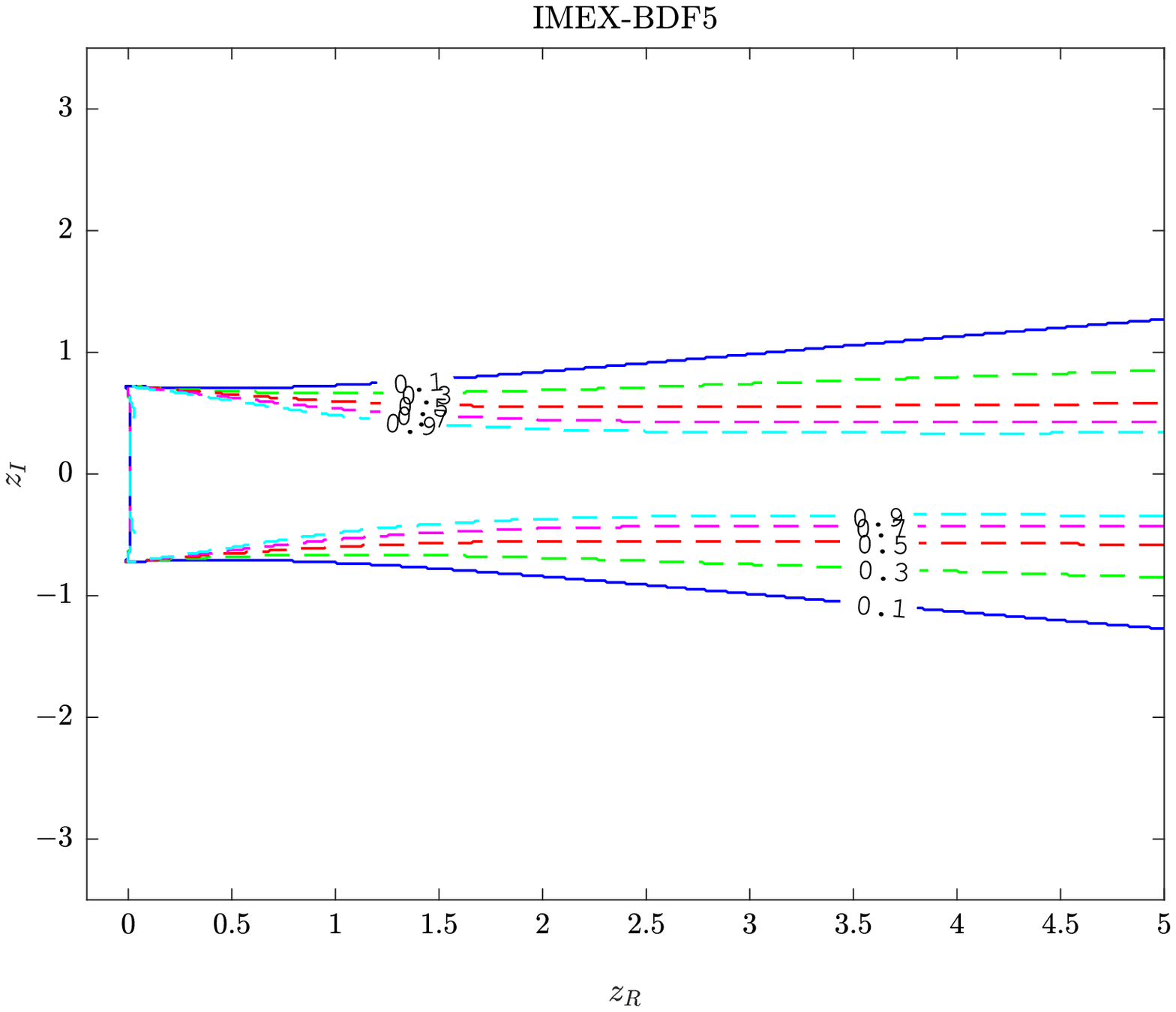}}
		\caption{AP-explicit methods. Stability regions of IMEX-BDF methods in terms of $z_I=i\xi\Delta t/\epsi^{\alpha}$ and $z_R=\Delta t/\epsi^{1+\alpha}$. The different contour lines correspond to different values of the scaling parameter $\epsi^\alpha$.}\label{fig:lambda_exp}
\end{figure}
In Figure \ref{fig:lambda_exp} we plot the stability regions of AP-explicit IMEX-BDF schemes with respect to the variable $z_R$ and $z_I$. The contour lines represent different values of the scaling parameter $\epsi^\alpha$. Note that, since we are assuming to use central differences, the stability regions are inversely proportional to the value of $\epsi^\alpha$, since $\epsi^{1-\alpha}$ measures the strength of the diffusive term in agreement with \eqref{eq:cel}. As expected, as the order of the methods increase the corresponding stability requirements become stronger and the various stability regions are reduced.
 

%
%

\begin{figure}\centering
	{\includegraphics[width=6.5cm]{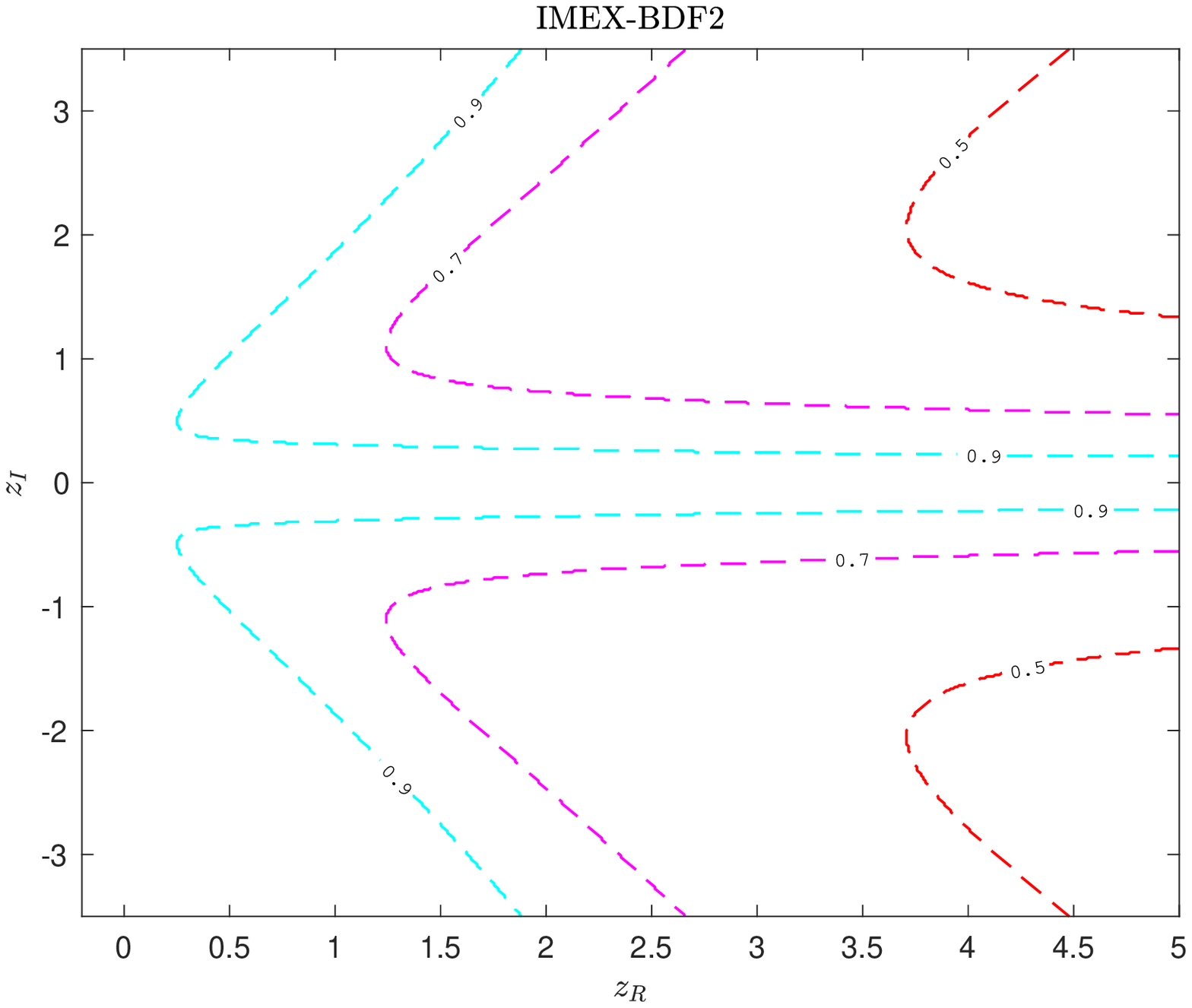}}
	\hspace{+0.15cm}
	{\includegraphics[width=6.5cm]{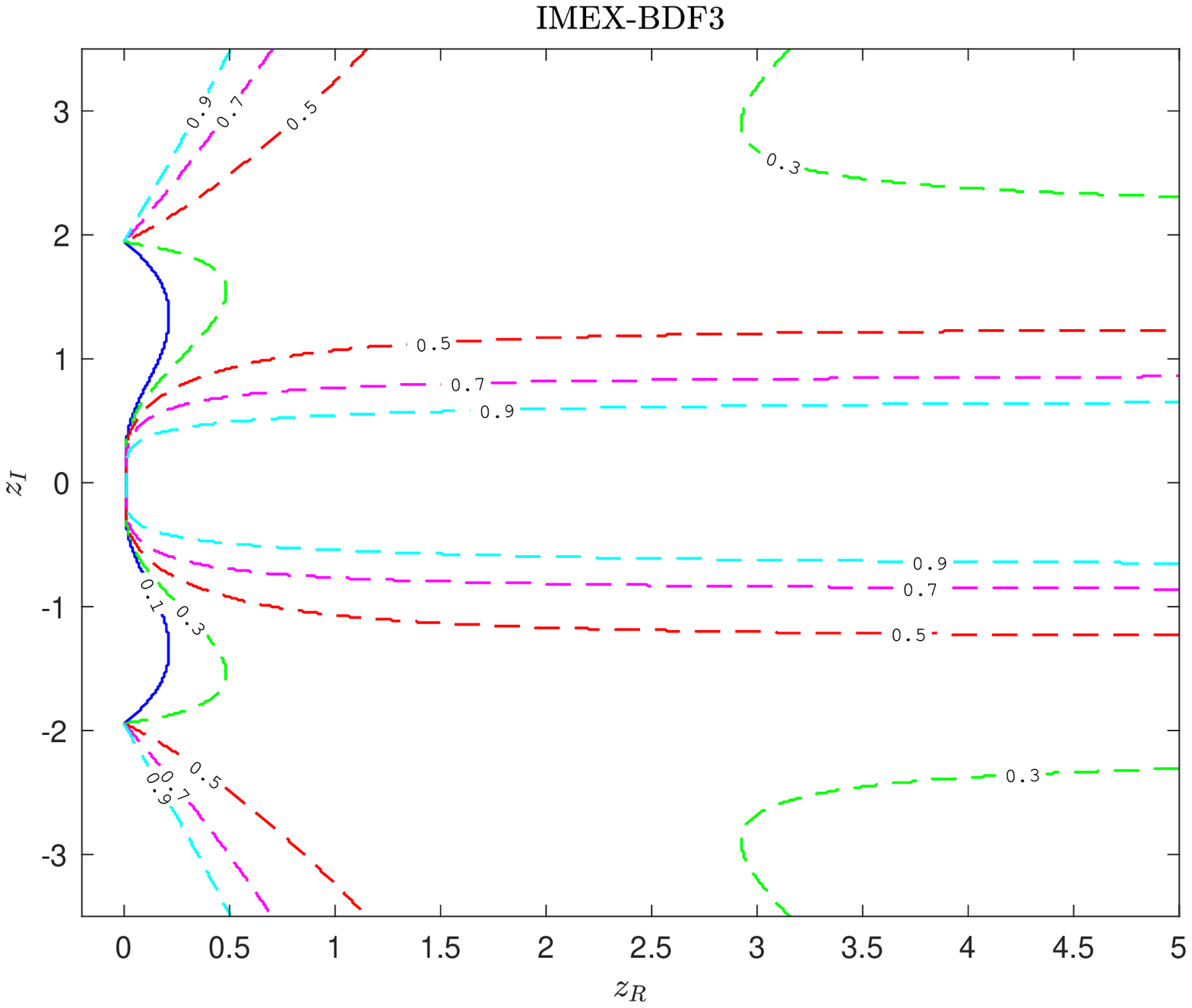}}
	\vspace{1cm}
	\\
	{\includegraphics[width=6.5cm]{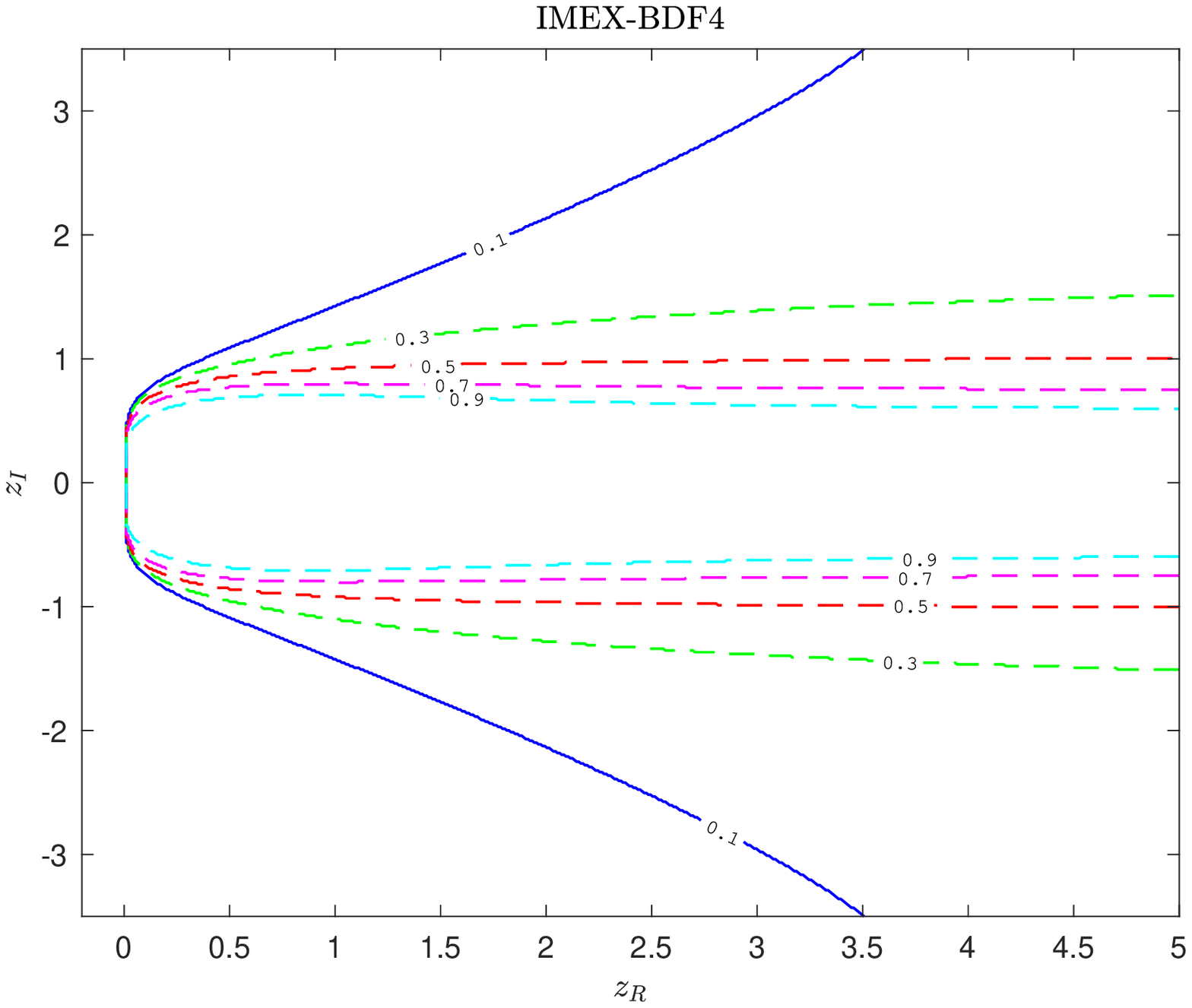}}
	\hspace{+0.15cm}
	{\includegraphics[width=6.5cm]{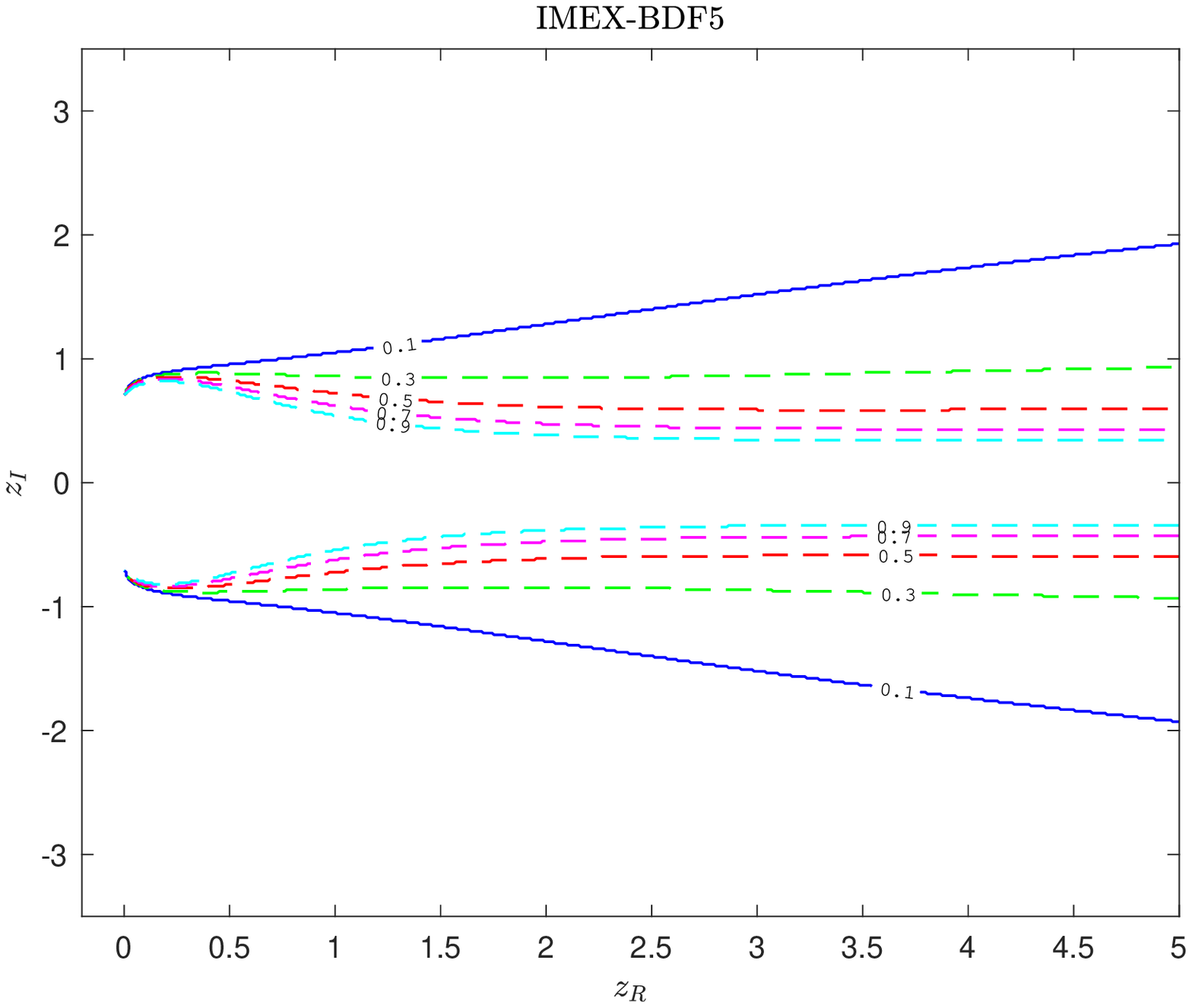}}
			\caption{AP-implicit methods. Stability regions of IMEX-BDF methods in terms of $z_I=i\xi\Delta t/\epsi^{\alpha}$ and $z_R=\Delta t/\epsi^{1+\alpha}$. The different contour lines correspond to different values of the scaling parameter $\epsi^\alpha$.}\label{fig:lambda_imp}
\end{figure}

\subsection{AP-implicit methods}
Next we apply an IMEX-LM method to \eqref{eq:tp1fd} in the AP-implicit form
\begin{eqnarray}
\label{eq:first2}
y^{n+1} &=& -\sum_{j=0}^{s-1} a_j y^{n-j} -\lambda_I \Delta t \sum_{j=-1}^{s-1} c_j z^{n-j}\\
\label{eq:second2}
z^{n+1} &=& -\sum_{j=0}^{s-1} a_j z^{n-j} +\lambda_R \varepsilon^{\alpha}  \Delta t \sum_{j=0}^{s-1} b_j y^{n-j} -\Delta t \sum_{j=-1}^{s-1} c_j (\lambda_R z^{n-j}+ \lambda_I y^{n-j}).
\end{eqnarray}
Thus, restricting to IMEX-BDF methods, by direct substitution of \eqref{eq:zn1} into the second equation \eqref{eq:second2} we get
\begin{eqnarray*}
	\left(y^{n+1} +\sum_{j=0}^{s-1} a_j y^{n-j}\right)\left(1+\lambda_R\Delta t c_{-1}\right) &=& -\sum_{j=0}^{s-1} a_j\left(y^{n-j} +\sum_{h=0}^{s-1} a_h y^{n-j-h-1}\right)\\
	&-&\lambda_I\lambda_R \Delta t^2 c_{-1}\varepsilon^{\alpha} \sum_{j=0}^{s-1} b_j y^{n-j}+(\lambda_I \Delta t c_{-1})^2 y^{n+1}, 
\end{eqnarray*}
or equivalently
\begin{eqnarray}
\nonumber
y^{n+1} &=& -\sum_{j=0}^{s-1} a_j y^{n-j}
- \frac1{1+\lambda_R\Delta t c_{-1}}\sum_{j=0}^{s-1} a_j\left(y^{n-j} +\sum_{h=0}^{s-1} a_h y^{n-j-h-1}\right)\\[-.15cm]
\label{eq:scstab2}
\\[-.3cm]
\nonumber
&-&\frac{\lambda_I\lambda_R \Delta t^2 c_{-1}\varepsilon^{\alpha}}{1+\lambda_R\Delta t c_{-1}} \sum_{j=0}^{s-1} b_j y^{n-j}
+\frac{(\lambda_I\Delta t c_{-1})^2}{1+\lambda_R\Delta t c_{-1}} y^{n+1}. 
\end{eqnarray}
Now, the above LMM method in the limit $\varepsilon\to 0$ for $\alpha=1$ leads to the scheme
\be
y^{n+1} = -\sum_{j=0}^{s-1} a_j y^{n-j} -\Delta t\, i\xi \sum_{j=0}^{s-1} b_j y^{n-j}-\Delta t \xi^2 c_{-1}y^{n+1},
\ee
which corresponds to an implicit-explicit IMEX-BDF scheme for the convection-diffusion equation.


The characteristic equation associated to scheme \eqref{eq:scstab2} takes the form
\be
\varrho(\zeta) + \frac1{1+z_R c_{-1}} \sigma_1(\zeta)+\frac{z_Iz_R\varepsilon^{\alpha} c_{-1}}{1+z_R c_{-1}}\sigma_2(\zeta)-\frac{z_I^2 c_{-1}^2}{1+z_R c_{-1}}\zeta^{s}=0.
\label{eq:char2}
\ee


We report in Figure \ref{fig:lambda_imp} the stability regions of the AP implicit IMEX-BDF methods with respect to the variable $z_R$ and $z_I$. The contour lines, as for the AP explicit case, represent different values of the scaling parameter $\epsi^\alpha$. The second order method is uniformly stable when $\epsi^\alpha<0.5$, all other methods show better stability properties compared to the AP-explicit case, in particular for large values of $z_R$. Again the stability regions diminish for increasing values of $\epsi^\alpha$, in agreement with the limit problem \eqref{eq:cel}, and as the order of the methods increases.

\section{Space discretization}

\revised{In this section we briefly discuss the space discretization adopted in the numerical examples. For the hyperbolic fluxes, we consider a WENO method of order five \cite{Shu} combined with a Rusanov flux. We stress that the space discretization is not constructed over the original discretized systems, namely \eqref{eq:SPs}, \eqref{eq:SPs_vec_a} and \eqref{eq:SPs_vec1}. Instead, we introduce the space discretization on the reformulated systems \eqref{eq:SPs_vec_ex} for the AP-explicit with $\alpha=1$, on \eqref{eq:SPs_vec_ex2_alpha} for the AP-explicit with $\alpha\in[0,1)$ and on \eqref{eq:SPs_vec_ex2_par} for the AP-implicit case. In fact, the adopted IMEX partitioning of the system, which guarantees boundedness of the eigenvalues, is of paramount importance to avoid instabilities of the fluxes and excessive numerical dissipation typical of diffusive scaling limits \cite{JL, NP, NP2}. As a consequence, the numerical diffusion is chosen accordingly to (\ref{eq:eigen_alpha}) in the numerical fluxes reported below.}

Given a generic flux function $F(Q)$ of $Q\in\mathbb{R}^n$, we first reconstruct the unknown values at the interfaces $Q_-,Q_+$ and successively we employ the numerical Rusanov flux defined as follows
\be\label{eq:Rusanov}
H(Q^-,Q^+)=\frac{1}{2}\left[F(Q^+)+F(Q^-)-\Theta(F'(q))S(Q^+-Q^-)\right], \ \Theta(F'(Q))= \max_{Q\in[Q^-,Q^+]}\{|\lambda(F'(Q))|\}
\ee
where $\max_{q\in[Q^-,Q^+]}\{|\lambda(F'(Q))|\}$ represents the maximum modulus of the eigenvalues of the Jacobian matrix $F'(Q)$ and $S\in\mathbb{R}^{n\times n}$ a transformation matrix. \revised{Hence, for systems \eqref{eq:SPs_vec_ex}, \eqref{eq:SPs_vec_ex2_alpha}, \eqref{eq:SPs_vec_ex2_par} and according to \eqref{eq:eigenvalues}, \eqref{eq:eigen_alpha},} this value will depend on the scaling parameter $\varepsilon$ and the choice of the discretization steps. In particular for the general hyperbolic system (\ref{I61}), \revised{independently on the scaling factor $\alpha$}, we have two unknowns $Q = (u,v)^T$ and three fluxes
\begin{align}
\hat f_{i+\frac{1}{2}} = \frac{1}{2}\left[(f(u^+_{i+\frac{1}{2}}) + f(u^-_{i+\frac{1}{2}})) - \Theta(u,v)(u^+_{i+\frac{1}{2}}-u^-_{i+\frac{1}{2}}) \right],\\ \hat v_{i+\frac{1}{2}} = \frac{1}{2}\left[(v^+_{i+\frac{1}{2}} + v^-_{i+\frac{1}{2}}) - \Theta(u,v)(u^+_{i+\frac{1}{2}}-u^-_{i+\frac{1}{2}}) \right]
\end{align}
and
\begin{align}
\hat p_{i+\frac{1}{2}} = \frac{1}{2}\left[(p(u^+_{i+\frac{1}{2}}) + p(u^-_{i+\frac{1}{2}})) - \Theta(u,v)(v^+_{i+\frac{1}{2}}-v^-_{i+\frac{1}{2}}) \right],
\end{align}
where to comply with \eqref{eq:Rusanov} we consider
\be
S=
\begin{bmatrix}
0 & 1\\
1 & 0
\end{bmatrix}
,\qquad 
\Theta(u,v) =\frac{1}{2}\left(\gamma(1-\theta_{\alpha})\pm\sqrt{\gamma^2(1-\theta_\alpha)^2+4\epsi^{-2\alpha}\theta_\alpha^2}\right), \ \gamma=f'(u).
\label{eq:bound}
\ee
The generic variables $w$ reconstructed at the grid interfaces $(i+\frac{1}{2})$ and $(i-\frac{1}{2})$ respectively on the right $w^-_{i+\frac{1}{2}}$ and on the left side $w^+_{i-\frac{1}{2}}$ are given by
\be w^-_{i+\frac{1}{2}}=\sum_{r=0}^{2}\omega_r w^{(r)}_{i+\frac{1}{2}}, \quad w^+_{i-\frac{1}{2}}=\sum_{r=0}^{2}\omega_r \tilde w^{(r)}_{i+\frac{1}{2}}
\ee
with weights
\be
\omega_r=\frac{\alpha_r}{\sum_{s=0}^{2}\alpha_s}, \ \alpha_r=\frac{d_r}{(\epsilon+\beta_r)^2}, \quad  \ \tilde \omega_r=\frac{\tilde \alpha_r}{\sum_{s=0}^{2}\tilde \alpha_s}, \ \tilde\alpha_r=\frac{\tilde d_r}{(\epsilon+\beta_r)^2},
\ee
and with standard smoothness indicators
\begin{align*}
\beta_0=\frac{13}{12}(w_i-2w_{i+1}+w_{i+2})^2+\frac{1}{4}(3w_{i}-4w_{i+1}+w_{i+2})^2\\
\beta_1=\frac{13}{12}(w_{i-1}-2w_{i}+w_{i+1})^2+\frac{1}{4}(w_{i-1}+w_{i+1})^2\\
\beta_2=\frac{13}{12}(w_{i-2}-2w_{i-1}+w_{i})^2+\frac{1}{4}(3w_{i-2}-4w_{i-1}+w_{i})^2,
\end{align*}
where $\epsilon=10^{-8}$, $d_0=3/10=\tilde d_2$, $d_1=3/5=\tilde d_1$ and $d_2=1/10=\tilde d_0$. Finally, the values $w^{(r)}_{i+\frac{1}{2}}$ and $w^{(r)}_{i+\frac{1}{2}}$ represent the third order reconstructions of the pointwise values $\bar w_i$. These are obtained through the formulas
\be w^{(r)}_{i+\frac{1}{2}}=\sum_{j=0}^{2}c_{rj} \bar w_{i-r+j}, \quad w^{(r)}_{-\frac{1}{2}}=\sum_{j=0}^{2}\tilde c_{rj} \bar w_{i-r+j}, \quad r=0,1,2
\ee
where $\bar w_{i-r+j}$ are the pointwise values of the unknown evaluated at the points $S_r(i)=\{x_{i-r},..,x_{i-r+2} \}$ $\ r=0,1,2$. Since, we use equispaced grid points, the coefficients $c_{rj}$ can be precomputed and their values are reported in Table \ref{tab:wenocoef}.
\begin{table}[ht]
	\centering
	\caption{Coefficients $c_{rj}$ for the $5$-th order WENO reconstruction on equispaced grid points.} \label{tab:wenocoef}
	\vspace{+0.25cm}
	{\footnotesize
		\begin{tabular}{|l|r|r|r|}
			\hline
			 r &{$j = 0$} &{$j= 1$}  &{$j= 2$}\\
			\hline
						 -1 & 11/6 & -7/6 & 1/3 \\
			\hline
						 0 & 1/3 & 5/6 & -1/6 \\
			\hline
						 1 & -1/6 & 5/6 & 1/3 \\
			\hline
						 2 & 1/3 & -7/6 & 11/6 \\
			\hline
			
	\end{tabular}}
\end{table}

\revised{In addition, we have a second order term $\partial_{xx} p(u)$ which, for the AP-explicit case, may be discretized by two consecutive application of the Rusanov flux with WENO reconstruction of the state variables and with numerical diffusion $\Theta(u,v)$ fixed equal to zero, or by a specific space discretization which is consistent with the limit problem. For example, by a sixth order finite difference formula 
\be
\partial_{xx} p(u(x_i))\simeq\frac{a p(u_{i-3})+b p(u_{i-2})+c p(u_{i-1})+dp(u_{i})+c p(u_{i+1})-b p(u_{i+2})+a p(u_{i+3})}{(\Delta x)^2}
\ee
with $a=1/90,b=-3/20,c=3/2,d=-49/18$. This latter approach has been adopted in the case of AP-implicit schemes, since the term $\partial_{xx} p(u)$ is implicit and we want to avoid nonlinearities induced by the WENO reconstructions.} 

\section{Numerical validation \& applications}\label{sec:numtest}

In this Section, we present different numerical tests to validate the analysis performed in the previous Sections. In particular, we report results for the IMEX linear multistep from order two up to order five \revised{both for the AP-explicit and for the AP-implicit formulations}. For the details about the IMEX-LM methods used we refer to Appendix \ref{app:imex_multistep}. We remark that other IMEX-LM methods can be included as well in the present formulation, see for example \cite{RSSZ1, RSSZ2}. \revised{In all test cases, the initial data has been chosen well prepared and the IMEX-LM methods have been initialized with a third order IMEX Runge-Kutta scheme (see \cite{BPR}) with a time step which satisfies the accuracy constraints.}  
\subsection{Test 1. \revised{Numerical convergence study for a linear problem.}}
We consider the following linear hyperbolic model with diffusive scaling for $\alpha=1$
\be
\left\{  
\begin{array}{l} 
\displaystyle  
\partial_t u + \partial_x v =0, \\
\displaystyle    
\partial_t v + \frac{1}{\epsi^2} \partial_x u = -\frac{1}{\epsi^{2}}(v - \gamma u),
\\
\end{array}
\right. 
\label{T1}
\ee
where $\gamma > 0$.
In the diffusive limit $\varepsilon \to 0 $ the second equation relaxes to the local equilibrium 
\[
\j = \gamma\rho - {\partial_x \rho},
\]
substituting into the first equation this gives the limiting advection-diffusion equation
\be
\partial_{t}\rho + \gamma \partial_{x} \rho = \partial_{xx} \rho. 
\label{T1relax}
\ee
\revised{
In particular, we consider the model \eqref{T1} solved on the domain $x\in[0,1]$, with $\gamma=1$ and periodic boundary conditions and with smooth initial data given by
\be
u(x,0) = \sin(2 x\pi), \ 
v(x,0) = \sin(2 x\pi)-\cos(2 x \pi).
\ee
Note that, the initial data is well prepared, in the sense that $v(x,0)=u(x,0)-\partial_x u(x,0)$.}
For this specific problem, we numerically estimate the order of convergence of the schemes for various values of $\varepsilon=1, 0.1, 0.01, 0.001$ by measuring the space and time $L_1$-error of the numerical solutions computed by using as reference solution the thinner grids. Namely, given a coarser grid $\Delta x_1=1/N$ with $N=128$ we consider 
\be \Delta x_{k+1}=\Delta x_{k}/2 \ \text{with} \ k=1,..,4.\ee 
\revised{The time step for AP-explicit methods is chosen as $\Delta t = \lambda\Delta x\max\left\{\varepsilon, \Delta x\right\}$, with $\lambda=0.25$.
Namely, the largest between the CFL condition imposed by the hyperbolic part and by the limiting parabolic part of the equations. Instead, for implicit schemes we choose $\Delta t = \lambda\Delta x\max\left\{\varepsilon, 1\right\}$ since the diffusion term is integrated implicitly in the limit.} 

The local truncation error is measured for the two components $u$ and $v$ as follows
\[
E^k_{\Delta x,\Delta t}(u) = |u^k_{\Delta x,\Delta t}(\cdot,T) - u^{k-1}_{\Delta x,\Delta t}(\cdot,T)|, \
E^k_{\Delta x,\Delta t}(v) = |v^k_{\Delta x,\Delta t}(\cdot,T) - v^{k-1}_{\Delta x,\Delta t}(\cdot,T)|,
\]
and the order of convergence is estimated by computing the rate between two $L_1$ errors of distinct numerical solutions 
\be\label{eq:conv_rate}
\texttt{p}_k(u) = \log_2(\| E^{k-1}_{\Delta x,\Delta t}(u)\|_{1}/ \|E_{\Delta x,\Delta t}^{k}(u) \|_{1}), \ 
\texttt{p}_k(v) = \log_2(\| E^{k-1}_{\Delta x,\Delta t}(v)\|_{1}/ \|E_{\Delta x,\Delta t}^{k}(v) \|_{1}).
\ee
The analysis is performed for several different IMEX linear multi-step schemes from second to fifth order (see Appendix A). In particular, we focus on the BDF and the TVB classes of IMEX multistep methods thanks to their favorable stability properties (see \cite{HR, HRS} for details and derivation).
\revised{We report in Table \ref{tab:order_u} the space-time $L_1$  errors and the relative rates of convergence for increasing size of the meshes considering $N=2^k$ points with $k = 8,..,11$ for the $u$ variable, while the space-time $L_1$  errors and the relative rates of convergence for the $v$ variable are shown in Table \ref{tab:order_v} for the AP-explicit schemes. In Table \ref{tab:order_u1} and \ref{tab:order_v1}, we report the corresponding $L_1$ errors and rates of convergence for the AP-implicit schemes. In all cases we can conclude that the expected orders of convergence are achieved by the schemes for the different values of the asymptotic parameter $\epsi$ and that the behavior of the schemes outperforms the corresponding IMEX Runge-Kutta methods for the non conserved quantity $v$ (see Table 2 in \cite{BPR17} for example). In particular, we observe the tendency of $4$-th order methods to achieve higher then expected convergence rates on the conserved quantity $u$ for moderate values of the stiffness parameter. The same tendency, on both variables $u$ and $v$, is observed for $5$-th order methods close to the diffusion limit particularly in the AP-implicit case.
On the contrary, the $4$-th order schemes in the AP-explicit implementation produce a slight deterioration on the non conserved quantity $v$ which is not observed in the AP-implicit setting. The SG($3$,$2$) scheme also suffers of a slight deterioration of accuracy close to the diffusion limit in the AP-explicit form. It should be noted that scheme in AP-explicit form close to the fluid limit have a smaller truncation error with respect to time thanks to the CFL condition $\Delta t = O(\Delta x^2)$. This can be observed by comparing the $L_1$-errors of the AP-explicit and AP-explicit formulations for $\varepsilon=0.001$.  
Finally, in Figure \ref{fig:T1_order} and \ref{fig:T1_order_impl}, we summarized in a plot the order of convergence for the mesh corresponding to 256 nodes as a function of $\varepsilon$ for the AP-explicit and AP-implicit methods. These plots emphasize that the convergence rate for all schemes is almost uniform.
}

\begin{table}[ht]
 \centering
 \caption{$L^1$ error and estimated convergence rates for $u$ \revised{in the AP-explicit case.}} \label{tab:order_u}
\vspace{+0.25cm}
{\footnotesize
 \begin{tabular}{|c|l|rr|rr|rr|rr|}
 	\hline
~& ~ & \multicolumn{2}{c|}{$\epsi = 1$} &\multicolumn{2}{c|}{$\epsi = 0.1$}  &\multicolumn{2}{c|}{$\epsi = 0.01$} & \multicolumn{2}{c|}{$\epsi = 0.001$} \\
 \hline
IMEX& \multicolumn{1}{c|}{$N$} & \multicolumn{1}{c}{ $\|E_{\Delta x,\Delta t}^{k} \|_{1}$} & \multicolumn{1}{c|}{Rate} & \multicolumn{1}{c}{$\|E_{\Delta x,\Delta t}^{k} \|_{1}$} & \multicolumn{1}{c|}{ Rate} & \multicolumn{1}{c}{$\|E_{\Delta x,\Delta t}^{k} \|_{1}$} & \multicolumn{1}{c|}{ Rate} & \multicolumn{1}{c}{$\|E_{\Delta x,\Delta t}^{k} \|_{1}$} & \multicolumn{1}{c|}{ Rate} \\
 \hline
\hline
\hline
\parbox[t]{4.5mm}{\multirow{3}{*}{\rotatebox[]{90}{\parbox[c]{1.2cm}{\footnotesize IMEX-SG(3,2) }}}}
& 128 & 9.9402e-05& -- & 8.2324e-05& -- &0.00019196& -- &1.0662e-07& -- \\
& 256 & 3.0565e-05& 1.7014& 2.5215e-05& 1.707& 5.7846e-05& 1.7305& 4.262e-08& 1.3229\\
& 512 & 8.3499e-06& 1.8721& 6.8755e-06& 1.8748& 1.5671e-05& 1.8841& 1.1871e-08& 1.8441\\
\hline
\hline
\parbox[t]{4.5mm}{\multirow{3}{*}{\rotatebox[]{90}{\parbox[c]{1.2cm}{\footnotesize IMEX-BDF2 }}}}
& 128 & 4.4737e-05& -- & 3.6789e-05& -- &8.9909e-05& -- &9.2432e-08& -- \\
& 256 & 1.4423e-05& 1.6331& 1.1759e-05& 1.6455& 2.8204e-05& 1.6726& 2.4676e-08& 1.9053\\
& 512 & 4.0033e-06& 1.8491& 3.2494e-06& 1.8555& 7.7212e-06& 1.869& 6.0694e-09& 2.0235\\
\hline
\hline
\parbox[t]{4.5mm}{\multirow{3}{*}{\rotatebox[]{90}{\parbox[c]{1.2cm}{\footnotesize IMEX-TVB(3,3) }}}}
& 128 & 7.3451e-08& -- & 1.2994e-07& -- &5.0813e-06& -- &3.1308e-08& -- \\
& 256 & 1.2953e-08& 2.5035& 1.8201e-08& 2.8358& 7.5974e-07& 2.7416& 3.8147e-09& 3.0369\\
& 512 & 1.8478e-09& 2.8094& 2.4055e-09& 2.9197& 1.0251e-07& 2.8898& 4.4859e-10& 3.0881\\
\hline
\hline
\parbox[t]{4.5mm}{\multirow{3}{*}{\rotatebox[]{90}{\parbox[c]{1.2cm}{\footnotesize IMEX-BDF3 }}}}
& 128 & 8.7902e-08& -- & 1.6714e-07& -- &6.6227e-06& -- &5.8099e-08& -- \\
& 256 & 1.5388e-08& 2.514& 2.4439e-08& 2.7738& 9.8974e-07& 2.7423& 6.0486e-09& 3.2638\\
& 512 & 2.1902e-09& 2.8127& 3.2815e-09& 2.8968& 1.3293e-07& 2.8964& 6.5821e-10& 3.2\\
\hline
\hline
\parbox[t]{4.5mm}{\multirow{3}{*}{\rotatebox[]{90}{\parbox[c]{1.2cm}{\footnotesize IMEX-TVB(4,4) }}}}
& 128 & 2.2234e-08& -- & 6.9585e-09& -- &1.2029e-06& -- &2.0263e-07& -- \\
& 256 & 8.4006e-10& 4.7262& 1.7448e-10& 5.3177& 9.1633e-08& 3.7145& 5.8212e-09& 5.1214\\
& 512 & 4.7195e-11& 4.1538& 3.7747e-12& 5.5305& 6.0883e-09& 3.9118& 2.2316e-10& 4.7052\\
\hline
\hline
\parbox[t]{4.5mm}{\multirow{3}{*}{\rotatebox[]{90}{\parbox[c]{1.2cm}{\footnotesize IMEX-BDF4 }}}}
& 128 & 2.3209e-08& -- & 7.4749e-09& -- &4.5554e-07& -- &2.449e-08& -- \\
& 256 & 6.8151e-10& 5.0898& 2.1477e-10& 5.1212& 3.3431e-08& 3.7683& 1.1948e-09& 4.3573\\
& 512 & 2.9394e-11& 4.5352& 5.499e-12& 5.2875& 2.2255e-09& 3.909& 6.3684e-11& 4.2297\\
\hline
\hline
\parbox[t]{4.5mm}{\multirow{3}{*}{\rotatebox[]{90}{\parbox[c]{1.2cm}{\footnotesize IMEX-TVB(5,5) }}}}
& 128 & 1.043e-08& -- & 4.1161e-08& -- &9.4337e-09& -- &1.9208e-07& -- \\
& 256 & 3.2924e-10& 4.9854& 1.2742e-09& 5.0137& 3.3621e-10& 4.8104& 2.4637e-09& 6.2848\\
& 512 & 1.4718e-11& 4.4835& 4.7629e-11& 4.7416& 1.0529e-11& 4.9969& 5.4809e-11& 5.4903\\
\hline
\hline
\parbox[t]{4.5mm}{\multirow{3}{*}{\rotatebox[]{90}{\parbox[c]{1.2cm}{\footnotesize IMEX-BDF5 }}}}
& 128 & 4.9951e-08& -- & 3.49e-08& -- &4.906e-09& -- &1.7525e-08& -- \\
& 256 & 1.5618e-09& 4.9992& 1.0894e-09& 5.0016& 1.5768e-10& 4.9594& 4.2623e-10& 5.3616\\
& 512 & 5.6577e-11& 4.7869& 3.2187e-11& 5.0809& 4.6748e-12& 5.076& 1.1341e-11& 5.232\\
\hline
 \end{tabular}}
 \end{table}

\begin{table}[ht]
	\centering
	\caption{$L^1$ error and estimated convergence rates for $v$ \revised{in the AP-explicit case.}} \label{tab:order_v}
	\vspace{+0.25cm}
	\footnotesize{
	\begin{tabular}{|c|l|rr|rr|rr|rr|}
		\hline
		~& ~ & \multicolumn{2}{c|}{$\epsi = 1$} &\multicolumn{2}{c|}{$\epsi = 0.1$}  &\multicolumn{2}{c|}{$\epsi = 0.01$} & \multicolumn{2}{c|}{$\epsi = 0.001$} \\
		\hline
		IMEX& \multicolumn{1}{c|}{$N$} & \multicolumn{1}{c}{ $\|E_{\Delta x,\Delta t}^{k} \|_{1}$} & \multicolumn{1}{c|}{Rate} & \multicolumn{1}{c}{$\|E_{\Delta x,\Delta t}^{k} \|_{1}$} & \multicolumn{1}{c|}{ Rate} & \multicolumn{1}{c}{$\|E_{\Delta x,\Delta t}^{k} \|_{1}$} & \multicolumn{1}{c|}{ Rate} & \multicolumn{1}{c}{$\|E_{\Delta x,\Delta t}^{k} \|_{1}$} & \multicolumn{1}{c|}{ Rate} \\
		\hline
\hline
\hline
\parbox[t]{4.5mm}{\multirow{3}{*}{\rotatebox[]{90}{\parbox[c]{1.2cm}{\footnotesize IMEX-SG(3,2) }}}}
& 128 & 4.1297e-05& -- & 0.0015868& -- &0.30746& -- &0.017621& -- \\
& 256 & 1.3046e-05& 1.6624& 0.00049005& 1.6951& 0.092824& 1.7278& 0.006956& 1.341\\
& 512 & 3.6092e-06& 1.8538& 0.00013416& 1.869& 0.025173& 1.8826& 0.0019759& 1.8158\\
\hline
\hline
\parbox[t]{4.5mm}{\multirow{3}{*}{\rotatebox[]{90}{\parbox[c]{1.2cm}{\footnotesize IMEX-BDF2 }}}}
& 128 & 5.3468e-05& -- & 0.0010541& -- &0.14541& -- &0.015341& -- \\
& 256 & 1.7643e-05& 1.5996& 0.00034308& 1.6194& 0.045795& 1.6669& 0.0040821& 1.91\\
& 512 & 4.9529e-06& 1.8328& 9.5659e-05& 1.8426& 0.012564& 1.8658& 0.001048& 1.9616\\
\hline
\hline
\parbox[t]{4.5mm}{\multirow{3}{*}{\rotatebox[]{90}{\parbox[c]{1.2cm}{\footnotesize IMEX-TVB(3,3) }}}}
& 128 & 7.2288e-07& -- & 8.2789e-06& -- &0.0083555& -- &0.0044838& -- \\
& 256 & 1.1051e-07& 2.7096& 1.2742e-06& 2.6999& 0.0012579& 2.7317& 0.00064242& 2.8031\\
& 512 & 1.5103e-08& 2.8712& 1.7412e-07& 2.8714& 0.0001703& 2.8849& 8.423e-05& 2.9311\\
\hline
\hline
\parbox[t]{4.5mm}{\multirow{3}{*}{\rotatebox[]{90}{\parbox[c]{1.2cm}{\footnotesize IMEX-BDF3 }}}}
& 128 & 6.3076e-07& -- & 8.0827e-06& -- &0.010784& -- &0.0087392& -- \\
& 256 & 9.8098e-08& 2.6848& 1.2676e-06& 2.6727& 0.0016207& 2.7342& 0.00099636& 3.1328\\
& 512 & 1.3503e-08& 2.861& 1.7436e-07& 2.862& 0.0002183& 2.8923& 0.00011713& 3.0885\\
\hline
\hline
\parbox[t]{4.5mm}{\multirow{3}{*}{\rotatebox[]{90}{\parbox[c]{1.2cm}{\footnotesize IMEX-TVB(4,4) }}}}
& 128 & 8.1836e-08& -- & 8.5455e-08& -- &0.0019902& -- &0.033128& -- \\
& 256 & 3.443e-09& 4.571& 7.8768e-09& 3.4395& 0.00015204& 3.7104& 0.00099626& 5.0554\\
& 512 & 1.5925e-10& 4.4343& 6.285e-10& 3.6476& 1.0128e-05& 3.9079& 4.1573e-05& 4.5828\\
\hline
\hline
\parbox[t]{4.5mm}{\multirow{3}{*}{\rotatebox[]{90}{\parbox[c]{1.2cm}{\footnotesize IMEX-BDF4 }}}}
& 128 & 6.8621e-08& -- & 2.6004e-08& -- &0.00076076& -- &0.0045931& -- \\
& 256 & 2.4028e-09& 4.8359& 2.607e-09& 3.3183& 5.5621e-05& 3.7737& 0.00022733& 4.3366\\
& 512 & 5.3325e-11& 5.4938& 2.1977e-10& 3.5683& 3.7004e-06& 3.9099& 1.2649e-05& 4.1676\\
\hline
\hline
\parbox[t]{4.5mm}{\multirow{3}{*}{\rotatebox[]{90}{\parbox[c]{1.2cm}{\footnotesize IMEX-TVB(5,5) }}}}
& 128 & 1.7977e-08& -- & 2.1366e-07& -- &9.0073e-06& -- &0.030281& -- \\
& 256 & 4.6141e-10& 5.2839& 6.7708e-09& 4.9798& 3.4995e-07& 4.6858& 0.00039906& 6.2457\\
& 512 & 6.7507e-11& 2.773& 1.182e-10& 5.84& 1.1336e-08& 4.9481& 9.4581e-06& 5.3989\\
\hline
\hline
\parbox[t]{4.5mm}{\multirow{3}{*}{\rotatebox[]{90}{\parbox[c]{1.2cm}{\footnotesize IMEX-BDF5 }}}}
& 128 & 1.2276e-07& -- & 1.533e-07& -- &1.6471e-06& -- &0.0022391& -- \\
& 256 & 3.8413e-09& 4.9981& 4.7927e-09& 4.9993& 5.7491e-08& 4.8405& 5.2309e-05& 5.4197\\
& 512 & 2.5154e-10& 3.9327& 1.9215e-10& 4.6405& 1.4584e-09& 5.3009& 1.3736e-06& 5.251\\
\hline
\end{tabular}}
\end{table}


\begin{table}[ht]
	\centering
	\caption{$L^1$ error and estimated convergence rates for $u$ \revised{in the AP-implicit case.}} \label{tab:order_u1}
	\vspace{+0.25cm}
	{\footnotesize
		\begin{tabular}{|c|l|rr|rr|rr|rr|}
			\hline
			~& ~ & \multicolumn{2}{c|}{$\epsi = 1$} &\multicolumn{2}{c|}{$\epsi = 0.1$}  &\multicolumn{2}{c|}{$\epsi = 0.01$} & \multicolumn{2}{c|}{$\epsi = 0.001$} \\
			\hline
			IMEX& \multicolumn{1}{c|}{$N$} & \multicolumn{1}{c}{ $\|E_{\Delta x,\Delta t}^{k} \|_{1}$} & \multicolumn{1}{c|}{Rate} & \multicolumn{1}{c}{$\|E_{\Delta x,\Delta t}^{k} \|_{1}$} & \multicolumn{1}{c|}{ Rate} & \multicolumn{1}{c}{$\|E_{\Delta x,\Delta t}^{k} \|_{1}$} & \multicolumn{1}{c|}{ Rate} & \multicolumn{1}{c}{$\|E_{\Delta x,\Delta t}^{k} \|_{1}$} & \multicolumn{1}{c|}{ Rate} \\
			\hline
			\hline
\hline
\parbox[t]{4.5mm}{\multirow{3}{*}{\rotatebox[]{90}{\parbox[c]{1.2cm}{\footnotesize IMEX-SG(3,2) }}}}
& 128 & 0.0001009& -- & 8.4701e-05& -- &0.00019279& -- &0.015388& -- \\
& 256 & 3.1168e-05& 1.6948& 2.6145e-05& 1.6958& 5.8207e-05& 1.7278& 0.0047058& 1.7093\\
& 512 & 8.5342e-06& 1.8687& 7.1564e-06& 1.8692& 1.5785e-05& 1.8826& 0.0012828& 1.8751\\
\hline
\hline
\parbox[t]{4.5mm}{\multirow{3}{*}{\rotatebox[]{90}{\parbox[c]{1.2cm}{\footnotesize IMEX-BDF2 }}}}
& 128 & 4.5983e-05& -- & 3.8769e-05& -- &9.0615e-05& -- &0.0074794& -- \\
& 256 & 1.4986e-05& 1.6175& 1.2626e-05& 1.6185& 2.8544e-05& 1.6666& 0.0023976& 1.6413\\
& 512 & 4.1817e-06& 1.8414& 3.5217e-06& 1.8421& 7.8321e-06& 1.8657& 0.00066348& 1.8535\\
\hline
\hline
\parbox[t]{4.5mm}{\multirow{3}{*}{\rotatebox[]{90}{\parbox[c]{1.2cm}{\footnotesize IMEX-TVB(3,3) }}}}
& 128 & 5.6622e-08& -- & 1.5426e-07& -- &5.1657e-06& -- &0.00024864& -- \\
& 256 & 9.4982e-09& 2.5756& 2.2969e-08& 2.7476& 7.7718e-07& 2.7326& 3.7993e-05& 2.7103\\
& 512 & 1.3103e-09& 2.8577& 3.121e-09& 2.8796& 1.0522e-07& 2.8849& 5.1683e-06& 2.878\\
\hline
\hline
\parbox[t]{4.5mm}{\multirow{3}{*}{\rotatebox[]{90}{\parbox[c]{1.2cm}{\footnotesize IMEX-BDF3 }}}}
& 128 & 7.3382e-08& -- & 1.8956e-07& -- &6.7012e-06& -- &0.0002871& -- \\
& 256 & 1.2326e-08& 2.5737& 2.8971e-08& 2.7099& 1.0064e-06& 2.7352& 4.443e-05& 2.6919\\
& 512 & 1.7158e-09& 2.8448& 3.9775e-09& 2.8647& 1.3555e-07& 2.8923& 6.0575e-06& 2.8747\\
\hline
\hline
\parbox[t]{4.5mm}{\multirow{3}{*}{\rotatebox[]{90}{\parbox[c]{1.2cm}{\footnotesize IMEX-TVB(4,4) }}}}
& 128 & 2.9453e-08& -- & 2.5367e-08& -- &1.2204e-06& -- &2.9521e-05& -- \\
& 256 & 1.2448e-09& 4.5645& 7.248e-10& 5.1292& 9.3626e-08& 3.7044& 2.2475e-06& 3.7153\\
& 512 & 5.0075e-11& 4.6356& 1.6962e-11& 5.4172& 6.2417e-09& 3.9069& 1.5677e-07& 3.8416\\
\hline
\hline
\parbox[t]{4.5mm}{\multirow{3}{*}{\rotatebox[]{90}{\parbox[c]{1.2cm}{\footnotesize IMEX-BDF4 }}}}
& 128 & 2.582e-08& -- & 7.3488e-09& -- &4.6271e-07& -- &1.2261e-05& -- \\
& 256 & 8.5367e-10& 4.9187& 2.0595e-10& 5.1572& 3.4155e-08& 3.7599& 8.3808e-07& 3.8709\\
& 512 & 1.6737e-11& 5.6726& 1.69e-12& 6.9291& 2.2839e-09& 3.9025& 5.5639e-08& 3.9129\\
\hline
\hline
\parbox[t]{4.5mm}{\multirow{3}{*}{\rotatebox[]{90}{\parbox[c]{1.2cm}{\footnotesize IMEX-TVB(5,5) }}}}
& 128 & 3.6171e-08& -- & 3.5169e-08& -- &9.4689e-09& -- &1.1944e-06& -- \\
& 256 & 1.1333e-09& 4.9962& 1.1385e-09& 4.9491& 3.3274e-10& 4.8308& 7.5898e-09& 7.298\\
& 512 & 5.1332e-11& 4.4646& 6.4841e-11& 4.1341& 1.2796e-11& 4.7006& 3.7672e-10& 4.3325\\
\hline
\hline
\parbox[t]{4.5mm}{\multirow{3}{*}{\rotatebox[]{90}{\parbox[c]{1.2cm}{\footnotesize IMEX-BDF5 }}}}
& 128 & 2.0507e-08& -- & 2.7439e-08& -- &4.9118e-09& -- &1.3221e-06& -- \\
& 256 & 6.2502e-10& 5.0361& 8.5434e-10& 5.0053& 1.5808e-10& 4.9575& 2.1912e-08& 5.915\\
& 512 & 2.5331e-11& 4.625& 3.7038e-11& 4.5277& 1.5295e-12& 6.6914& 1.0639e-09& 4.3643\\
\hline

	\end{tabular}}
\end{table}

\begin{table}[ht]
	\centering
	\caption{$L^1$ error and estimated convergence rates for $v$ \revised{in the AP-implicit case.}} \label{tab:order_v1}
	\vspace{+0.25cm}
	\footnotesize{
		\begin{tabular}{|c|l|rr|rr|rr|rr|}
			\hline
			~& ~ & \multicolumn{2}{c|}{$\epsi = 1$} &\multicolumn{2}{c|}{$\epsi = 0.1$}  &\multicolumn{2}{c|}{$\epsi = 0.01$} & \multicolumn{2}{c|}{$\epsi = 0.001$} \\
			\hline
			IMEX& \multicolumn{1}{c|}{$N$} & \multicolumn{1}{c}{ $\|E_{\Delta x,\Delta t}^{k} \|_{1}$} & \multicolumn{1}{c|}{Rate} & \multicolumn{1}{c}{$\|E_{\Delta x,\Delta t}^{k} \|_{1}$} & \multicolumn{1}{c|}{ Rate} & \multicolumn{1}{c}{$\|E_{\Delta x,\Delta t}^{k} \|_{1}$} & \multicolumn{1}{c|}{ Rate} & \multicolumn{1}{c}{$\|E_{\Delta x,\Delta t}^{k} \|_{1}$} & \multicolumn{1}{c|}{ Rate} \\
			\hline
			\hline
\hline
\parbox[t]{4.5mm}{\multirow{3}{*}{\rotatebox[]{90}{\parbox[c]{1.2cm}{\footnotesize IMEX-SG(3,2) }}}}
& 128 & 4.4648e-05& -- & 0.00089276& -- &0.30567& -- &0.11134& -- \\
& 256 & 1.3846e-05& 1.6891& 0.00027568& 1.6953& 0.092283& 1.7278& 0.032643& 1.7702\\
& 512 & 3.7984e-06& 1.866& 7.5475e-05& 1.8689& 0.025026& 1.8826& 0.0087961& 1.8918\\
\hline
\hline
\parbox[t]{4.5mm}{\multirow{3}{*}{\rotatebox[]{90}{\parbox[c]{1.2cm}{\footnotesize IMEX-BDF2 }}}}
& 128 & 2.3413e-05& -- & 0.00042114& -- &0.14368& -- &0.060252& -- \\
& 256 & 7.6904e-06& 1.6062& 0.00013723& 1.6177& 0.045257& 1.6667& 0.0174& 1.7919\\
& 512 & 2.1539e-06& 1.8361& 3.8288e-05& 1.8416& 0.012418& 1.8657& 0.0046647& 1.8992\\
\hline
\hline
\parbox[t]{4.5mm}{\multirow{3}{*}{\rotatebox[]{90}{\parbox[c]{1.2cm}{\footnotesize IMEX-TVB(3,3) }}}}
& 128 & 2.0289e-07& -- & 1.2756e-06& -- &0.0081773& -- &0.001989& -- \\
& 256 & 3.279e-08& 2.6294& 1.9536e-07& 2.707& 0.0012317& 2.731& 0.00026696& 2.8974\\
& 512 & 4.5301e-09& 2.8557& 2.6669e-08& 2.8729& 0.00016679& 2.8845& 3.4663e-05& 2.9451\\
\hline
\hline
\parbox[t]{4.5mm}{\multirow{3}{*}{\rotatebox[]{90}{\parbox[c]{1.2cm}{\footnotesize IMEX-BDF3 }}}}
& 128 & 2.3141e-07& -- & 1.3649e-06& -- &0.010612& -- &0.0023402& -- \\
& 256 & 3.8064e-08& 2.6039& 2.1348e-07& 2.6766& 0.0015951& 2.7339& 0.00032067& 2.8674\\
& 512 & 5.2868e-09& 2.848& 2.9369e-08& 2.8617& 0.00021487& 2.8921& 4.2074e-05& 2.9301\\
\hline
\hline
\parbox[t]{4.5mm}{\multirow{3}{*}{\rotatebox[]{90}{\parbox[c]{1.2cm}{\footnotesize IMEX-TVB(4,4) }}}}
& 128 & 7.1151e-08& -- & 9.5652e-08& -- &0.0019469& -- &0.00028374& -- \\
& 256 & 2.6963e-09& 4.7218& 2.8851e-09& 5.0511& 0.0001488& 3.7098& 1.6152e-05& 4.1348\\
& 512 & 1.0634e-10& 4.6642& 3.3063e-10& 3.1253& 9.9148e-06& 3.9076& 1.0431e-06& 3.9527\\
\hline
\hline
\parbox[t]{4.5mm}{\multirow{3}{*}{\rotatebox[]{90}{\parbox[c]{1.2cm}{\footnotesize IMEX-BDF4 }}}}
& 128 & 6.5438e-08& -- & 1.7232e-08& -- &0.00074557& -- &9.082e-05& -- \\
& 256 & 2.1846e-09& 4.9047& 4.6192e-10& 5.2213& 5.4526e-05& 3.7733& 5.8265e-06& 3.9623\\
& 512 & 3.0347e-11& 6.1697& 3.0551e-11& 3.9183& 3.6281e-06& 3.9097& 3.7459e-07& 3.9592\\
\hline
\hline
\parbox[t]{4.5mm}{\multirow{3}{*}{\rotatebox[]{90}{\parbox[c]{1.2cm}{\footnotesize IMEX-TVB(5,5) }}}}
& 128 & 8.0907e-08& -- & 1.4931e-07& -- &8.7961e-06& -- &3.1526e-05& -- \\
& 256 & 2.3712e-09& 5.0926& 4.7065e-09& 4.9875& 3.4175e-07& 4.6859& 2.7811e-07& 6.8247\\
& 512 & 1.4198e-10& 4.0619& 5.5157e-11& 6.415& 1.1008e-08& 4.9563& 9.1116e-09& 4.9318\\
\hline
\hline
\parbox[t]{4.5mm}{\multirow{3}{*}{\rotatebox[]{90}{\parbox[c]{1.2cm}{\footnotesize IMEX-BDF5 }}}}
& 128 & 4.6717e-08& -- & 1.0372e-07& -- &1.6077e-06& -- &9.3079e-06& -- \\
& 256 & 1.4517e-09& 5.0081& 3.1885e-09& 5.0237& 5.6036e-08& 4.8425& 1.6743e-07& 5.7968\\
& 512 & 8.9081e-12& 7.3484& 1.0165e-10& 4.9711& 1.4002e-09& 5.3226& 7.4625e-09& 4.4878\\
\hline
	\end{tabular}}
\end{table}

\begin{figure}\centering
\includegraphics[width=7.5cm]{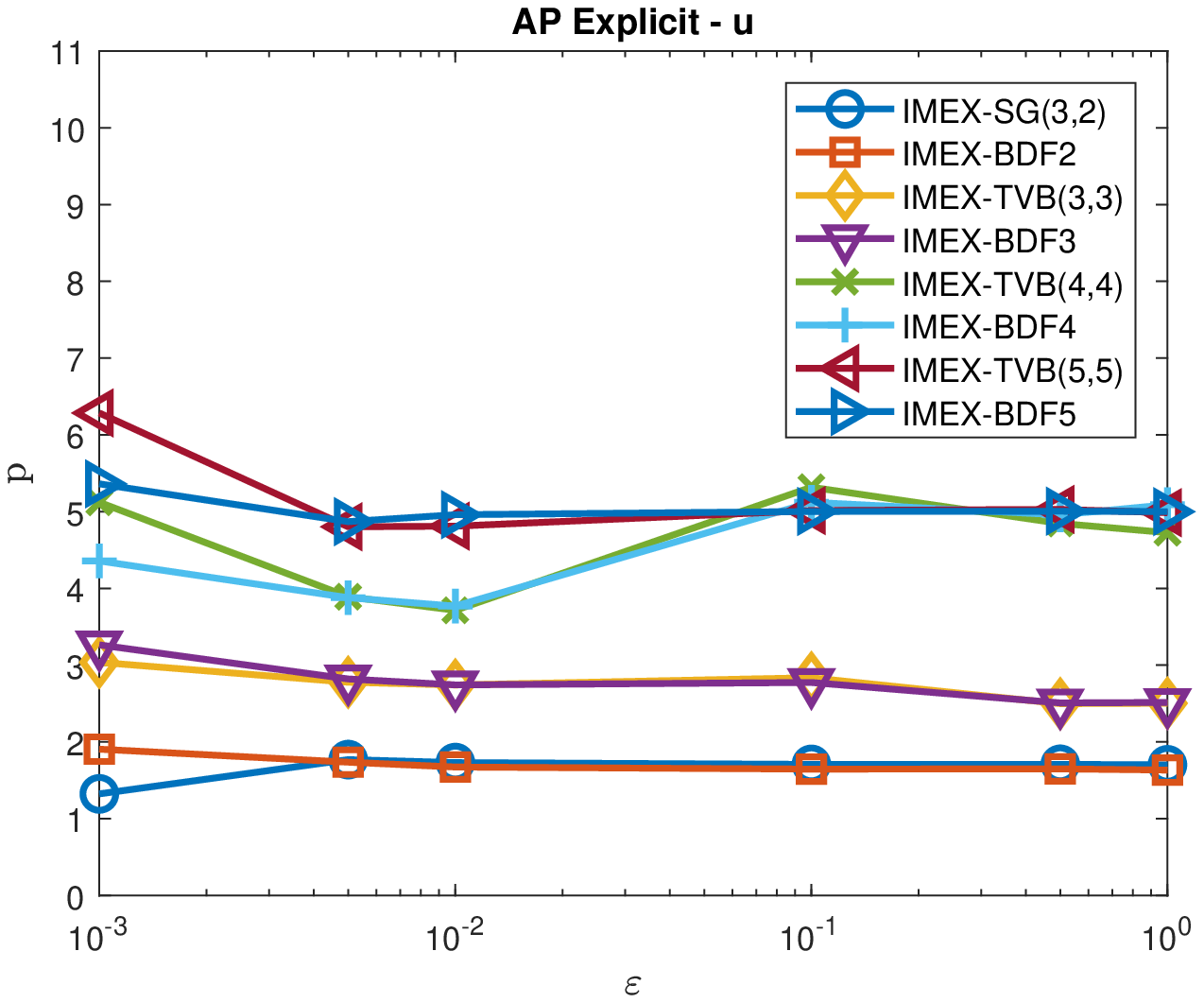}
\hspace{+0.5cm}
\includegraphics[width=7.5cm]{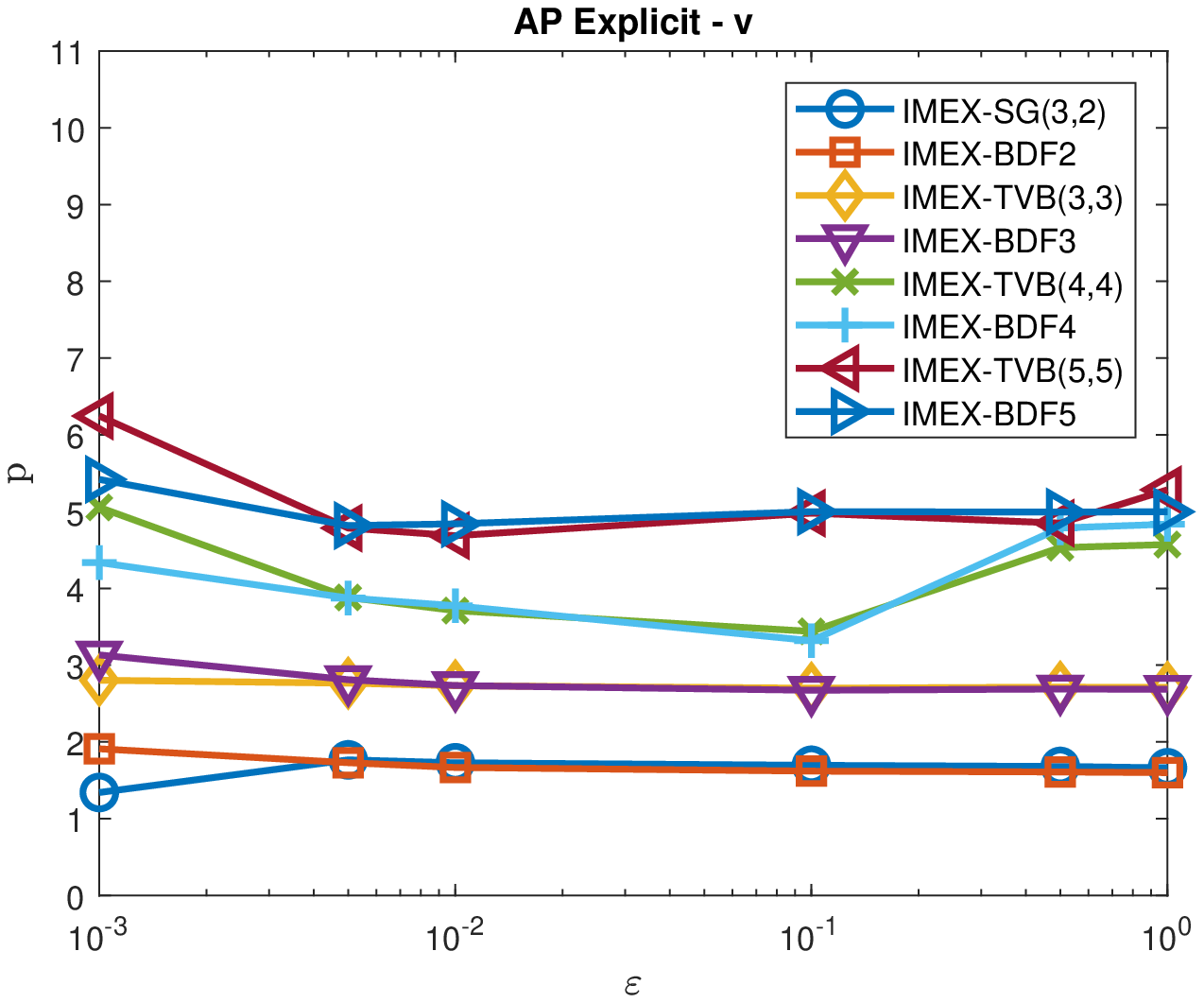}
\caption{Order of convergence of AP-explicit methods in the $L_1-$norm using \eqref{eq:conv_rate} for the $u$ variable (left) and the $v$ variable (right).}\label{fig:T1_order}
\end{figure}

\begin{figure}\centering
	\includegraphics[width=7.5cm]{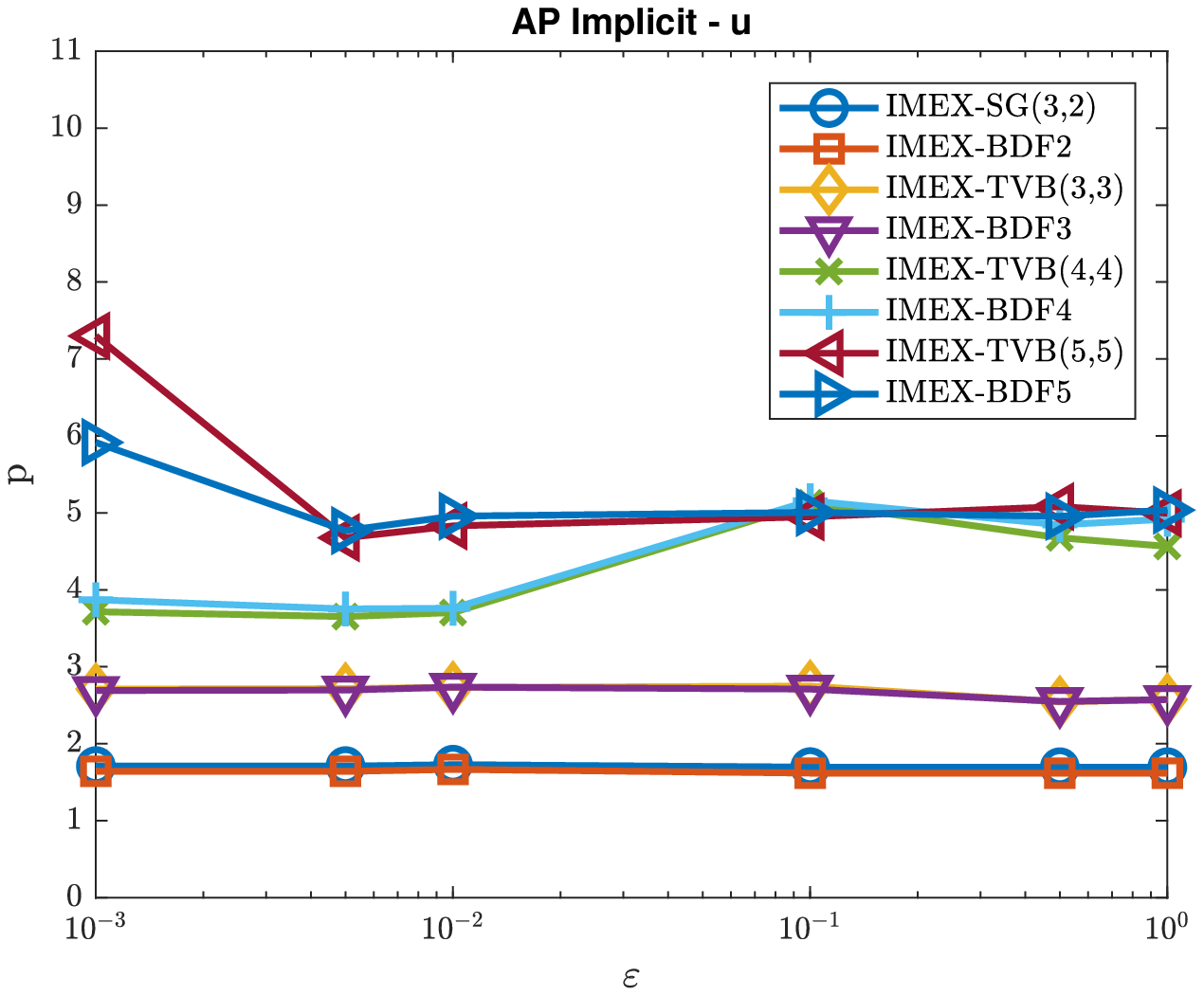}
	\hspace{+0.5cm}
	\includegraphics[width=7.5cm]{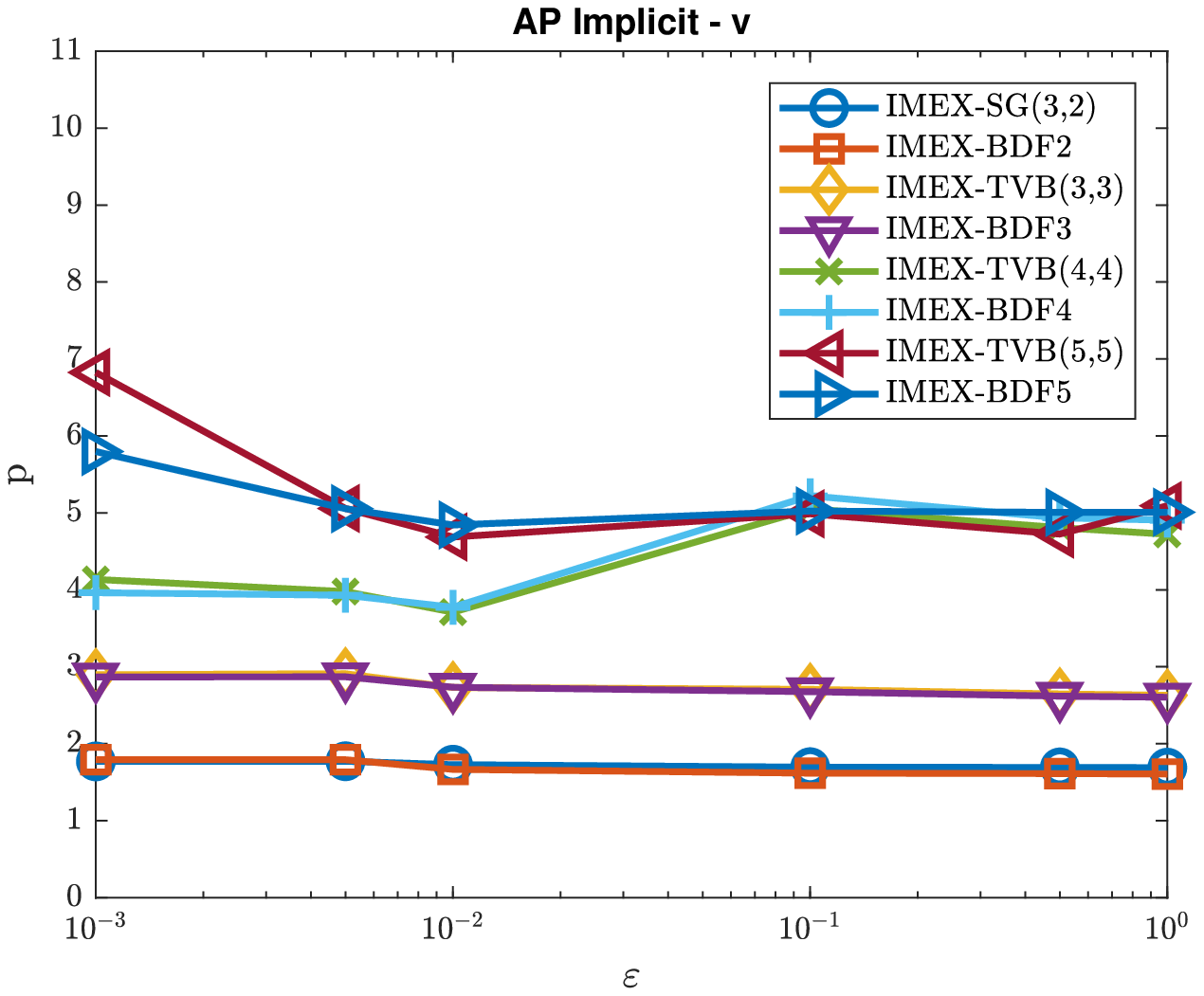}
	\caption{Order of convergence of AP-implicit methods in the $L_1-$norm using \eqref{eq:conv_rate} for the $u$ variable (left) and the $v$ variable (right).}\label{fig:T1_order_impl}
\end{figure}

\subsection{Test 2: Riemann problem for the linear model.}
Next, we consider a Riemann problem defined on the space interval $[0,4]$ with discontinuous initial data as follows 
\be
\left\{
\begin{array}{ll}
	u_L = 4.0,\quad v_L=0,\quad & 
	0\leq x \leq 2,
	\\
	u_R = 2.0,\quad v_R=0,\quad &
	2 < x \leq 4.
	\\
\end{array} \right. \label{RWdata}
\ee
\revised{For the above initial data, the linear hyperbolic system in the form \eqref{I61} with zero-flux boundary conditions is solved and comparisons with different values of the relaxation parameter $\varepsilon$ are shown. The same problem has been studied in \cite{BPR17} using IMEX Runge-Kutta schemes.}

In the limit $\varepsilon\to 0$, the exact solution of the corresponding advection-diffusion equation is known and it reads 
\be
u(x,t) = \frac{1}{2}(\rho_L + \rho_R) + \frac{1}{2}(\rho_L - \rho_R){\rm erf}\left(\frac{t-x+2}{2\sqrt{t}}\right),
\label{eq:exact_T2a}
\ee
with $\textrm{erf}(x)$ the error function. We report in Figure \ref{fig:T2_linRiem} the numerical solution for $u$ at final time $T = 0.25$ computed using two different time integration schemes, namely the IMEX-BDF2, and the IMEX-TVB(4,4). We choose $N=80$ points in space, and compare the AP-explicit approach  \eqref{eq:SPs_vec_ex2_alpha}, with the AP-implicit one \eqref{eq:SPs_vec_ex2_par}. The reference solution is computed with the IMEX-BDF5 scheme using $\Delta t_{\textrm{ref} }= \Delta t/10$ and $N=200$ space points. \revised{In order to preserve the CFL conditions the time steps for the different regimes of $\varepsilon$ are selected as follows: in the AP-explicit case $\Delta t = \lambda\Delta x\max\{\varepsilon, \Delta x\}$ and in the AP-implicit case $\Delta t = \lambda\Delta x\max\{\varepsilon, 1\}$ with $\lambda=0.4$.}
%
\revised{In Figure \ref{fig:T2_linRiem}, for the diffusive limit, we observe that the different schemes agree well with the exact solution \eqref{eq:exact_T2a}. In the hyperbolic regime, $\varepsilon=0.5$, the shock is correctly captured compared to the reference solution by both schemes with a slightly better resolution in the case of IMEX-TVB(4,4) for the AP-implicit approach. Note that, no spurious oscillations are observed for IMEX-BDF2 scheme, even if it does not satisfy any specific TVB stability property in the hyperbolic regime.} 

\begin{figure}\centering
	{\includegraphics[width=7.5cm]{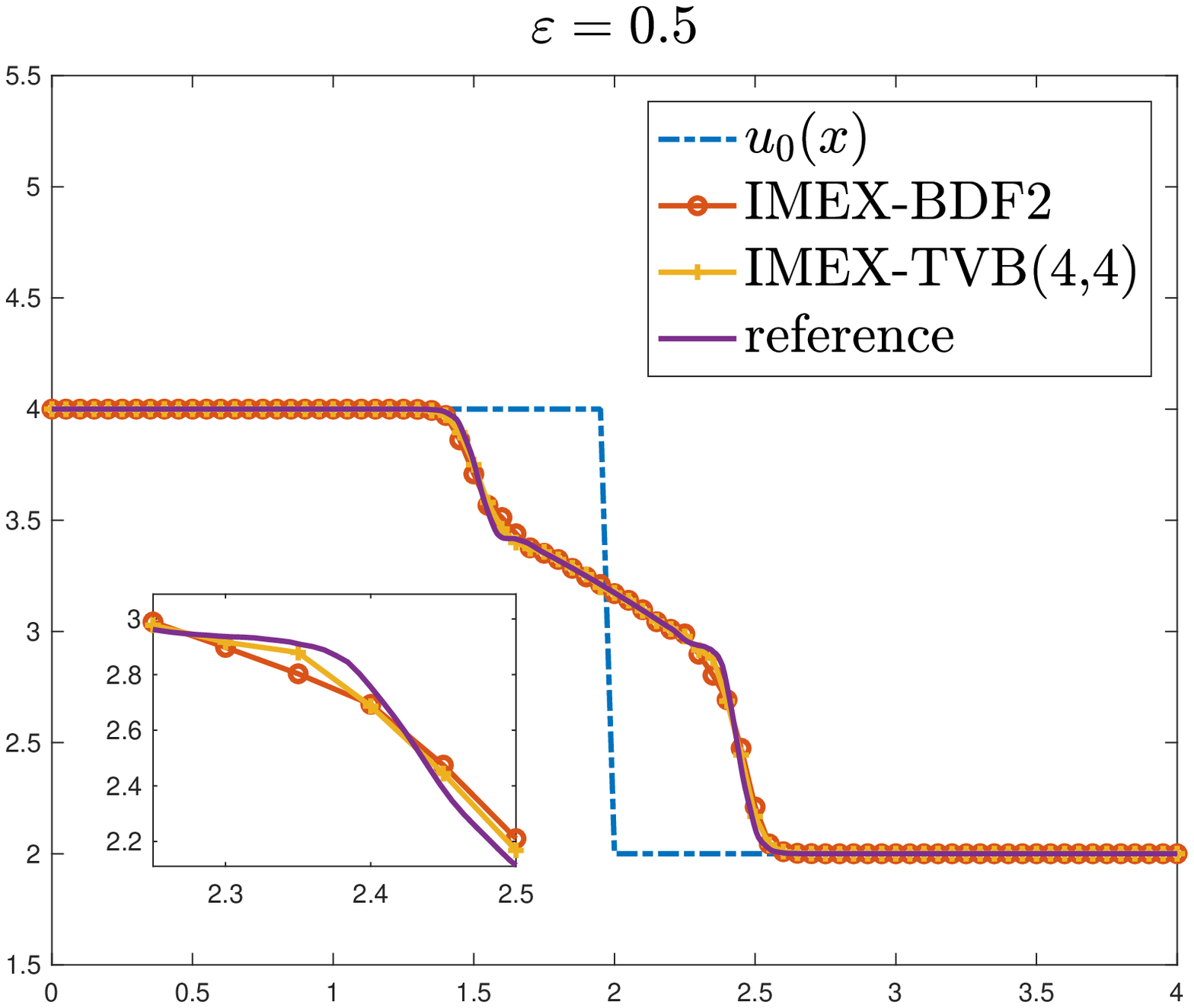}}
	\hspace{+0.5cm}
	{\includegraphics[width=7.5cm]{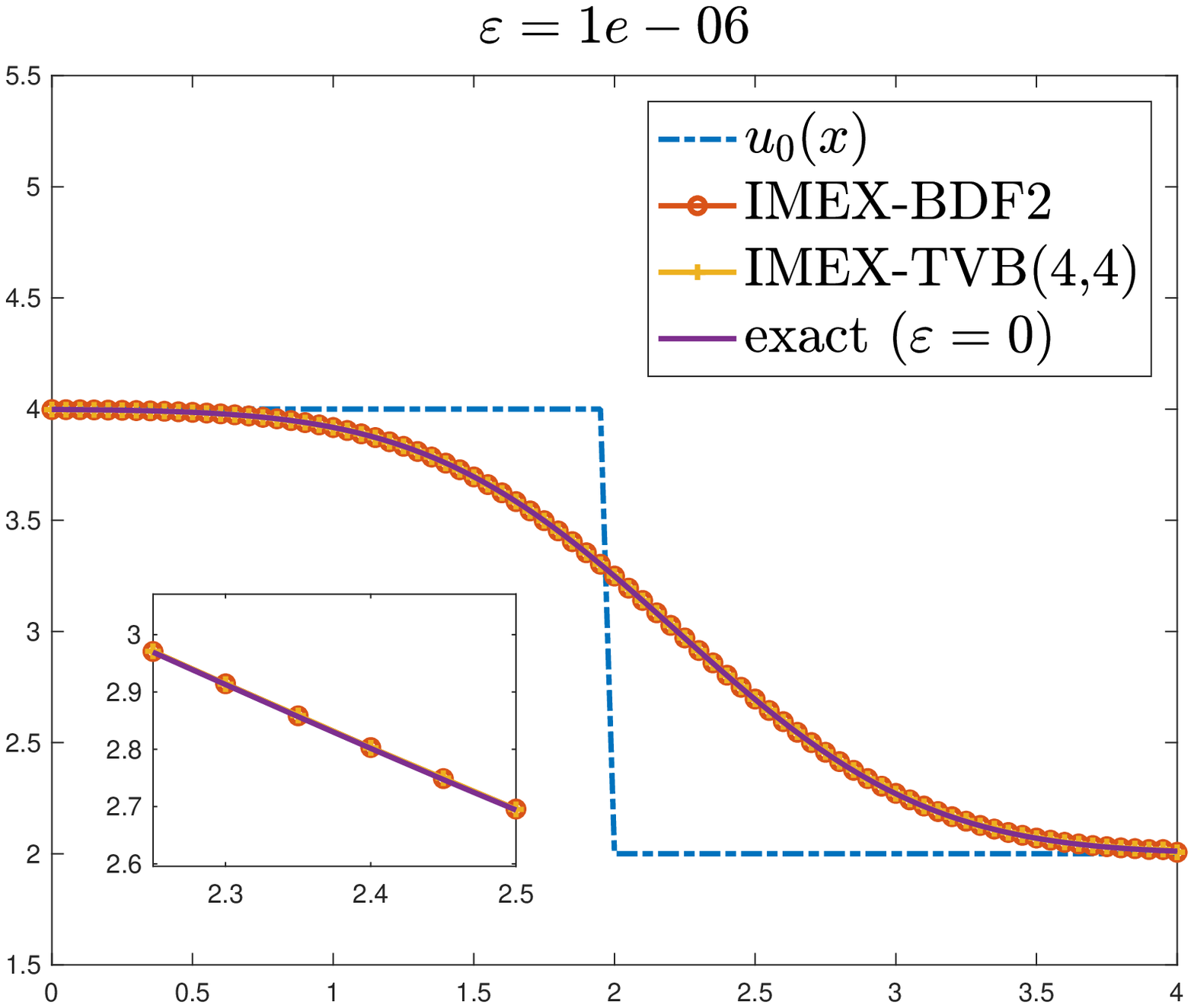}}
	\\
	{\includegraphics[width=7.5cm]{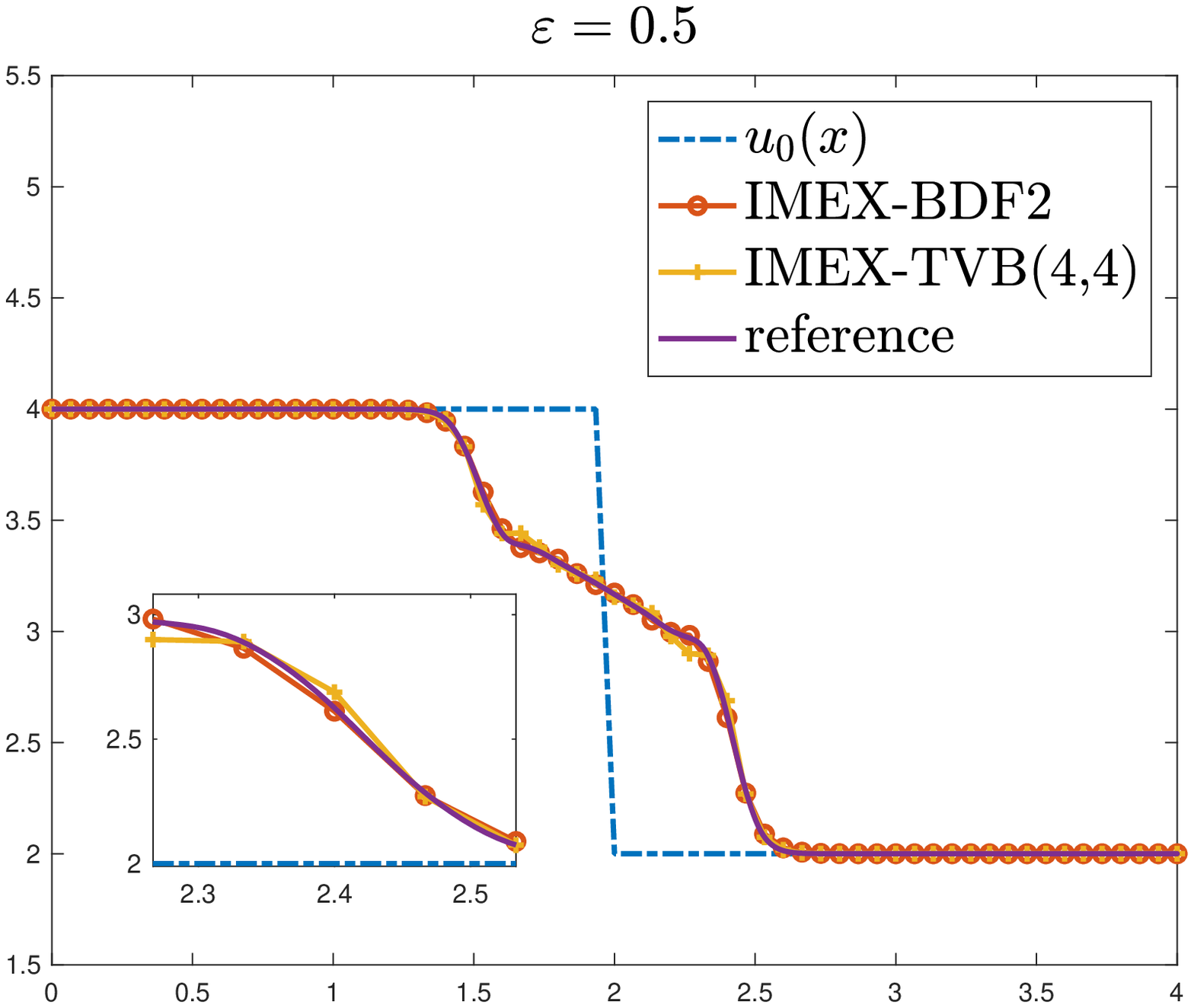}}
	\hspace{+0.5cm}
	{\includegraphics[width=7.5cm]{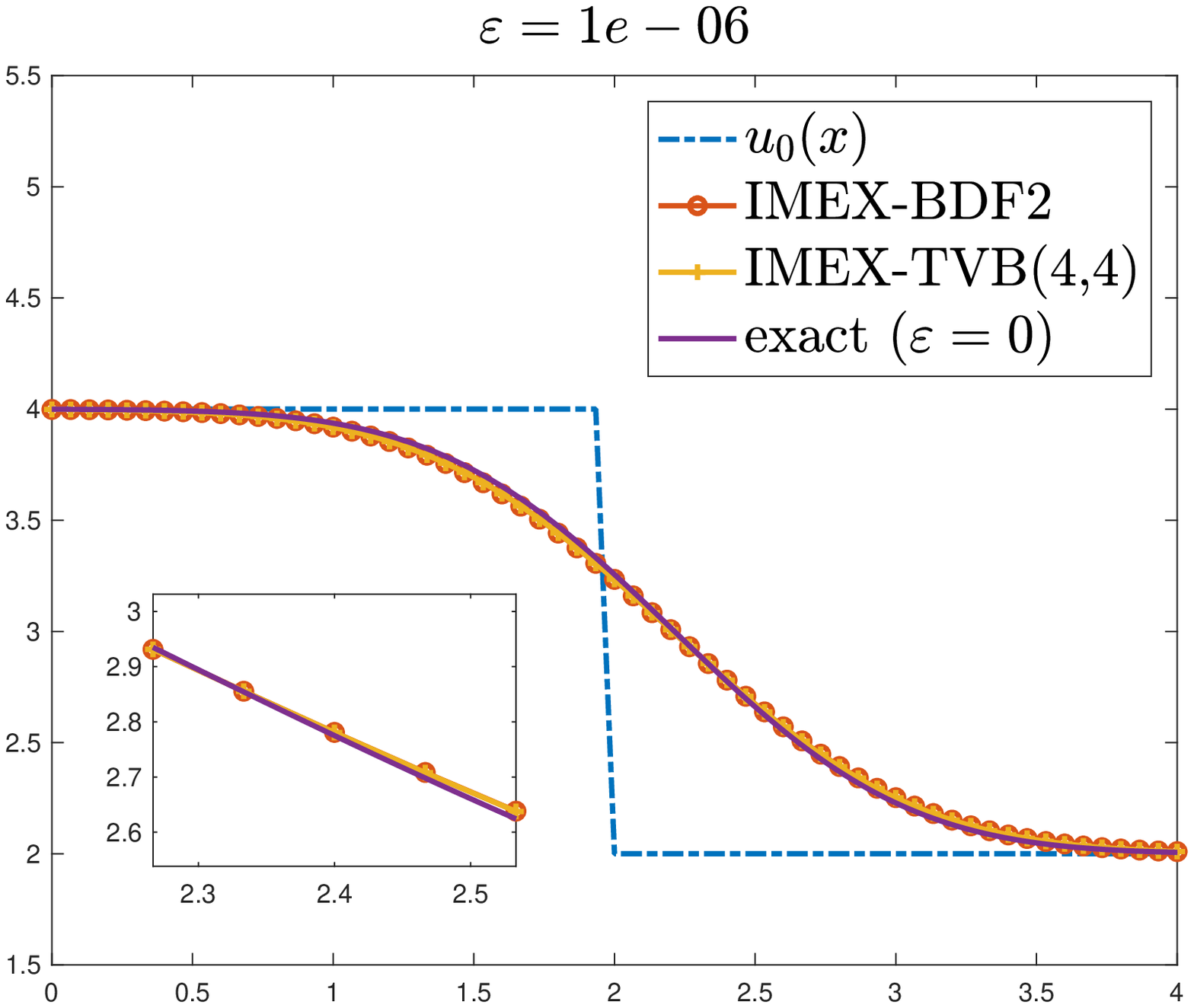}}
	\caption{Test 2. Riemann problem for linear system with diffusive scaling \eqref{T1} at final time $T = 0.25$. Hyperbolic regime with $\varepsilon = 0.5$ (left) and diffusive regime with $\varepsilon = 10^{-6}$ (right). Top row AP-explicit approach, bottom row AP-implicit approach. 
	}\label{fig:T2_linRiem}
\end{figure}

\subsection{Test 3: Barenblatt solution for the porous media equation.}
We then consider a nonlinear diffusion limit, by considering the following hyperbolic system with diffusive relaxation
\be
\begin{cases}
	\displaystyle  
	\partial_t\rho + \partial_x \j =0, \\
	\displaystyle    
	\partial_t \j + \frac{1}{\epsi^2} \partial_x \rho = -\frac{1}{\epsi^{2+\alpha}} k(\rho) \j.
\end{cases}
\label{E1}
\ee
\revised{This problem has been previously studied in \cite{JPT, NP2}.} Note that, when $k(\rho)=1$, the system \eqref{E1} is a model of relaxing heat flow and, as $\epsi \to 0$, it relaxes towards the heat equation
\be  
\partial_t\rho = \epsi^\alpha \partial_{xx} \rho,
\qquad 
\j = -\partial_x \rho.
\label{E2}
\ee
On the other hand, by choosing $k(\rho)=(2\rho)^{-1}$, the limiting equation for this model results in the porous media equation
\be
\displaystyle  
\partial_t\rho = \epsi^\alpha \partial_{xx} \rho^2, \qquad
\j = - 2 \rho \partial_x \rho.
\label{porous1}
\ee
\revised{For this problem we consider the IMEX-BDF2 and the IMEX-BDF5 schemes for $\alpha=0$ by computing the numerical solution 
 with $N=80$ mesh points, using the AP-explicit approach and selecting the time step as $\Delta t = \lambda\Delta x \max\{\varepsilon,\Delta x\}$ with $\lambda=0.4$.}
 
The numerical solution in the limit $\varepsilon\to 0$ is compared with the analytical Barenblatt solution for the porous media equation \cite{Bar1} 
\be
\rho (x,t)= \left\{
\begin{array}{ll}
	\displaystyle
	\frac{1}{r(t)} \left[ 1-\left( \frac{x}{r(t)} \right)^2 \right], &
	|x|\leq r(t),
	\\
	\displaystyle
	0, & |x|>r(t),
	\\
\end{array}
\right.
\label{bare}
\ee
where $r(t)=[12(t+1)]^{1/3},~~t \geq 0$ and $x \in [-10,10]$. Note that the above solution defines also the initial state of the system where in addition we considered $v(x,t=0)=0$. In Figure \ref{fig:T2_exactBarenblatt}, this comparison is shown 
using $\epsi = 10^{-6}$ at time $T = 3$. \revised{The sharp front of the Barenblatt solution is very well captured by both schemes without observing any relevant different by the increase order of accuracy of the IMEX-BDF5 method. This is somewhat expected, since the CFL condition $\Delta t=O(\Delta x^2)$, in practice, gives to the IMEX-BDF2 scheme the same accuracy of a fourth order method. Therefore, the accuracy barrier here is represented by the $5$-th order space discretization method.}


\begin{figure}\centering
	\includegraphics[width=8.00cm]{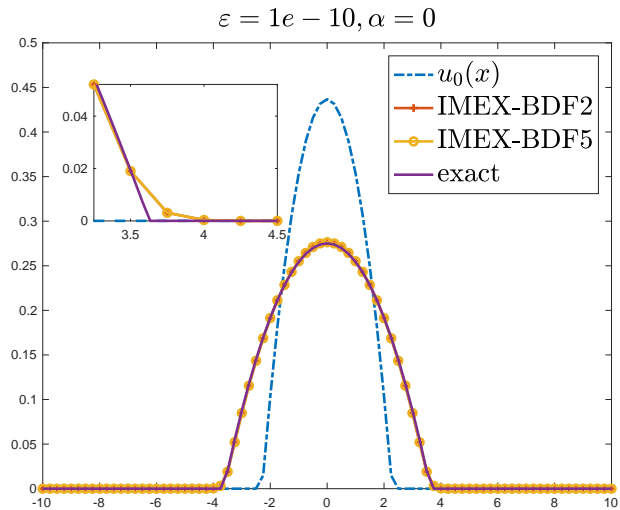}
	\caption{Test 3. Comparison between the Barenblatt analytical solution of system (\ref{porous1}), and the numerical solution obtained from system (\ref{E1}) with $k(u)=(2u)^{-1}$ and $\epsi = 10^{-6}$ at time $T = 3$.}\label{fig:T2_exactBarenblatt}
\end{figure}


\subsection{\revised{Test 4: Applications to the Ruijgrook--Wu model.}} 
\revised{We consider, in this last test case, an application of the schemes to the so-called Ruijgrook--Wu model of rarefied gas dynamic \cite{RW, GP, GPR, NP2, BPR17, JPT}
\be
\left\{
\begin{array}{l} 
\displaystyle  
M \partial_t f^+ + \partial_x f^+=-\frac{1}{Kn} ( a f^+ - b f^- - c f^+ f^-),\\[+.2cm]
\displaystyle
M \partial_t f^- - \partial_x f^-=\frac{1}{Kn} ( a f^+ - b f^- - c f^+ f^- ), 
\\
\end{array}
\right. 
\label{E3}
\ee
where $f^+$ and $f^-$ denote the particle density distribution at time $t$, position $x$ and with velocity $+1$ and $-1$ respectively. Here $Kn$ is the Knudsen number, $M$ is the Mach number of the system and $a$,$b$ and $c$ are positive constants which characterize the microscopic interactions . The local (Maxwellian) equilibrium  $f_{\infty}^{\pm}$ is characterized by 
\be
f^{+}_{\infty} = \frac{bf^{-}_{\infty}}{a-cf^{-}_{\infty}}.
\ee
The macroscopic variables for the model are the density $\rho$
and momentum $v$ defined by
\be
\rho = f^++f^-,\quad v=(f^+ -f^-)/M.
\ee
The nondimensional multiscale problem is obtained taking $M = \epsi^\alpha$ and $Kn = \epsi$, the Reynolds number of the system is then defined as usual according to $Re = M/Kn = 1/\epsi^{1-\alpha}$. In macroscopic variables taking $a=b=1/2$, $c=M=\epsi^{\alpha}$ this can be written as \cite{GPR}
\be
\begin{cases}
	\displaystyle  
	\partial_t \rho + \partial_x \j =0 \,, \\
	\displaystyle
	\partial_t \j + \frac{1}{\epsi^{2\alpha}}\partial_x \rho
	= -{{1}\over{\epsi^{1+\alpha}}} \left[\j - \frac{1}{2}\left( \rho^2 - \epsi^{2\alpha} \j^2\right) \right].
\end{cases}
\label{RW}
\ee
The model has the nice feature to provide nontrivial limit behaviors for several values of $\alpha$ including the corresponding compressible Euler ($\alpha=0$) limit and the incompressible Euler ($\alpha\in (0,1)$) and Navier-Stokes ($\alpha=1$) limits (see \cite{GP, GPR}).
}
\revised{In \cite{NP2} the above model has been used under a similar but different scaling. If we denote with the pair $\talpha,\tepsi$ the scaling parameters in \cite{NP2}, these correspond to take $M=\tepsi$, $Kn=\tepsi^{1+\talpha}$. The two scaling can be made equivalent for $\alpha\neq 0$, if we map the pair $\alpha,\epsi$ into $\talpha,\tepsi$ as follows
\be
\talpha = \frac{1-\alpha}{\alpha},\qquad
\tepsi = \e^{\alpha},\qquad \alpha \neq 0.
\label{eq:map}
\ee
Note, that the non-linearity on the source term depends both on $\rho$ and $\j$, so that we have $f=f(u,v)$ in \eqref{I61}. Nevertheless, following the same strategy described in the previous Sections, our methods can be applied in a straightforward way also in this situation. 
}

\revised{For the Ruijgrook--Wu model \eqref{RW} it can be shown, via Chapman-Enskog expansion, that for $\alpha \in (1/3,1]$ so that $2\alpha > 1-\alpha$ and small values of $\epsi$ we have
\be
\displaystyle \j = \frac{1}{2} \rho^2 -\varepsilon^{1-\alpha} \partial_x \rho + {\mathcal O}(\epsi^{2\alpha}).
\label{eq:vbur}
\ee
Then, the solution behavior is characterized by the viscous Burgers equation
\be
\begin{array}{l} 
	\displaystyle \partial_t \rho + \partial_x \left(\frac{\rho^2}{2}\right) = \varepsilon^{1-\alpha} \partial_{xx}\rho+ {\mathcal O}(\epsi^{2\alpha}).
\end{array}
\label{eq:burger}
\ee
}
\comment{
In the sequel, for {\em Test 4a} we apply the AP-explicit approach with $\Delta t = 
\lambda\Delta x \max\{\epsi, \min\{1,\Delta x/\epsi^{1-\alpha}\}\}$, $\lambda=0.4$ which again corresponds to the maximum between the CFL of the hyperbolic part and the one originated by the limiting viscous Burger equation, whereas for {\em Test 4b} and {\em Test 4c} we use the AP-implicit approach with $\Delta t = \lambda\Delta x \max\{\epsi, 1\}$ and $\lambda=0.1$. Therefore, for $\epsi \leq 1$ the AP-implicit approach uses the same CFL condition in all test cases independently of $\epsi$ and $\alpha$.
The reference solution is computed using the IMEX-BDF5 scheme with time step $\Delta t_{\textrm{ref} }= \Delta t/10$ and $N_{\textrm{ref} }=400$ space points. The initial data for the non conserved quantity $v$ has been taken well prepared accordingly to \eqref{eq:vbur}.
In the sequel we will restrict to $\talpha \in [0,1]$ as in \cite{NP2}, therefore, from \eqref{eq:map}, we have $\alpha\in [0.5,1]$. Note, however, that the schemes can be applied for any value of $\alpha > 0$. We refer to \cite{GPR} for results in the fluid limit $\alpha=0$.}

\subsubsection*{\revised{Test 4a. Riemann problem}}
We select the Riemann problem characterized by \eqref{RWdata} as initial data. The numerical parameters are also the same of the previous test case, therefore we fix the final time $T = 0.25$, and we compute the solution $u(x)$ using two different time AP-explicit integration schemes, BDF2 and TVB(4,4) with $N=100$ points in space.}

\begin{figure}\centering
	{\includegraphics[width=7.5cm]{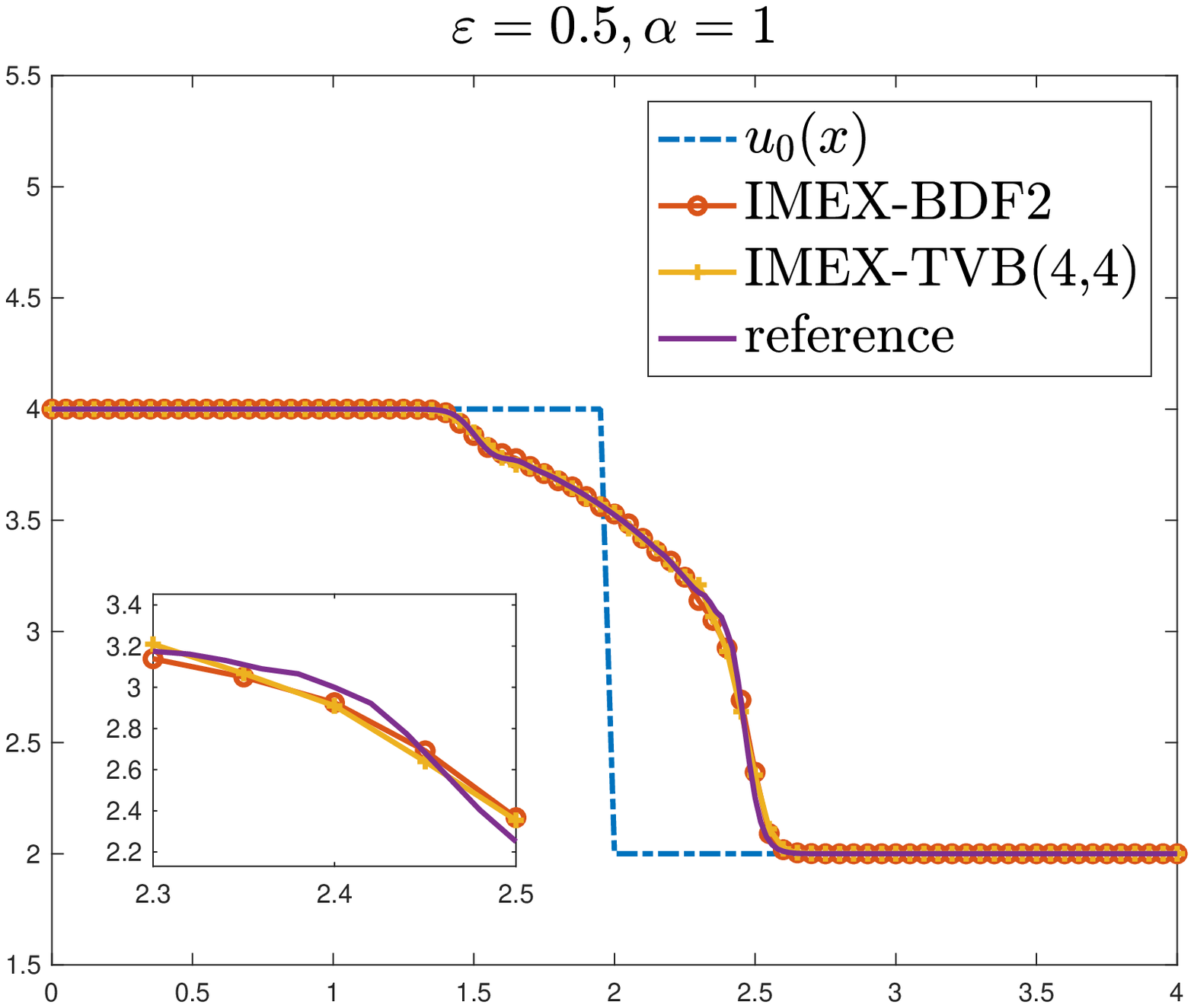}}
	\hspace{+0.5cm}
	{\includegraphics[width=7.5cm]{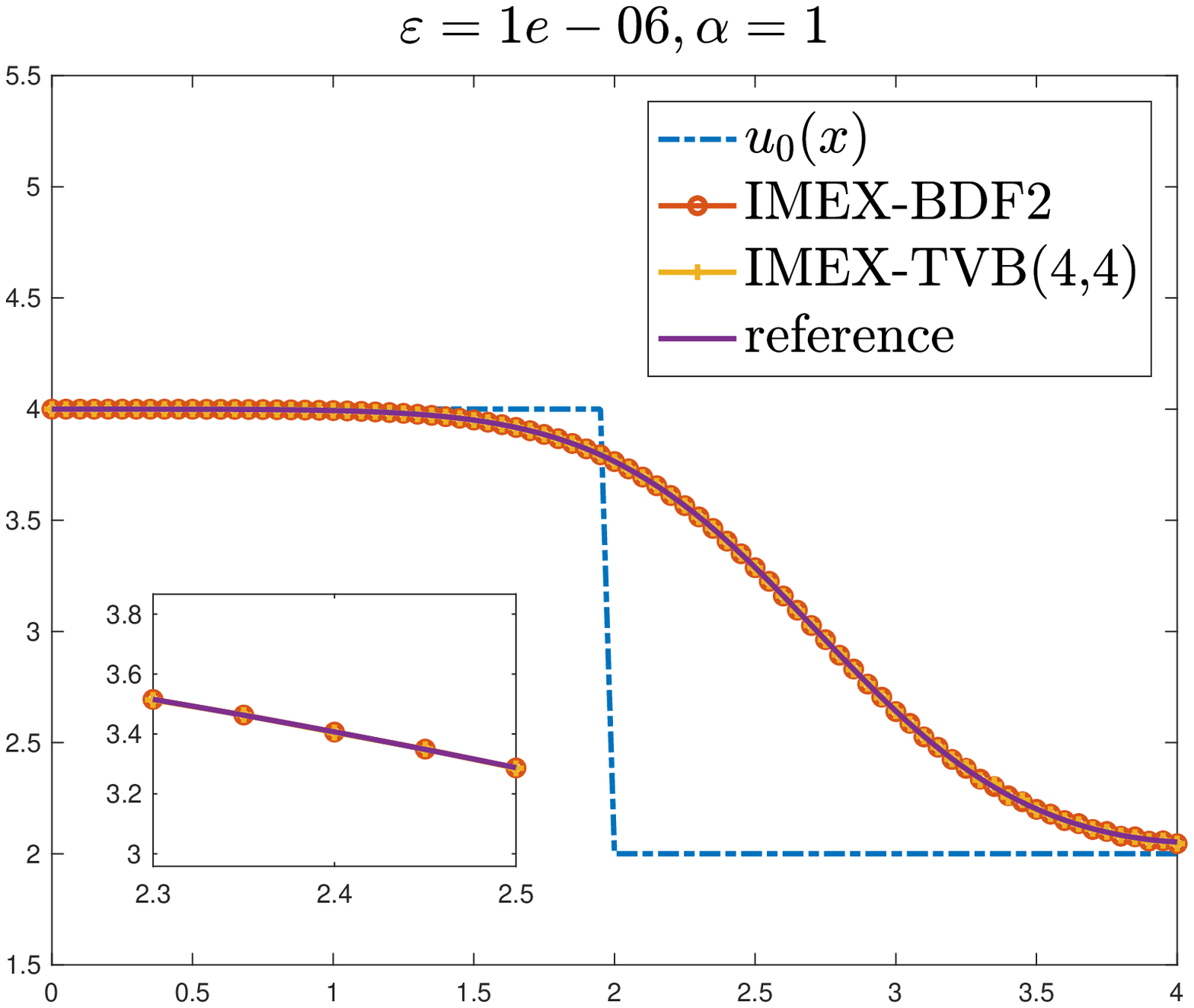}}
	\\\vspace{+0.25cm}
	{\includegraphics[width=7.5cm]{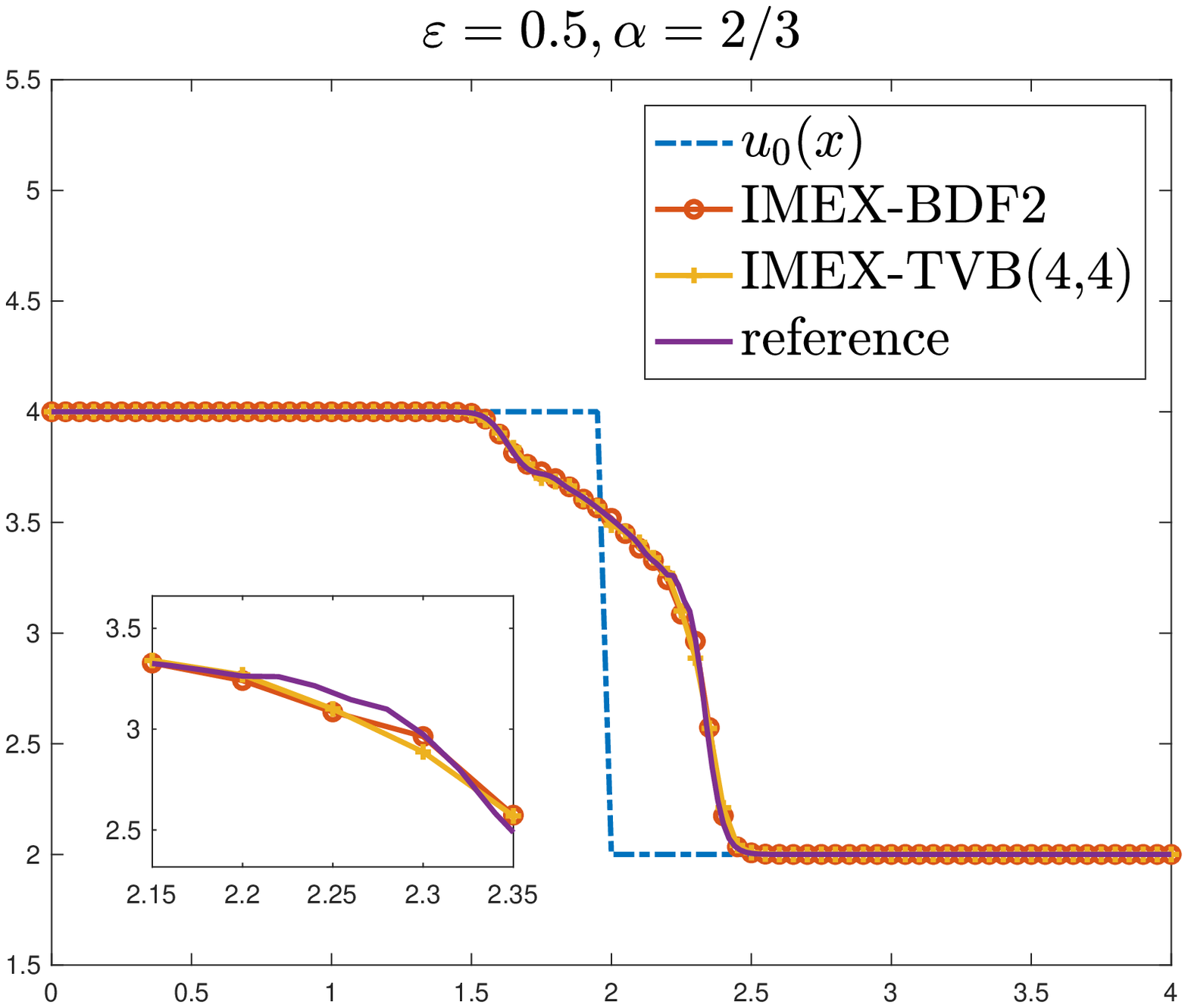}}
	\hspace{+0.5cm}
	{\includegraphics[width=7.5cm]{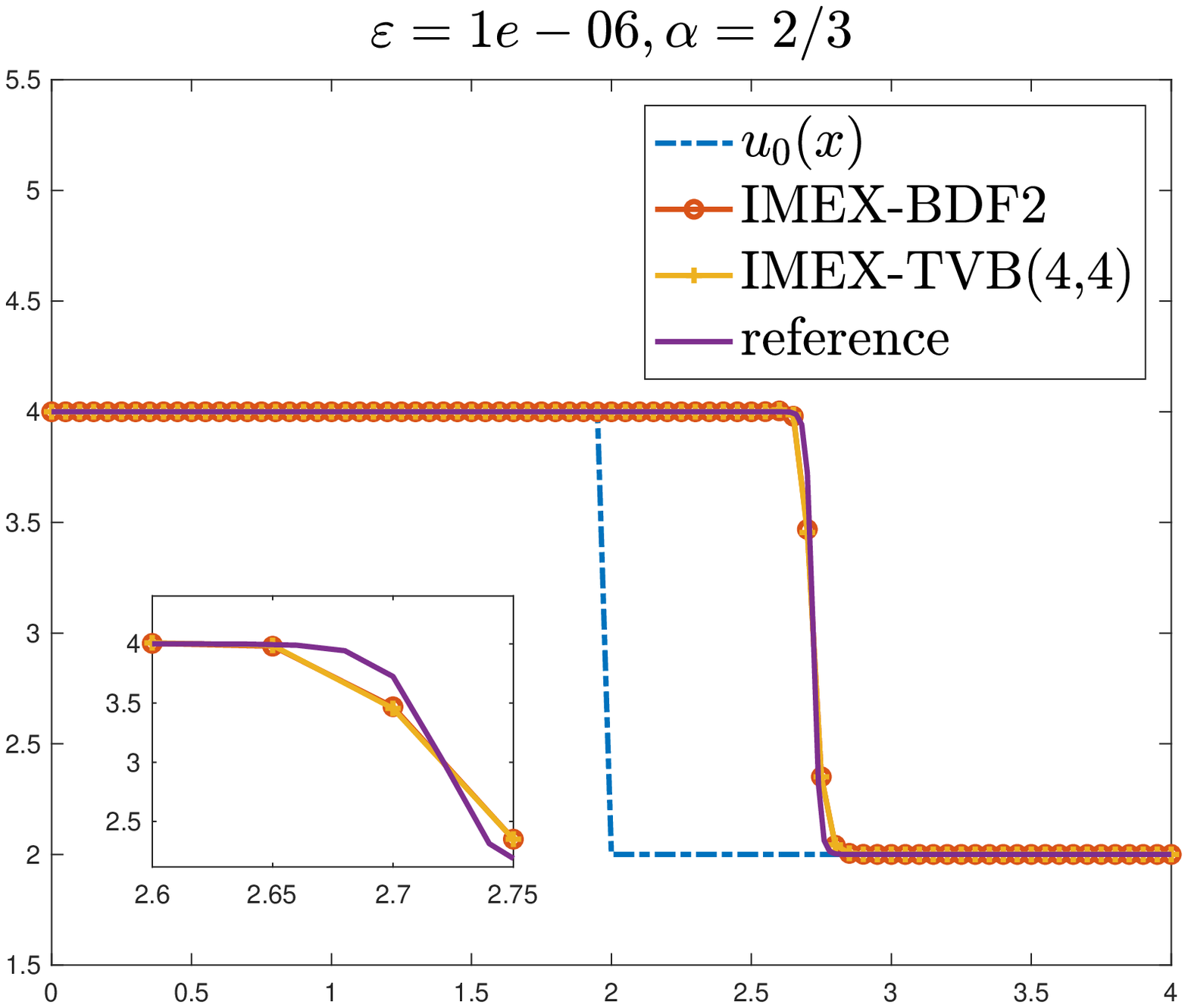}}
	\caption{{\em Test 4a}. Riemann problem for the non linear Ruijgrook--Wu model \eqref{RW} at final time $T = 0.25$ with AP-explicit methods. Rarefied regime with $\varepsilon = 0.5$ (left) and limiting behavior for $\varepsilon = 10^{-6}$ (right). Top row $\alpha = 1$, bottom row $\alpha =2/3$.}\label{fig:T2c_nonlinRiem}
\end{figure}

\revised{In Figure \ref{fig:T2c_nonlinRiem} we report in the left column the behavior in the rarefied (non stiff) regime $\varepsilon = 0.5$ and in the right-column the limit behavior $\varepsilon = 10^{-6}$, whereas top row depicts the behavior for $\alpha=1$ and bottom row $\alpha = 2/3$. Even if the rarefied solutions for the different values of $\alpha$ are similar, we remark the difference of the solution profiles in the limit $\epsi\to 0$. For $\alpha=1$ we obtain the classical viscous Burger equation where the discontinuity is smeared out by the diffusion term, while for $\alpha=2/3$ we have the sharp shock formation of the inviscid Burger equation. In both cases the methods yield an accurate description of the dynamic without spurious oscillations or excessive numerical dissipation.} 

\subsubsection*{\revised{Test 4b. Propagation of a square wave}}
\begin{figure}\centering
	{\includegraphics[width=4.75cm]{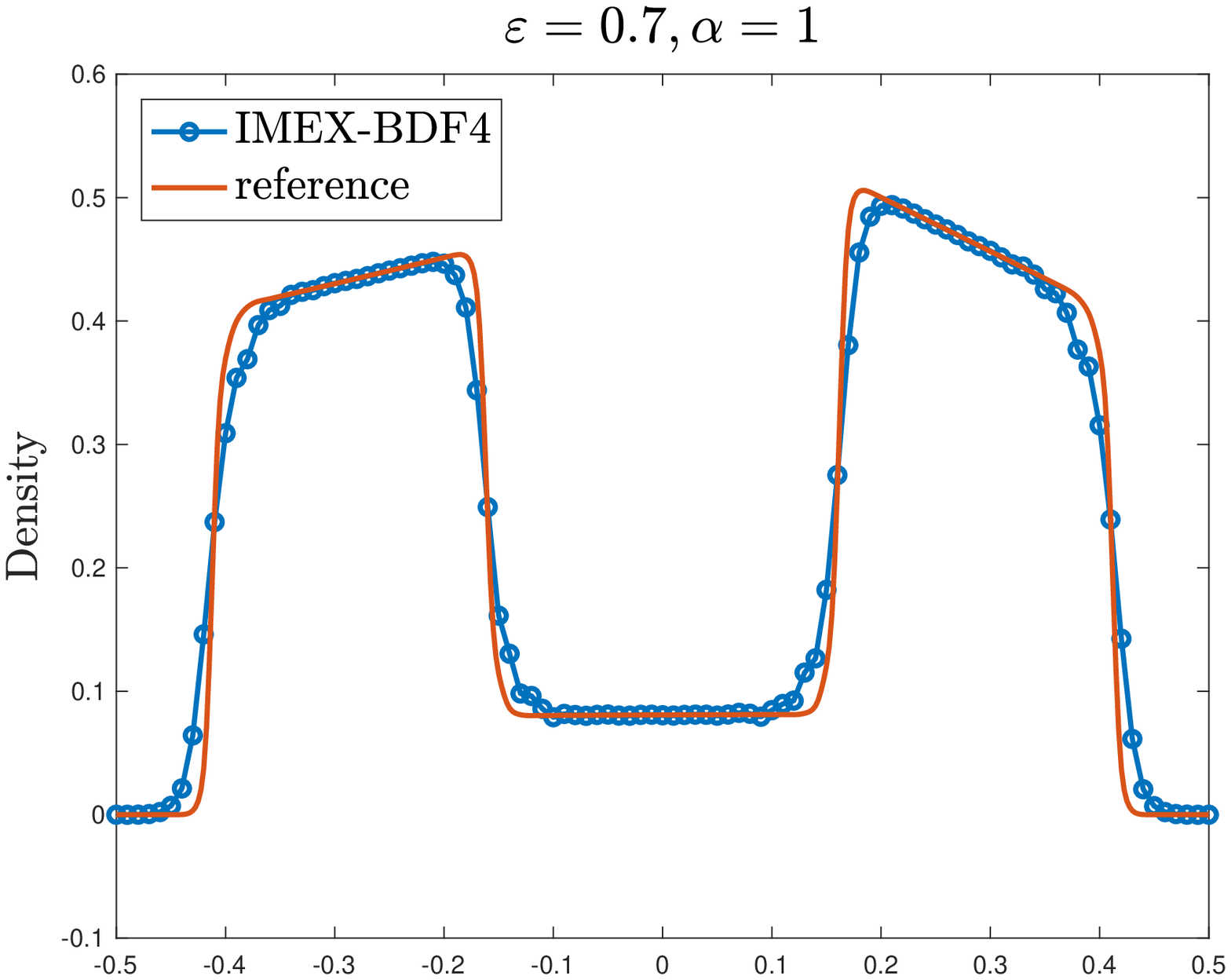}}
	\hspace{+0.25cm}
	{\includegraphics[width=4.75cm]{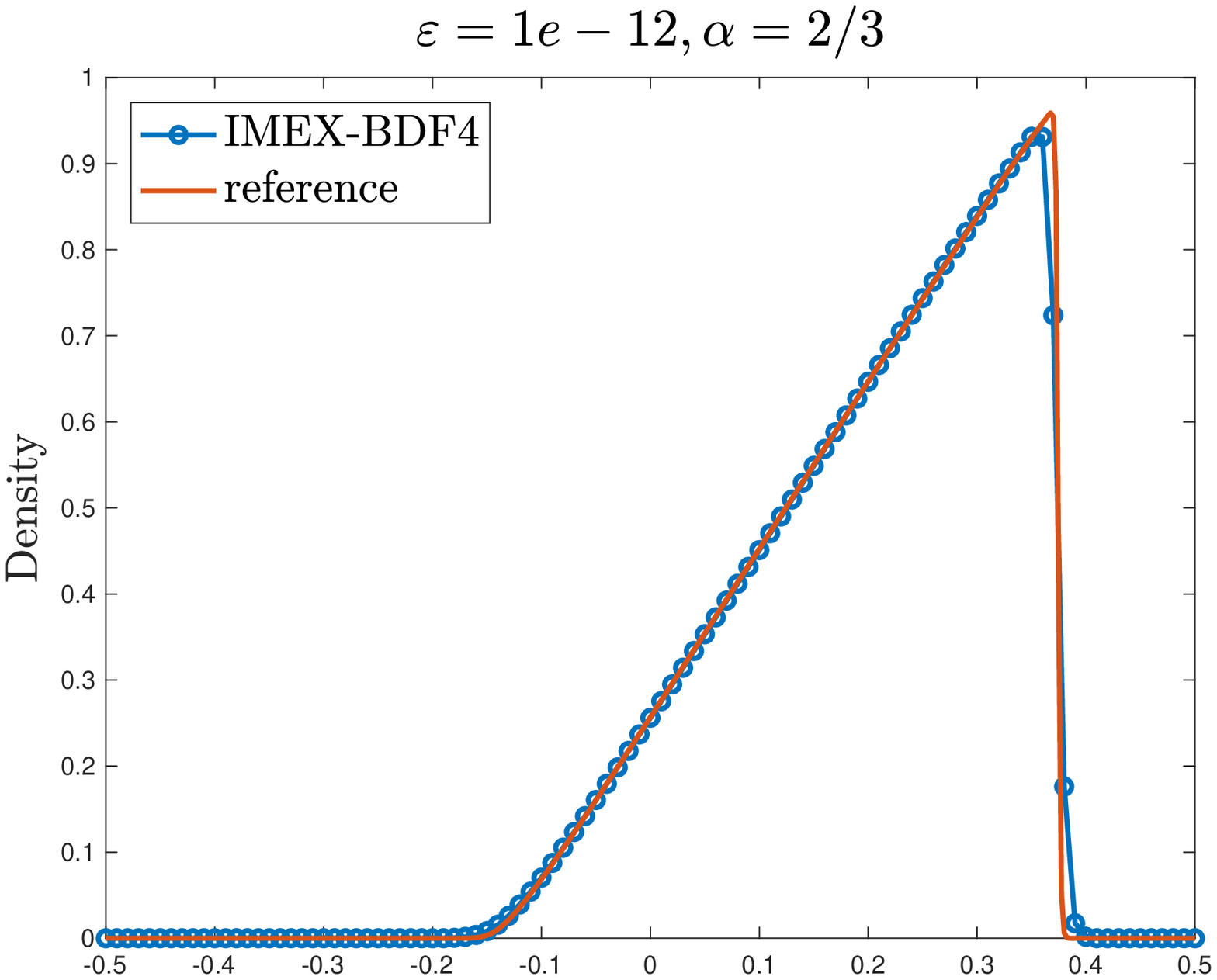}}
	\hspace{+0.25cm}
	{\includegraphics[width=4.75cm]{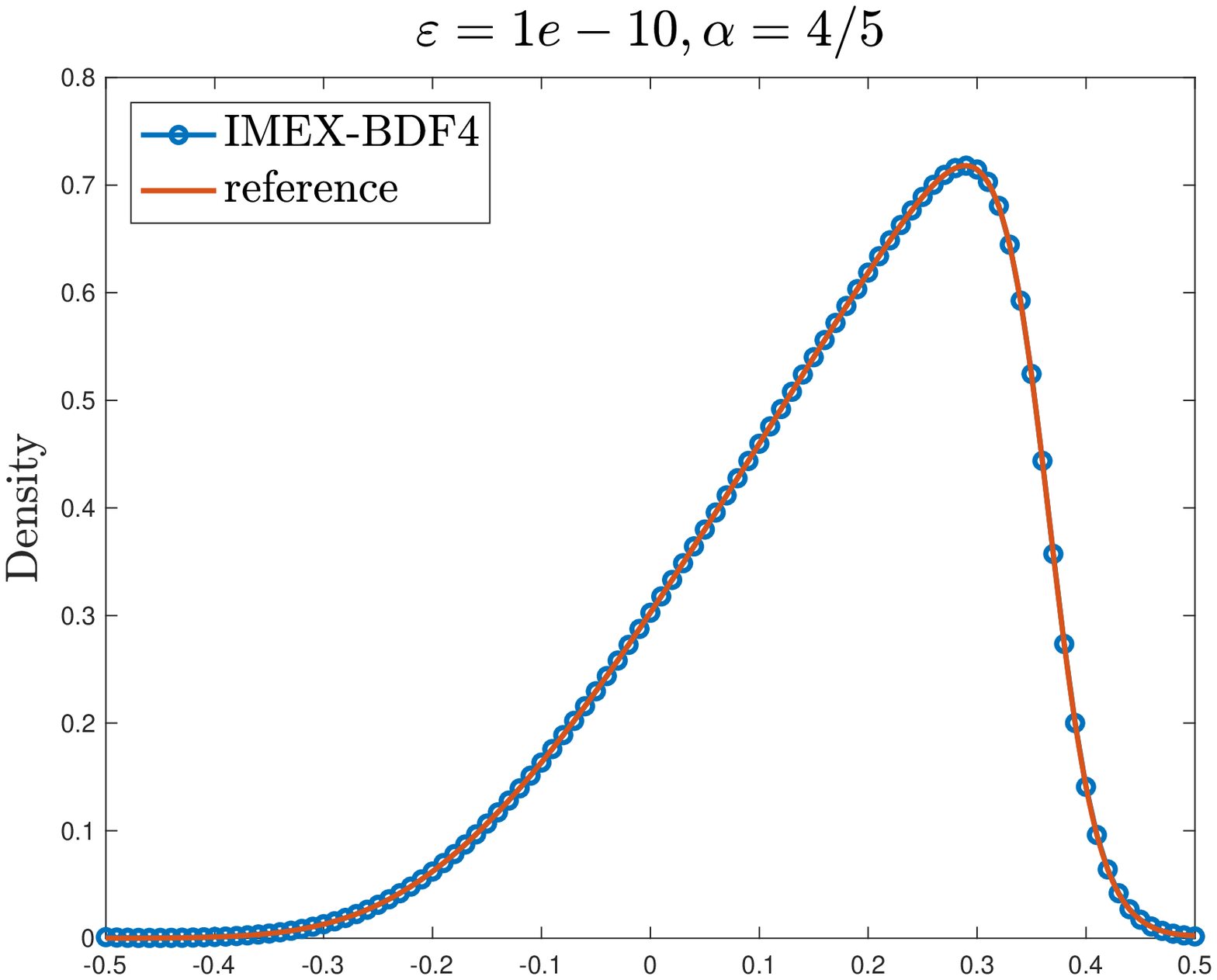}}
	\\\vspace{+0.25cm}
	{\includegraphics[width=4.75cm]{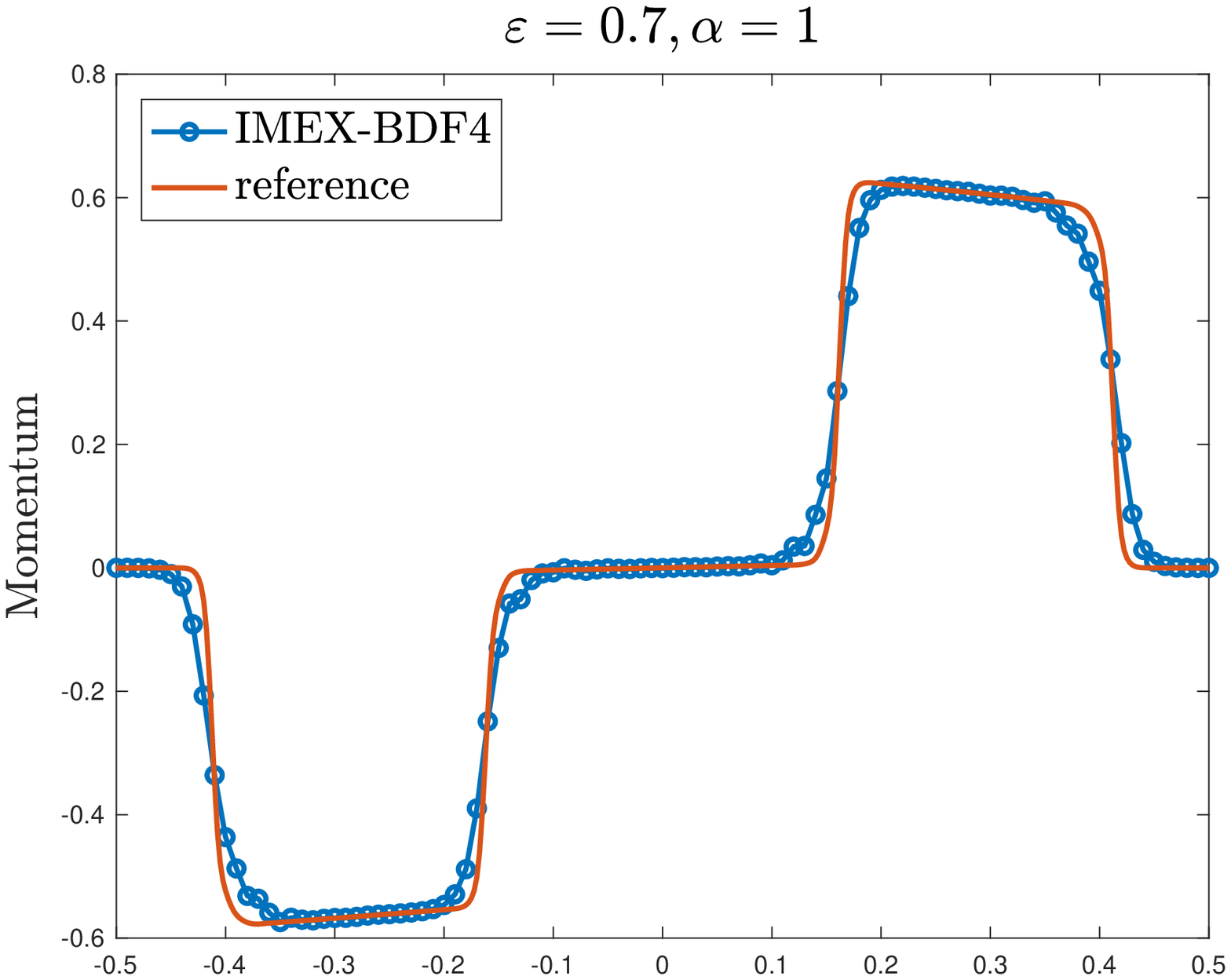}}
	\hspace{+0.25cm}
	{\includegraphics[width=4.75cm]{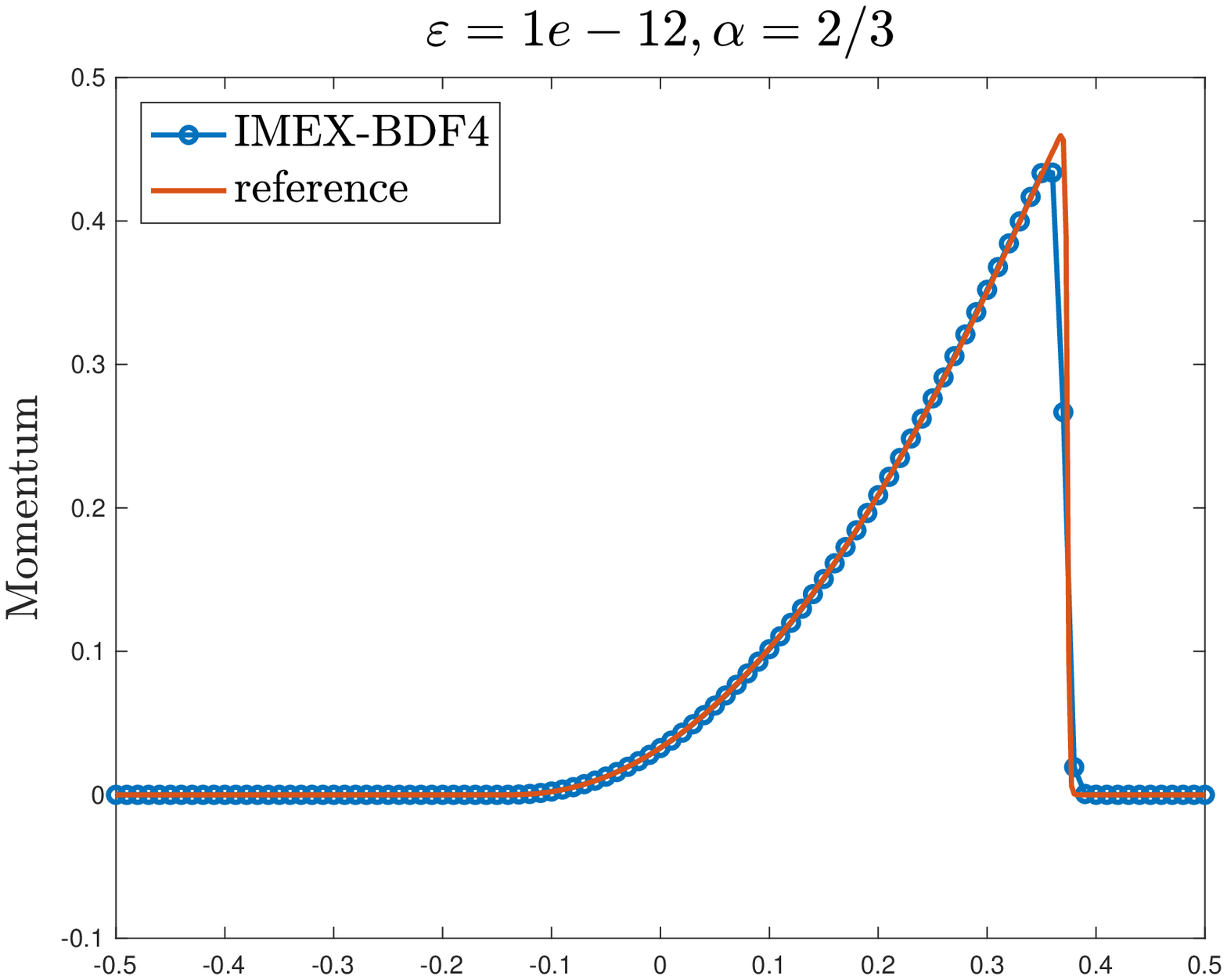}}
	\hspace{+0.25cm}
	{\includegraphics[width=4.75cm]{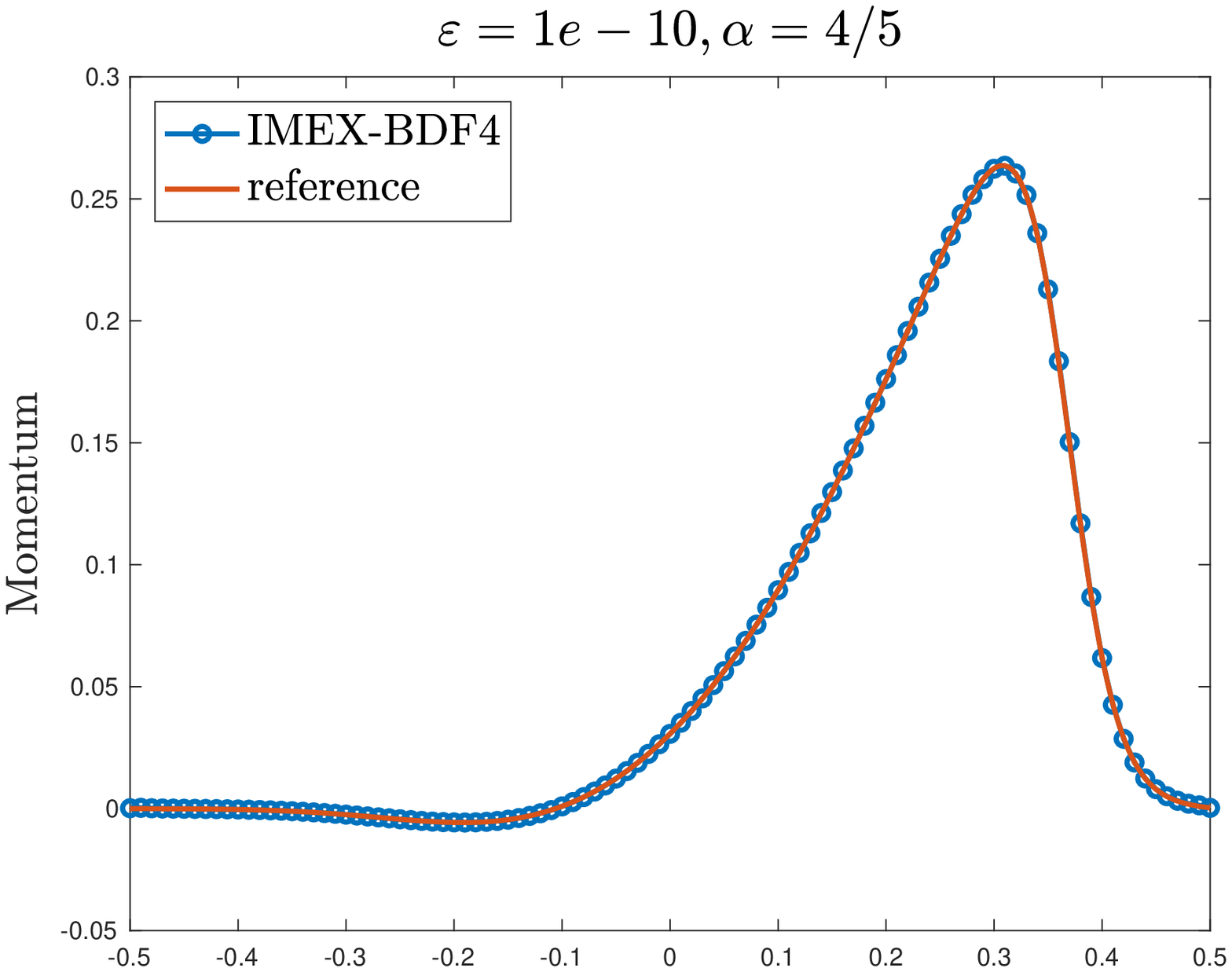}}
	\caption{{\em Test 4b}. Non linear Ruijgrook--Wu model \eqref{RW} with discontinuous initial data \eqref{RW0}, top column reports the density, bottom column the momentum. Left row $\epsi=0.7, \alpha=1$ and final time $T=0.2$, middle row $\epsi=10^{-12},\alpha=2/3$ final time $T=0.5$, right row $\epsi = 10^{-10},\alpha=4/5$, with final time $T=0.5$. The AP-implicit formulation has been used.}\label{fig:T3_isotropic}
\end{figure}

Next, we consider the Ruijgrook--Wu model, in the space interval $[-0.5,0.5]$ with initial data defined as follows
\begin{align}\label{RW0}
u_0(x)=
\begin{cases}
1\qquad \textrm{ if }\, |x|\leq 1/8\\
0\qquad \textrm{ otherwise},
\end{cases}
\qquad  v_0(x) = 0,
\end{align}
where we account for reflecting boundary conditions, i.e. $v=0, \partial_x u = 0$ on the boundaries  $x=\pm0.5$.
\comment{
This test problem has been previously studied in \cite{BPR17, NP2} with IMEX Runge-Kutta methods.
We study the solution to \eqref{RW} for three different regimes of the parameters $\alpha$ and $\epsi$ and we solve the model with $N=100$ space grid points using the AP-implicit IMEX-BDF4 method. 
We report in Figure \ref{fig:T3_isotropic} the evolution of the solution for the density $u$ and the momentum $v$, respectively on the top and bottom rows. The first column represents the {rarefied regime}, for $\epsi=0.7, \alpha=1$. In this regime, the transport part dominates and the initial data propagates in the directions of the particles. This behavior is well described by the method without spurious numerical oscillations.  In the second column we have the {hyperbolic limit} for $\epsi=10^{-12}, \alpha=2/3$, corresponding to the inviscid Burger equation with a shock propagating in the right direction. Even in this case the shock profile is well captured. Finally, in the last column  we report the {parabolic limit} for $\epsi = 10^{-10}, \alpha = 4/5$, corresponding to a viscous Burger equation. As expected, the shock profile is regularized by the presence of the diffusive term.}

\subsubsection*{\revised{Test 4c. Anisotropy of the multiscale parameter $\alpha$.}}
In the last test case, we solve the model \eqref{RW} considering a multiscale parameter $\alpha$ to be a function of the space $x$, whereas the relaxation parameter is fixed to $\epsi=10^{-8}$. 
\comment{This test aims at reproducing a realistic situation where the scaling terms may depend on the physical quantities and vary in the different regions of the computational domain.  
We report the results obtained using the AP-implicit scheme with IMEX-BDF3 and the same discretization parameters of the previous test case.
The initial date is defined in the space interval $[-0.5,0.5]$ as follows
\begin{align}\label{RW0c}
u_0(x)=
\begin{cases}
1\qquad \textrm{ if }\, |x|\leq 1/8\\
0.5\qquad \textrm{ otherwise},
\end{cases}
\qquad  v_0(x) = 0.
\end{align}
}

\begin{figure}\centering
	{\includegraphics[width=4.75cm]{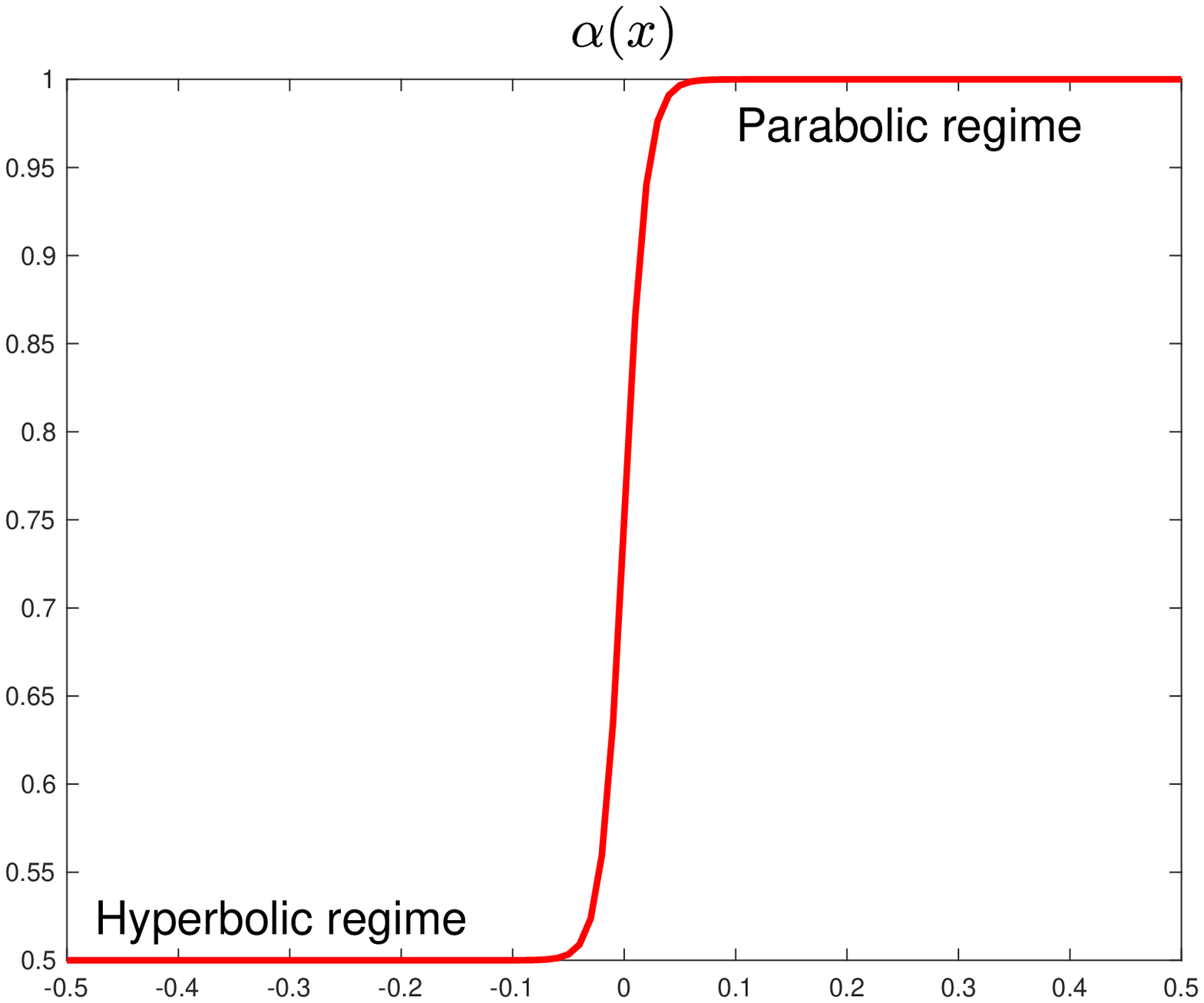}}
	\hspace{+0.25cm}
	{\includegraphics[width=4.75cm]{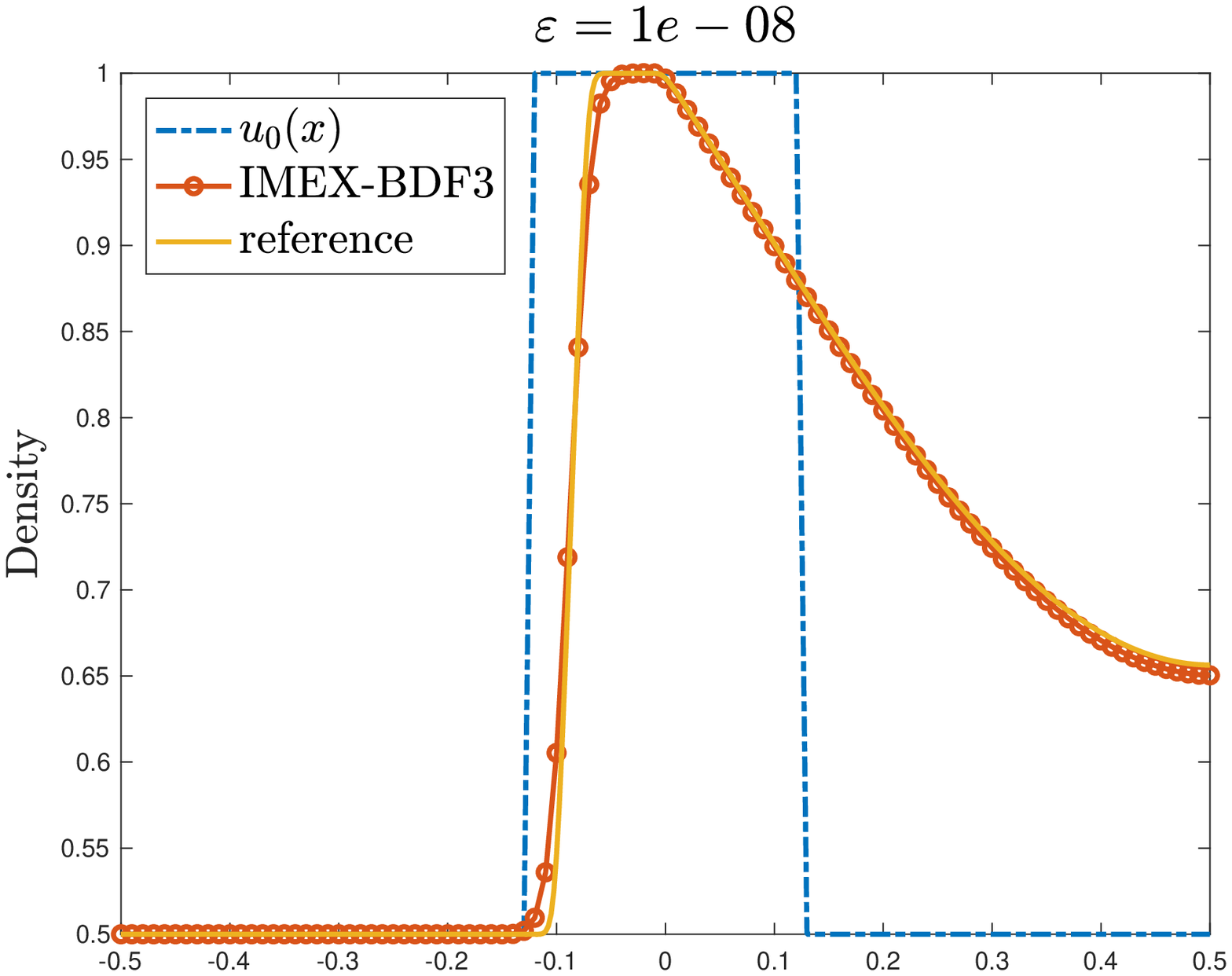}}
	\hspace{+0.25cm}
	{\includegraphics[width=4.75cm]{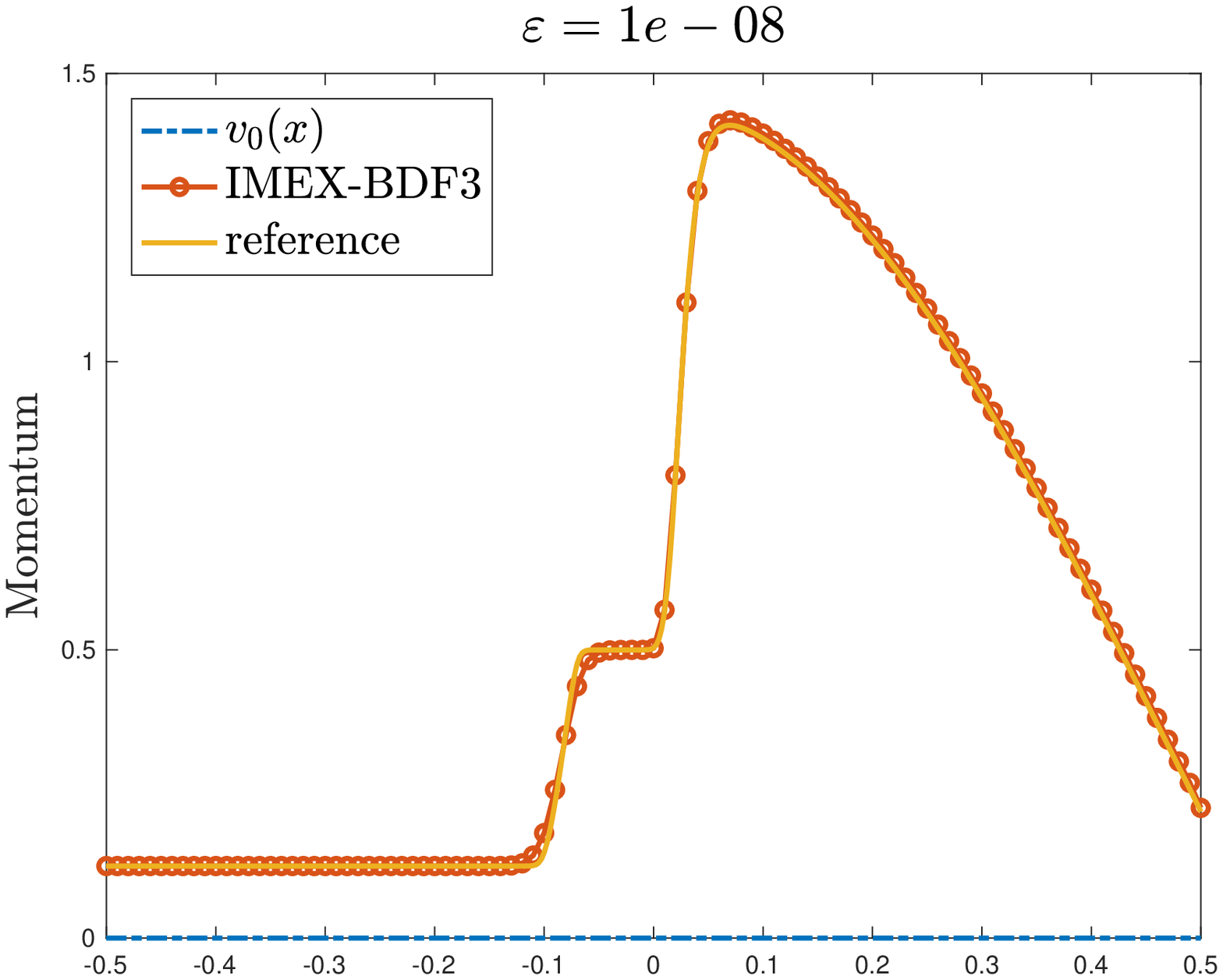}}
	\\\vspace{+0.25cm}
	{\includegraphics[width=4.75cm]{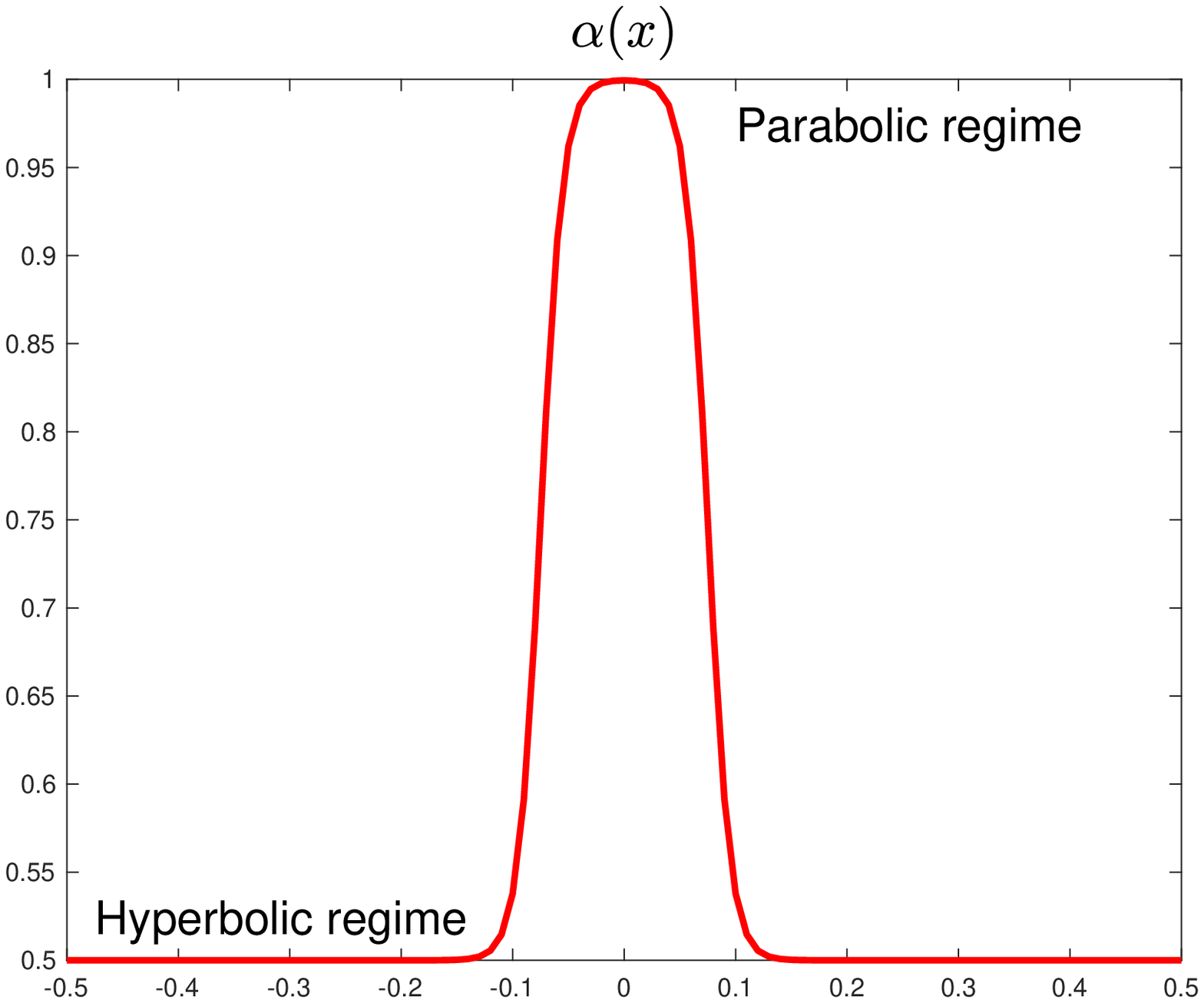}}
	\hspace{+0.25cm}
	{\includegraphics[width=4.75cm]{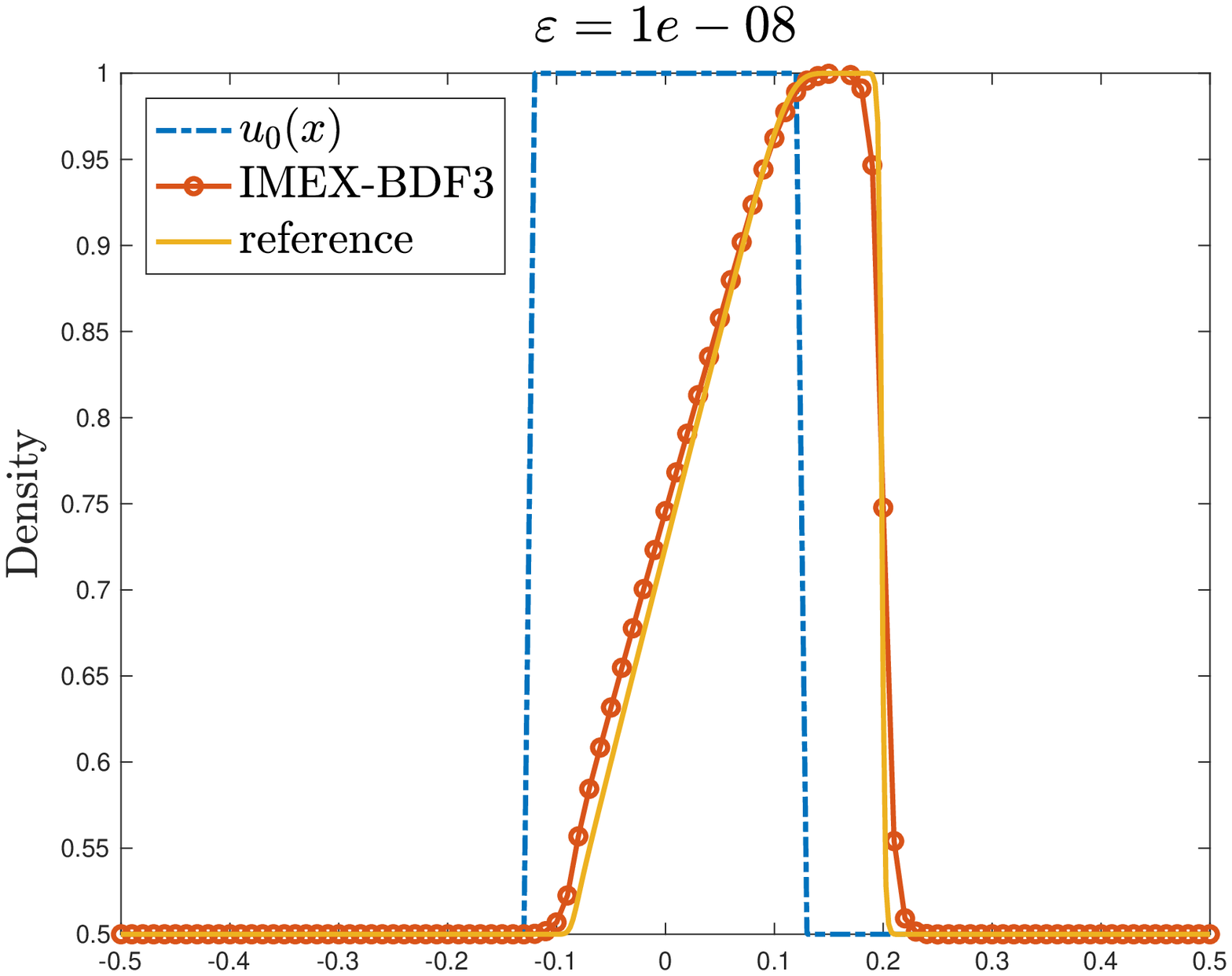}}
	\hspace{+0.25cm}
	{\includegraphics[width=4.75cm]{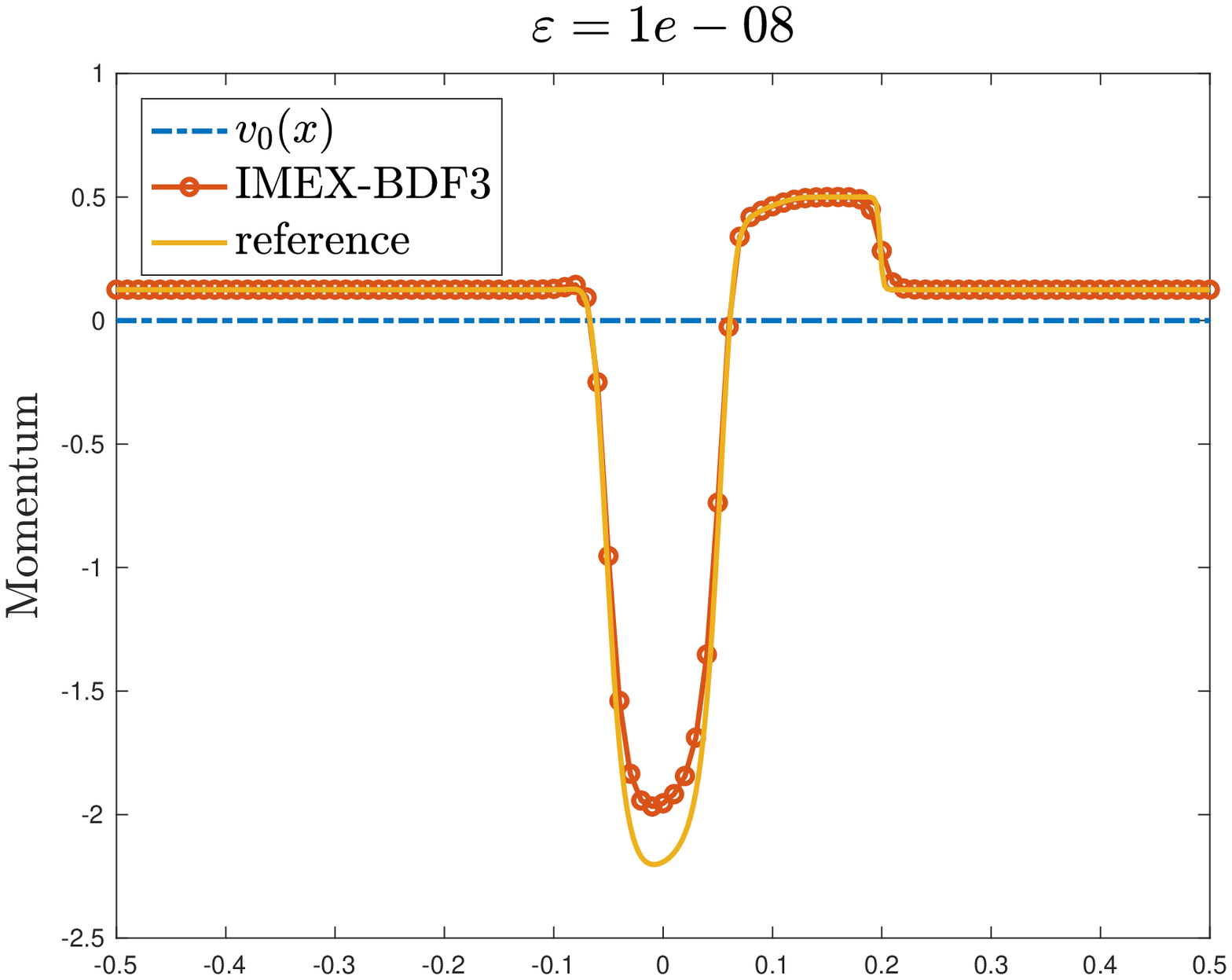}}
	
	\caption{{\em Test 4c}. \comment{Non linear Ruijgrook--Wu model with space dependent $\alpha$ and $\epsi=10^{-8}$. Left column shows the space variation of the multiscale parameter, whereas the central and right columns depict the evolution of the density $u$ and momentum $v$ from the initial data $u_0(x),v_0(x)$ to final time. Top row accounts a single variation from hyperbolic to parabolic and final time $T=0.05$, whereas second row has two transitions and final time $T=0.1$.}}\label{fig:T3_anisotropic}
\end{figure}

In the left column of Figure \ref{fig:T3_anisotropic}, we report the value of function $\alpha(x)$ as a function of the space, varying between $0.5$ (hyperbolic regime) and $1$ (parabolic regime). The central and right columns depict respectively the evolution of  $u(x,t)$ and $v(x,t)$ showing the initial and final profiles, \revised{a similar test case has been presented in \cite{BPR17} for IMEX Runge-Kutta methods.}

In the first row, we account a single variation of the regime from the hyperbolic to the parabolic 
\[
\alpha(x) = 1-\frac{1}{2}H(x), \qquad \textrm{ where } \quad H(x) = \frac{1}{1+\exp(x/\delta)}, \quad \delta = 0.01.
\]
\revised{At final time $T=0.05$ we observe  a rarefaction wave moving to the left (in the hyperbolic regime) and a smooth profile on the right (in the parabolic regime). Note that, the method captures well the complicated shock structure even at the interface between the two regions.}

Second row considers two variations of the the regime from hyperbolic to parabolic,  choosing $\alpha(x)$ as 
\[
\alpha(x) =  \frac{1}{2}-\frac{1}{2}\left(H(x+0.075)-H(x-0.075)\right),
\]
\revised{at final time is $T=0.1$ we observe that the discontinuous initial data $u_0(x)$ is blunted within the parabolic regime, whereas a shock and rarefaction waves emerge in the hyperbolic one, both well described by the numerical method.}

\section{Conclusions}
In this work we have developed a unified IMEX multistep approach for hyperbolic balance laws under different scalings. These problems, inspired by the classical hydrodynamical limits of kinetic theory \cite{CIP}, are challenging for numerical methods since the nature of the asymptotic behavior is not known a-priori and depends on the scaling parameters. Therefore, the schemes should be able to capture correctly asymptotic limits characterized by hyperbolic conservation laws as well as diffusive parabolic equations. A major difficulty in the schemes construction is represented by the unbounded growth of the characteristic speeds of the system in diffusive regimes. For these problems, we developed two different kind of approaches, originating a problem reformulation with bounded characteristic speeds, and which lead respectively to explicit (AP-explicit)  or explicit-implict (AP-implicit) time discretizations of the asymptotic limit.
Several numerical results for linear and non linear hyperbolic relaxation systems have confirmed that the IMEX multistep methods are capable to describe correctly the solution for a wide range of relaxation parameters and for different values of the scaling coefficient $\alpha$.
Compared to the IMEX Runge-Kutta approach developed in \cite{BPR17} the IMEX multistep schemes here constructed present several advantages. In particular, it is possible to achieve easily high order accuracy and in general a more uniform behavior of the error with respect to the scaling parameters is observed. In addition, when dealing with computationally challenging problems such as the case of kinetic equations with stiff collision terms its is possible to strongly reduce the number of evaluations of the most expensive part of the computation represented by the source term.
Future research will be in the direction of extending the present results to the more difficult case of diffusion limits for non linear kinetic equations and, more in general, \revised{to the case of low Mach number limits and all Mach number flows \cite{BPR, Klar, Klar2, Lem2, HJL}.}

\appendix

\section{\revised{Order conditions for IMEX-LM methods and examples}}\label{app:imex_multistep}
In this appendix we give the details of the particular IMEX-LM methods used in the manuscript. Let us recall that an order $p$,  $s$-step, IMEX-LM scheme is obtained provided that
\be
\label{eq:SPs_conditions}
\begin{aligned}
&1+ \sum_{j=0}^{s-1}a_j  = 0, \\
&1- \sum_{j=0}^{s-1}ja_j  = \sum_{j=0}^{s-1}b_j=  \sum_{j=-1}^{s-1}c_j, \\
&\frac{1}{2}+ \sum_{j=0}^{s-1}\frac{j^2}{2}a_j  = -\sum_{j=0}^{s-1}jb_j=  c_{-1}-\sum_{j=0}^{s-1}jc_j.
\\
&\qquad\vspace{+4cm}~\vdots
\\
&\frac{1}{p!}+ \sum_{j=0}^{s-1}\frac{(-j)^p}{p!}a_j  = \sum_{j=0}^{s-1}\frac{(-j)^{p-1}}{(p-1)!}b_j =  \frac{c_{-1}}{(p-1)!}+\sum_{j=0}^{s-1}\frac{(-j)^{p-1}}{(p-1)!}c_j, \\
\end{aligned}
\ee 
Moreover the following theorem holds true \cite{Ascher2}
\begin{theorem}
For an $s$-step IMEX scheme we have
\begin{enumerate}
\item If $p\leq s$ the $2p+1$ constraints of \eqref{eq:SPs_conditions} are linearly independent, therefore there exist $s$-step IMEX-LM schemes of order $s$.
\item An $s$-step IMEX-LM scheme has accuracy at most $s$.
\item The family of $s$-step IMEX-LM schemes of order $s$ has $s$ parameters. 
\end{enumerate}
\end{theorem}

Listed below the IMEX-LM methods analyzed along the paper, for further details and additional methods we refer to \cite{Ak1, Ak2,DPLMM, HR,RSSZ1, RSSZ2}.

\begin{table}[ht]
 \centering
 \caption{Examples of IMEX-LM methods. The numbers between brackets denote respectively: $s$ number of steps, $p$ order of the method. For BDF methods the number of steps is equivalent to the order of the method.} \label{tab2}
\vspace{+0.25cm}
\footnotesize
 \begin{tabular}{|p{1.1cm}|p{0.35cm}|p{6.2cm}|p{0.45cm}|p{6.cm}|  }
 \hline
IMEX & \multicolumn{2}{c|}{Explicit} & \multicolumn{2}{c|}{Implicit}\\
 \hline
 \hline
 \multirow{2}{*}{SG(3,2)} 
 & $a^T$
 &$(-\frac{3}{4},0,-\frac{1}{4},0,0)$ 
 &  $c^T$ 
 & $(0,0,\frac{1}{2},0,0)$ 
 \\
   \cline{2-5}
 & $b^T$&$(\frac{3}{2},0,0,0,0)$     &  $c_{-1}$&$1$   \\
  \hline
 \multirow{2}{*}{BDF2}
 & $a^T$ &$(-\frac{4}{3}, \frac{1}{3}, 0, 0, 0)$ 
 & $c^T$ &$(0,0,0,0,0)$   
 \\
 \cline{2-5}
 & $b^T$ &$(\frac{4}{3}, -\frac{2}{3}, 0, 0, 0)$     
 &  $c_{-1}$&$\frac2 3$\\
   \hline
    \hline
 \multirow{2}{*}{TVB(3,3)}
 & $a^T$&$(-\frac{ 3909}{ 2048}, \frac{1367}{1024}, -\frac{873}{2048}, 0, 0)$
  &  $c^T$& $(-\frac{1139}{12288}, -\frac{367}{6144}, \frac{1699}{12288}, 0, 0)$ 
 \\
 \cline{2-5}
 & $b^T$&$(\frac{18463}{12288}, -\frac{1271}{768}, \frac{8233}{12288}, 0, 0)$
 &  $c_{-1}$&$1089/2048$\\
\cline{1-5}
 \multirow{2}{*}{BDF3}
 & $a^T$&$(-\frac{18}{11}, \frac{9}{11}, -\frac{2}{11}, 0, 0)$ 
 &  $c^T$ &$(0,0,0,0,0)$  
 \\
  \cline{2-5}
 & $b^T $&$(\frac{18}{11}, -\frac{18}{11}, \frac{6}{11}, 0, 0)$     
 &  $c_{-1}$&$\frac 6 {11}$\\
   \hline
\hline
 \multirow{2}{*}{TVB(4,4)}
 & $a^T$&$(-\frac{ 21531}{8192}, \frac{22753}{8192}, -\frac{12245}{8192}, \frac{2831}{8192}, 0)$ 
 & $c^T $&$ (-\frac{3567}{8192}, -\frac{697}{24576}, \frac{4315}{24576}, -\frac{41}{384}, 0)$ 
 \\
 \cline{2-5}
 & $b^T$&$(\frac{ 13261}{8192}, -\frac{75029}{24576}, \frac{54799}{24576}, -\frac{15245}{24576}, 0)$     
 &  $c_{-1}$&$\frac{4207}{8192}$\\
\cline{1-5}
 \multirow{2}{*}{BDF4} 
 & $a^T$&$(-\frac{48}{25}, \frac{36}{25}, -\frac{16}{25}, \frac{3}{25}, 0)$ &  
 $c^T$&$(0, 0, 0, 0, 0)$ 
  \\\cline{2-5}
 & $b^T $&$(\frac{48}{25}, -\frac{72}{25}, \frac{48}{25}, -\frac{12}{25}, 0)$     
 &  $c_{-1}$&$\frac{12}{25}$   \\
  \hline
  \hline
  \multirow{2}{*}{TVB(5,5)}
 & $a^T$&$(-\frac{ 13553}{4096}, \frac{38121}{8192}, -\frac{7315}{2048}, \frac{6161}{4096}, -\frac{2269}{8192})$ 
 & $c^T$&$ (-\frac{ 4118249}{5898240}, \frac{768703}{2949120}, \frac{47849}{245760}, -\frac{725087}{2949120}, \frac{502321}{5898240})$ 
 \\
    \cline{2-5}
 & $b^T $&$(\frac{10306951}{5898240},-\frac{13656497}{2949120},\frac{1249949}{245760},-\frac{7937687}{2949120},\frac{3387361}{5898240})$     
 &  $c_{-1}$&$\frac{4007}{8192}$\\
 \cline{1-5}
 \multirow{2}{*}{BDF5}
& $a^T$&$(-\frac{ 300}{137}, \frac{300}{137}, -\frac{200}{137}, \frac{75}{137}, -\frac{12}{137})$ 
&  $c^T$&$(0, 0, 0, 0, 0)$   
 \\\cline{2-5}
 & $b^T$&$(\frac{ 300}{137},- \frac{600}{137}, \frac{600}{137}, -\frac{300}{137}, \frac{60}{137})$     
 &  $c_{-1}$&$\frac{60}{137}$\\
  \hline
 \end{tabular}
 \end{table}


\revised{\section*{Acknowledgments}
The authors are grateful to the unknown referees for having reported a series of inaccuracies in the first version of the manuscript and for the indications that led to this improved revision.}

\end{document}